\documentclass[12pt]{article}

\usepackage[]{pdfsync} 
\usepackage[francais,english]{babel} 
\usepackage[]{latexsym} 
\usepackage{amsmath} 
\usepackage{graphicx} 
\usepackage{amsfonts} 
\usepackage{amssymb}
\usepackage[T1]{fontenc}

\font \db = msbm10 at 12 pt
\font \sdb = msbm7 at 8.5 pt 
\def\iff{if and only if }

  \def \o{\omega}      
\def \e{\varepsilon}  \def \f{\varphi}    \def \rr{\varrho}

\font \db = msbm10 at 12 pt
\font \sdb = msbm7 at 8.5 pt       
\def \R{\mbox{\db R}}  \def \H{\mbox{\db H}} \def \E{\mbox{\db E}} \def \P{\mbox{\db P}} 
\def \N{\mbox{\db N}}      
\def \sR{\mbox{\sdb R}}    
\def \sN{\mbox{\sdb N}}

\def \Tr{{\rm tr}}   

\newcommand \parsn{\par \smallskip \noindent }  \newcommand \parmn{\par \medskip \noindent } 
\newcommand \pars{\par \smallskip }  \newcommand \parm{\par \medskip } 
\newcommand \parb{\par \bigskip }  

\def \S{\mbox{\db S}}   \def \C{\mbox{\db C}}   \def \Z{\mbox{\db Z}}  
      \def \sZ{\mbox{\sdb Z}} 

           \def \EE{{\cal E}}
\def \FF{{\cal F}}        \def \II{{\cal I}}   
\def \LL{{\cal L}}     \def \NN{{\cal N}}    \def \O{{\cal O}}    
  \def \RR{{\cal R}}         
      \def \WW{{\cal W}}

     \def \ub{\underbar}        \def\lmt{\longmapsto}
  \def \Lra{\Longrightarrow}   \def \LRa{\Leftrightarrow} 

                      \def \sea{\searrow}

   \def \ds{\displaystyle}   \def \ts{\textstyle}       \def \ss{\scriptstyle}
\def \rt1{\sqrt{-1}\,\,}  \def \1{^{-1}}            \def \2{^{-2}}             \def \5{{\ts \frac{1}{2}}}
\def\moins{\rotatebox{30}{\raise0.3pt\hbox{$\,\scriptstyle\setminus$}\,}} 

\def\cotg{{\rm cotg}\,}   \def\tg{{\rm tg}\,}   \def\coth{{\rm coth}\,}        \def\th{{\rm th}\,}
\def\ch{{\rm ch}\,}       \def\sh{{\rm sh}\,}   \def\arctg{{\rm arctg}\,}   

\def \parn{\par\noindent }  \def\indf{\leavevmode\indent }  
\def \parsn{\par\smallskip \noindent } \def \parmn{\par\medskip \noindent } 
\def \parm{\par \medskip } \def \pars{\par \smallskip } 

\begin{document}

\newtheorem{defi}[subsubsection]{Definition} 
\newtheorem{theo}[subsubsection]{Theorem} 
\newtheorem{prop}[subsubsection]{Proposition}  
\newtheorem{cor}[subsubsection]{Corollary} 
\newtheorem{lem}[subsubsection]{Lemma} 
\newtheorem{rem}[subsubsection]{Remark} 
\newtheorem{subs}[subsubsection]{} 

\newtheorem{Cor}[subsection]{Corollary} 
\newtheorem{Lem}[subsection]{Lemma} 
\newtheorem{Theo}[subsection]{Theorem} 
\newtheorem{Prop}[subsection]{Proposition} 
\newtheorem{Rem}[subsection]{Remark} 

\newcommand \beq{\begin{equation}} \newcommand \eeq{\end{equation}} 

\newcommand \bthe{\begin{theo}}  \newcommand \ethe{\end{theo}}   
\newcommand \bpro{\begin{prop}}  \newcommand \epro{\end{prop}}   
\newcommand \bcor{\begin{cor}}    \newcommand \ecor{\end{cor}}     
\newcommand \blem{\begin{lem}}   \newcommand \elem{\end{lem}}   
\newcommand \brem{\begin{rem}}  \newcommand \erem{\end{rem}} 
\newcommand \bdefi{\begin{defi}}   \newcommand \edefi{\end{defi}}  
\newcommand \bsub{\begin{subs} \  }   \newcommand \esub{\end{subs}}  

\newcommand \bCor{\begin{Cor}}    \newcommand \eCor{\end{Cor}}     
\newcommand \bLem{\begin{Lem}}   \newcommand \eLem{\end{Lem}}   
\newcommand \bThe{\begin{Theo}}  \newcommand \eThe{\end{Theo}}   
\newcommand \bPro{\begin{Prop}}  \newcommand \ePro{\end{Prop}}   
\newcommand \bRem{\begin{Rem}}  \newcommand \eRem{\end{Rem}}

\title{\bf Exact Small Time Equivalent for the density of the Circular Langevin Diffusion}
\author{{\Large J. Franchi }\\ {\normalsize  IRMA, Universit\'e de Strasbourg et CNRS}, \\ {\normalsize 7 rue Ren\'e Descartes, 67084 Strasbourg, France}}
\date {Octobre 2015} 
\maketitle

\begin{abstract}
A small time equivalent of the density is obtained for the circular analogue of the Langevin diffusion,
which is strictly hypoelliptic (and non-Gaussian), hence of a different nature as the known sub-Riemannian case. The singular case, analogous to the case of conjugate points (the cut-locus problem) in the sub-Riemannian framework, is totally handled too, though much more difficultly.  
\end{abstract} 

\parm  
\centerline{ \small {\bf  CONTENT}}
 \vskip 1mm  \noindent 
\ref{introd}. \  Introduction \hfill page \pageref{introd}   \parn 
\ref{s.Dudley1}. \  The circular Langevin diffusion and the main result \hfill page \pageref{s.Dudley1}  \parn
\ref{methFPGpinned}. \  Pinned Fourier transform \hfill page \pageref{methFPGpinned} \par  \quad  \ref{sec.Bbd}) \  Brownian bridge disintegration  \hfill page \pageref{sec.Bbd} \par 
\indent \quad  \ref{sec.FPGF}) \  Fourier-Plancherel Formula \hfill page \pageref{sec.FPGF} \parsn 
\ref{sec.changbsuit1scomp}. \  Analysis and domination of $\,P_\e(\xi',\xi)$  \hfill page \pageref{sec.changbsuit1scomp}  \par 
\indent \quad  \ref{handlerrort})  \   Analysis of the phase vector $\;{\ds U^\e}$   \hfill page \pageref{handlerrort} \par \quad \ref{sec.IPPT2}) \   Applying Theorem 2.1 of [T2]   \hfill page \pageref{sec.IPPT2} \par 
\indent \quad \ref{sec.changbsuit1sapple[T2]}) \  Domination of the decreasing integral   \hfill page \pageref{sec.changbsuit1sapple[T2]} \parn 
\ref{sec.asymptpe}. \  Small-time equivalent for $\,p_\e\,$, provided $\,w\not=0$  \hfill page  \pageref{sec.asymptpe}  \parn 
\ref{sec.casew=0}. \   The singular case $\,w=0$  \hfill page  \pageref{sec.casew=0} \par 
\indent \quad  \ref{sec.Analw=0}) \  Analysis of the phase vector $\,\tilde U^\e$ and second use of [T2]   \hfill page  \pageref{sec.Analw=0} \par 
\indent \quad  \ref{sec.Analogw=0}) \  Domination of the decreasing integral in the singular case  \hfill page  \pageref{sec.Analogw=0} \par 
\indent \quad  \ref{sec.KApprox0}) \  Intermediate small-time equivalent in the singular case  \hfill page  \pageref{sec.KApprox0} \parn 
\ref{sec.calculFonctQuad}.  \  Computation of a quadratic Laplace transform \hfill page  \pageref{sec.calculFonctQuad}  \parn 
\ref{sec.CompLoi(om,om^2)}. \  Evaluation of the oscillatory integral $\, {\ds \II_\e(y,z) }$  \hfill page  \pageref{sec.CompLoi(om,om^2)} \par
\indent \quad  \ref{sec.UsCalcQuadr}) \  Reduction to a finite-dimensional oscillatory integral \hfill  page \pageref{sec.UsCalcQuadr} \par
\indent \quad  \ref{sec.z=0>y1}) \  First sub-case\,: $\, y\le0\,,\, z=0$  \hfill  page \pageref{sec.z=0>y1} \par
\indent \quad  \ref{sec.znot=0<y}) \  Second sub-case\,: $\, y>0\,,\, z=0$   \hfill  page \pageref{sec.znot=0<y}\par
\indent \quad  \ref{sec.znot=0<yz}) \  Third sub-case\,: $\, z\not=0$   \hfill  page \pageref{sec.znot=0<yz} 
\parn 
\ref{Ref}. \  References \hfill page  \pageref{Ref}  \par 

\section{Introduction} \label{introd} \indf 
   The problem of estimating the heat kernel, or the density of  a diffusion, particularly as time goes to zero, has been extensively studied for a long time, firstly in the elliptic case, and then largely solved and understood in the sub-Riemannian case too. Let us mention only the articles [V], [A], [BA1], [BA2], [L], 
and the existence of other works on that subject by Azencott, Molchanov and Bismut, quoted in  [BA1]. \par

       To summary roughly, a very classical question addresses the asymptotic behavior (as $\,s\sea 0$) of the density $\,p_s(x,y)$ of the diffusion $(x_s)$ solving a Stratonovich stochastic differential equation 
\begin{equation*}
x_s = x + \sum_{j=1}^k \int_0^s V_j(x_\tau)\circ dW^j_\tau + \int_0^s V_0(x_\tau)\, d\tau \, , 
\end{equation*}
where the smooth vector fields $\,V_j\,$ are supposed to satisfy a H\"ormander condition. \par
   The elliptic case being very well understood for a long time ([V], [A]), the studies focussed then on the sub-elliptic case, that is to say, when the strong H\"ormander condition (that the Lie algebra generated by the fields $\,V_1,\ldots, V_k\,$ has maximal rank everywhere) is fulfilled. In that case these fields generate a sub-Riemannian distance $\,d(x,y)$, defined as in control theory, by considering only $C^1$ paths 
whose tangent vectors are spanned by them. Then the wanted asymptotic expansion tends to have the following Gaussian-like form\,: 
\begin{equation} \label{f.GaussLAsymp} 
p_s(x,y) = s^{-d/2}\, \exp\!\big( - d(x,y)^2/(2 s)\big) \bigg( \sum_{\ell=0}^n \gamma_\ell(x,y)\, s^\ell + \O(s^{n+1})\bigg)  
\end{equation}
for any $\,n\in\N^*$, with smooth $\,\gamma_\ell$'s and $\gamma_0>0\,$, provided $\,x,y$ are not conjugate points (and uniformly within any compact set which does not intersect the cut-locus). See in particular ([BA1], th\'eor\`eme 3.1). Note that the condition of remaining outside the cut-locus is here necessary, as showed in particular by [BA2]. \par
   The methods used to get this or a similar result have been of different nature. In [BA1], G. Ben Arous proceeds by expanding the flow associated to the diffusion (in this direction, see also [Ca]
) and using a Laplace method applied to the Fourier transform of $\,x_s\,$, then inverted by means of Malliavin's calculus (with a deterministic Malliavin matrix). \par 
   The strictly hypoelliptic case, i.e., when only the weak H\"ormander condition (requiring the use of the drift vector field $V_0$ to recover the full tangent space) is fulfilled, remains much more problematic, and then rarely addressed. There is a priori no longer any reason that in such case the asymptotic behavior of $\,p_s(x,y)$ remains of the Gaussian-like type (\ref{f.GaussLAsymp}), all the less as a natural candidate for replacing the sub-Riemannian distance $\,d(x,y)$ is missing. Indeed this already fails for the mere (Gaussian) Langevin process $\big(\omega_s\,,\int_0^s\omega_\tau\, d\tau\big)$: the missing distance is replaced by a time-dependent distance $\,d_s(x,y)$ which presents some degeneracy in one direction, namely $\,d(x,y)^2/(2 s)$ is replaced by 
$$\frac{\hbox{\footnotesize 6}}{s^3}\,\Big|(x-y) - \frac{s}{2}\, (\dot x-\dot y)\Big|^2 + \frac{\hbox{\footnotesize 1}}{2s}\,\big|\dot x-\dot y\big|^2 = \frac{\hbox{\footnotesize 1}}{2s}\Big(\big|\dot x-\dot y\big|^2 +\frac{\hbox{\footnotesize 12}}{s^2}\,\Big|(x-y) - \frac{s}{2}\, (\dot x-\dot y)\Big|^2\, \Big) . $$ 
   See also [DM] for a more involved (non-curved, strictly hypoelliptic, perturbed) case where Langevin-like estimates hold (without precise asymptotics), roughly having the following Li-Yau-like form\,:
\begin{equation} \label{F.DM} 
C\1\, s^{-N}\,e^{-C\, d_s(x_s,y)^2} \,\le \, p_s(x,y)\, \le \, C\, s^{-N}\,e^{-C\1 \, d_s(x_s,y)^2} \,, \quad \hbox{for }\; 0<s<s_0\,. 
\end{equation} 

    In [F] was computed the small time asymptotics of a toy model, namely a diffusion in the second Wiener chaos, simplest case after the Gaussian setting, in the specific off-diagonal regime of a dominant normalized Gaussian contribution. The exponential term appeared as given by the same time-dependent distance as in the Langevin case, the strictly second chaos coordinate appearing only in the off-exponent term, as a perturbative contribution.  \par
       A stronger interest lies on a significant strictly hypoelliptic diffusion, namely the relativistic diffusion, first constructed in Minkowski's space (see [Du], [F-LJ2]). It makes sense on a generic smooth Lorentzian manifold as well, see [F-LJ1]. In the simplest case of Minkowski's space $\R^{1,d}$, it consists in the pair $(\xi_s, \dot\xi_s)\in \R^{1,d}\times\H^d$ (parametrized by its proper time $s$, and analogous to a Langevin process), where the velocity $(\dot\xi_s)$ is a hyperbolic Brownian motion. Note that even there, a curvature constraint must be taken into account, namely that of the mass shell $\,\H^d$, at the heart of this framework. \par 
 
       This (Dudley) relativistic diffusion, even restricted to 3 dimensions which already contain the essence of the difficulty, constitutes a significant example, altogether explicit, physical and not too much  complicated, allowing a priori to progress towards the understanding of a more generic, but less accessible to begin with, strictly hypoelliptic case. 
However even this apparently simple example is not easy at all to analyze, regarding the small time asymptotics.  \par  
     The present work handles a near example, but in which the curvature is non-negative and the fiber is compact, namely $\R^2\times \S^1\equiv T^1_+\R^2$, endowed with its natural strictly hypoelliptic diffusion $(\xi_s, \dot\xi_s)$, analogue of both the Langevin and the Dudley diffusions, we name ``circular Langevin diffusion''. In this setting we obtain the exact small time ($\e\sea 0$) equivalent for the density (heat kernel) $p_\e(x_0;x)$. To the author's knowledge, this is the first example of a complete result of this type in a strictly hypoelliptic framework. It reveals a different nature from the known sub-elliptic framework, see Section \ref{s.Dudley1} below. \par 
       

       The present work was strongly influenced by the beautiful article [BA1], which decisively handled the off-cut locus generic sub-Riemannian framework, as far as the small-time asymptotics of the heat kernel is considered. Thus the strategy adopted below roughly resembles the strategy followed by G. Ben Arous, at least for the non-degenerate case. However the present purpose is to deal with a strictly hypoelliptic situation, to which [BA1] does not apply, and to the best of the author's knowledge, remained unsolved. A main obstacle to handle a strictly hypoelliptic framework is the lack of sub-Riemannian distance. For that reason, the strategy used here only partially follows the method of [BA1], and stresses on the Brownian bridge instead of the Brownian motion.       
Another sensitive though more technical reason is that the Malliavin matrix which will intervene below is no longer deterministic as in [BA1], so that we shall resort to the treatment of infinite-dimensional oscillatory integrals (as in [T2], up to a modification).   \par
        As the author is not yet ready to handle a generic strictly hypoelliptic framework, and actually doubts that a generic unified treatment (or even, unified generic result) be possible (different strictly hypoelliptic frameworks could produce different types of results\,; the present one already differs notably from the classical Gaussian Langevin case, which even ignores degenerate cases), the focus is here put on a simple first example, which allows to concentrate on the heart of a method to get beyond the sub-Riemannian framework. In this way some difficulties handled in [BA1] don't have their counterpart here, in particular the asymptotic development is here limited to the first order term (i.e., only provides an equivalent as time goes to 0). Another simplification is due to the boundedness of the drift vector field $V_0$, which does not hold in the relativistic (Dudley) case. It is not clear whether the latter could be handled similarly as is done here. On the contrary, the choice of considering a 3-dimensional case is unessential, but avoids even heavier notation and computations which higher dimensions would cause. Otherwise, focussing on the present (relatively) simple example allows to also handle the case of conjugate points\,; which is delicate, even in the sub-Riemannian framework\,: the cut locus constitutes a real difficulty in that matter, see for example [BA2], and its case does not seem to be generically solved in a sub-elliptic framework. Furthermore the choice of a relatively simple particular framework allows to fully explicit all coefficients of the wanted equivalents, which will likely be out of reach in an even slightly more generic framework. \pars

       The content is organized as follows. \parn
In Section \ref{s.Dudley1} are described the strictly hypoelliptic diffusion under consideration and the central result\,: Theorem \ref{th.mainres}.   \parn   
Section \ref{methFPGpinned} develops the leading strategy of the proof, which resembles that of [BA1], with some difference\,: whereas two main tools are as in [BA1], namely a Fourier-Parseval expression for the density of the heat kernel under consideration and the Girsanov transform, the Brownian component is here pinned, and the non-existing geodesic tube of [BA1] is here replaced by some partially pinned ``pseudo-geodesic'', which happens here to be directed by a mere segment. 
This leads to a somewhat simplified  expression for the density kernel, which however contains an oscillatory integral  which is not directly computable. \parn
In Section \ref{sec.changbsuit1scomp} is analyzed the infinite-dimensional oscillatory integral which results from the preceding. 
As in [BA1] it is necessary to resort to Malliavin calculus, but in a more involved way since the Malliavin matrix is  here no longer deterministic. A key argument will rely on the analysis of the decay of such infinite-dimensional oscillatory integrals, as performed in particular in [T2].  This requires technical estimates. Another key argument is an estimation of a variance from below. As a consequence, a computable expression of the density is deduced in Section \ref{sec.asymptpe}, in the non-degenerate case $\,w\not=0\,$. \parn
Section \ref{sec.casew=0} is devoted to the singular case $w=0\,$, which is much more delicate, and is analogous to the study at a sort of cut-locus, relating to some absent or hidden metric. In this section is performed a reduction to a simpler oscillatory integral, in a way resembling the non-degenerate case, but with a more delicate estimation of the variance. \parn 
In Section \ref{sec.calculFonctQuad} the Fourier-Laplace transform of somme quadratic Brownian-bridge functional is computed, which is needed to solve the singular case $w=0\,$. \parn
Finally in the last section \ref{sec.CompLoi(om,om^2)} is analyzed the delicate finite-dimensional oscillatory integral, which ultimately yields the wanted asymptotics of the singular case $w=0\,$; which can have two sensibly different natures, depending on the target-point. 

\section{The circular Langevin diffusion and the main result} \label{s.Dudley1} 
\begin{figure}
\centering
\includegraphics[scale=1.9]{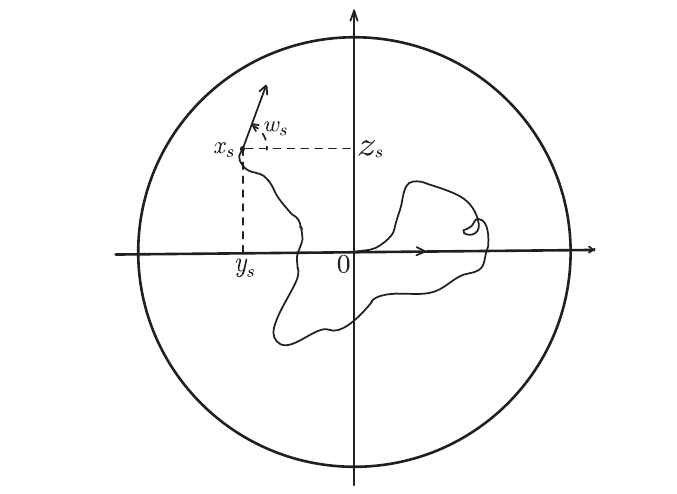}
\caption {A circular Langevin trajectory from $0\equiv (0,0,0)$ to $x_s=(w_s,y_s,z_s)$}
\label{CircLangev.pdf}
\end{figure}

\indf 
   The circular Langevin diffusion reads \  $x_s =(\o_s, y_s, z_s)$, \   with a standard real Brownian motion $(\o_s)$  and \pars
 \centerline{  ${\ds y_s:= \int_0^s\cos\o_\tau\, d\tau\, , \quad  z_s:= \int_0^s\sin\o_\tau\, d\tau\,}$. }\parsn  
We have here parametrized $\,T_+^1\R^2\equiv \S^{1}\times\R^2$ by $(\omega, y, z)\in\R^3$
A typical trajectory is depicted in Figure \ref{CircLangev.pdf}. \vspace{-1mm} 

  $\!\!$  Consider the scaled diffusion $\,{\ds s\mapsto x^\e_s = \Big( \sqrt{\e}\,\o_s\,,\, \e\! \int_0^{s}\! \cos\!\big[\sqrt{\e}\,\o_{\tau}\big] d\tau \,,\, \e\! \int_0^{s}\! \sin\!\big[\sqrt{\e}\,\o_{\tau}\big] d\tau\Big)}$, \parn which has the same law as $\,{\ds [s\mapsto x^1_{\e s}= x_{\e s}]}$ and satisfies the stochastic differential equation \parm 
\centerline{$d x^\e_s = \sqrt{\e}\,V_1(x^\e_s)\, d \omega_s + \e\, V_0(x^\e_s)\, ds\;$ (analogous to (2.1) in [BA1]),} \parsn 
with \quad  ${\ds V_1= \partial_\o\,,\, V_0= \cos\o\,\partial_y+ \sin\o\,\partial_z = \big[V_1,[V_1,V_0]\big] ,\;  [V_0,V_1] = \sin\o\,\partial_y - \cos\o\,\partial_z\,}$, \parn
which span $T\R^3$ at any point. Note that the vector fields $V_j$ and their derivatives are here bounded. The H\"ormander hypoellipticity criterion ensures the existence of a smooth density $\,p_\e(\cdot , \cdot)$ with respect to he Lebesgue measure for the relativistic diffusion $\,x_\e\,$. The density of $\,x^\e_1\,$ is $\,p_\e\,$ as well. We are interested in  small times $\e$. \   According to Remark \ref{rem.ActGroup} below, by homogeneity we can restrict to the starting point 0. \par      
   We thus fix $\,(w,y,z)\in\R^3$, and look for the equivalent as time $\,\e\sea 0\,$ of the generic value $\,p_\e\big(0\, ; (w,y ,z)\big)$ of the density of $\,x^\e_1\,$.   The idea is that $\,\sqrt{\e}\,\o_\cdot\,$ should concentrate around the linear $\,s\mapsto ws\,$, at least in the non-degenerate case $\,w\not=0\,$ which we shall consider first. This non-degeneracy condition happens to be necessary in particular in the proof of  Proposition \ref{pro.estintdecre} below, and, seen from the origin $(0,0,0)$, seems to correspond to an off-cut locus regime, whereas cut locus problems are significant, see [BA2]. This lets understand why the singular case $\,w =0\,$ is sensibly more delicate, as appears in our result, Theorem \ref{th.mainres} below. \pars

   The result to which this article is devoted is the following. 
\bThe \label{th.mainres} \  Consider a standard real Brownian motion $(\o_s)$, the circular Langevin diffusion \  
${\ds  x_s =\left(\o_s\,, \int_0^s\cos\o_\tau\, d\tau\, ,     \int_0^s\sin\o_\tau\, d\tau\right)_{s\ge 0}}$, and denote by $\,p_s\big(0\, ; \cdot \big) $ the density of the law of $\,x_s\,$ started from $0\equiv (0,0,0)$. It admits the following equivalents as $\,\e\sea 0$. \pars 
$(i)$ Non-degenerate case $\,w\not=0\,$. \  For any $(w,y,z)\in\R^*\times\R^2$ we have\,:  
$$ p_\e\big(0\, ; (w, y , z)\big) =\, \frac{\big(2+o(1)\big)\,w^2}{\pi\,\e^{3}\,\sqrt{2\pi\,\e \, \Delta(w)}} \times \exp \!\left[ -\,\frac{w^2}{\e} \times\! \left(\frac{1}{2}+\frac{\psi(w,y/\e,z/\e) }{\Delta(w) }\right) \right]\! ,  \vspace{-2mm} $$ 
with, to express the energy minimizing action functional,  the following positive determinant\,: \vspace{-1mm} 
$$ \Delta(w) := 1- \frac{\sin^2\! w + 4(1-\cos w)}{w^2} + \frac{4 (1-\cos w)\sin w }{w^3} \, \raise1.9pt\hbox{,}  $$ and the following positive (i.e., non-negative non-degenerate) quadratic form $ \psi(w,\cdot,\cdot)$:   \vspace{-1mm}
{
$$ \psi(w,y,z) =  {\ts \left(1+ \frac{\sin(2w)}{2w} - 2\big(\frac{\sin w}{w }\big)^2\right)\! \left(\frac{\sin w}{w} - y\right)^2 +  \left(1- \frac{\sin(2w)}{2w} - 2 \big(\frac{1-\cos w}{w}\big)^2\right)\! \left(\frac{1- \cos w}{w} - z\right)^2 } \vspace{-1mm} $$
$$ {\ts -\,4 \left(\frac{\sin^2\! w}{2w} - \frac{(1-\cos w)\sin w}{w^2}\right)\! \left(\frac{\sin w}{w} - y\right)\!\left(\frac{1- \cos w}{w} - z\right) }  \raise0.4pt\hbox{.} $$ }
   $(ii)$ First degenerate case. \  For any $\,y\le 0$, we uniformly have\,:  
$$ p_\e\big(0\, ; (0,y, 0)\big) =\, \exp\!\left[-4\pi^2\,\frac{\e -y}{\e^2}\right]\times\frac{2\sqrt{2}\,e+o(1)}{\e^{3}\,\sqrt{\e-y}} \times \int_0^\infty\, \frac{\sin\!\left(\frac{\pi}{8} + x - \frac{1}{2}\,\arctg  x\right)}{(x^2+1)^{1/4}}\, dx\, . \vspace{-0mm} $$ 
$(iii)$ Second degenerate case. \  For any $\,(y,z)\in \R^2$ such that $\,z\not=0\,$ or $\,y>0$, we have\,:  
$$  p_\e\big(0\, ; (0,y, z)\big) =\, \frac{\exp\!\left[-\pi^2\,C_\e(y,z)\right]}{\e^4\, C_\e(y,z)^{3/4}} \times \frac{(2\pi/\sh\pi)^{1/4}}{ \sqrt{\pi -2\,\th\!{\frac{\pi}{2}}}} \int_0^\infty\, \frac{\sin\!\left(\frac{3\pi}{16} + x\right)}{x^{1/4}}\, dx\, \, \big(1+o(1)\big) , \vspace{-2mm} $$ 
with \   $\,{\ds C_\e(y,z) := \,\frac{\pi\,z^2}{2\left(\pi -2\,\th\!{\frac{\pi}{2}}\right)\e^3} + \frac{y-\e}{\e^2}\,}$\raise1.9pt\hbox{.} 
\eThe
\bRem \label{rem.ActGroup}  \  {\rm  The statement in Theorem \ref{th.mainres} is written when starting from 0. Actually this covers the generic initial value as well, by homogeneity, owing to the obvious action of the underlying affine rotation group. Namely, we merely have\,:  \vspace{-1mm} 
$$ p_\e\big((w_0 , y_0 , z_0) ; (w , y , z)\big) 
$$ 
$$ =\, p_\e\big((0 , 0 , 0) ; (w-w_0\, , (y-y_0)\cos w_0 + (z-z_0)\sin w_0 \, , (z-z_0)\cos w_0 - (y-y_0)\sin w_0 )\big) . $$ 
}\eRem 
\bRem \label{rem.coeffwnot=0}  \  {\rm     The coefficients $\Delta(w)$ and $\psi (w,y,z)$ that appear in the non-degenerate case of the statement come essentially from some inverse Malliavin matrix  $(DU^0)\1$ ; \  in particular (see Remark \ref{rem.detDV0} and Section \ref{sec.asymptpe}  below), \  we have \quad  $\Delta(w) = 4w^4 \det(DU^0) $ \big($= \frac{w^6}{2160} + \O(w^8)$ for small $|w|$, which lets appear the singularity of the case $w=0$\big) \  and \  \parsn
\centerline{$\psi (w,y,z)= 2w^2 \left(y-\frac{\sin w}{w}\,\raise1.3pt\hbox{,}\,\frac{1- \cos w}{w} - z\right)\times DU^0 \times^t\!\!\left(y-\frac{\sin w}{w}\,\raise1.3pt\hbox{,}\,\frac{1- \cos w}{w} - z\right)$,} \parsn  
which vanishes if and only if \  $\frac{\sin w}{w} - y = \frac{1- \cos w}{w} - z=0$. \parn  
In polar coordinates \  $y=\rr\cos\alpha\,,\, z= \rr\sin\alpha\,$, this reads 
{\small 
$$ {\ts \psi (w,y,z)= \rr^2 \left[1-\frac{\sin w}{w}\,\cos(w-2\alpha) - \frac{2(1-\cos w)}{w^2} \big(1-\cos(w-2\alpha)\big) \right] \hskip20mm \vspace{-2mm} } $$
$$ {\ts \hskip20mm -\,2\rr \left[1-\frac{\sin w}{w}\right] \frac{\sin(w-\alpha) + \sin\alpha}{w} + 2 \left[1-\frac{\sin w}{w}\right] \frac{1-\cos w}{w^2} }\,\raise1.7pt\hbox{.} \vspace{-1mm} $$ } 
}\eRem 
\bRem \label{rem.supportborné}  \  {\rm  The expressions given by Theorem \ref{th.mainres} are only asymptotically valid.  Indeed, starting from 0, by the very definition of $(x_\e)$ for any positive time $\e\,$ we have $\,y_\e^2+z^2_\e \le \e^2\,$, so that $\,p_\e\big(0\, ; (w,y,z)\big)$ must vanish for $\,{y^2+z^2}\ge \e^2$. This is coherent with Theorem \ref{th.mainres} \big(as  $\lim_{\e\to 0}\limits\, p_\e(0\,; \cdot)= 0$\big), but forbids to think of these asymptotics as more widely valid. 
}\eRem 
\bRem \label{rem.commentexposants}  \  {\rm  The initial null speed of the third component $\,z_s\,$ of $\,x_s\,$ explains the stronger energy in the exponent $(iii)$ when $\,z\not=0$. Similarly, the initial positive speed of the second component $\,y_s\,$ makes more difficult that it rapidly reaches a non-positive value $y$, than a positive one. This explains the above difference between the exponents of $\,p_s\,$, regarding both cases $(ii)$ and \big[$(iii)$, $z=0$\big]. Finally the presence of $(y-\e)$ in the exponents is natural, owing to the deterministic contribution at time $\,\e$ (corresponding to $\,w\,$ remaining at 0).  
}\eRem 
\bRem \label{rem.comment-actions-cutlocus}  \  {\rm  Though there is no natural underlying metric in this strictly hypoelliptic framework, the exponents given by Theorem \ref{th.mainres} describe some action-energy  functionals, and  furthermore the singular case $\,w=0\,$ lets think of a cut-locus which would be relative to some absent or hidden metric.  \par Moreover, these exponents show important differences with respect to the sub-elliptic case, generically handled in [BA1] for the non-degenerate case, and then partially \big(precisely, on the diagonal, which corresponds to $w=y=z=0$ here\big) in [BA2] for the degenerate case. Indeed, on the one hand the energy functional in Theorem \ref{th.mainres}$(i)$ above cannot be expressed as $d^2(0,x)/2\e$ (as in the sub-elliptic case) and even hardly as $d_\e^2(0,x)/2\e$ (as in the flat Langevin case), and on the other hand the decay factor is polynomial in the diagonal sub-elliptic case, while exponential in Theorem \ref{th.mainres}$(ii)(y=0)$ above. 

}\eRem   \vspace{-3mm} 

\section{Pinned Fourier transform} \label{methFPGpinned} \indf 
   Recall that  $\,p_\e\,$ denotes the density of  
${\ds \Big( \sqrt{\e}\,\o_1 ,\, \e\! \int_0^{1} \cos\!\big(\sqrt{\e}\,\o_{\tau}\big)\, d\tau \,,\,  \e\! \int_0^{1} \sin\!\big(\sqrt{\e}\,\o_{\tau}\big)\, d\tau\Big)}$.  \vspace{-2mm} 
\subsection{Brownian bridge disintegration} \label{sec.Bbd}  \indf 
   Let us condition on $\,\o_1=\frac{w}{\sqrt{\e}}\,$\raise1.7pt\hbox{,} and use the Brownian bridge $(\o_s\,, 0\le s\le 1)$ from 0 to $\frac{w}{\sqrt{\e}}\,$\raise1.5pt\hbox{,}
whose law will be denoted by $\,\P_0^{w/\sqrt{\e}}$, to first disintegrate the Brownian law. 
\blem \label{lem.Bbdes} \  For each real $\,w$,  \  {${\ds \sqrt{2\pi\,\e}\,e^{w^2/2\e}\, p_\e\big(0;(w,\cdot,\cdot)\big)}\,$ is the density of \parn  
${\ds \Big( \e\! \int_0^{1} \cos\!\big(w\tau\!+\!\sqrt{\e}\,\o_{\tau}\big)\, d\tau \,,\, \e\! \int_0^{1} \sin\!\big(w\tau\!+\!\sqrt{\e}\,\o_{\tau}\big)\, d\tau\Big)}$} under the Brownian bridge law $\,\P_0^{0}\,$. 
\elem 
\ub{Proof} \quad  For all test functions $\,\f,\Phi$, denoting by $\P_0$ the standard Brownian law we have\,:  
$$ \int_{\sR} \left[\int_{\sR^2} \Phi(u,v)\, p_\e\big(0;(w,u,v)\big)\, du\, dv\right] \f(w)\, dw $$ 
$$ = \E_0\Big[ \f \Big( \sqrt{\e}\,\o_1\Big) \Phi\Big( \e\! \int_0^{1} \cos\!\big(\sqrt{\e}\,\o_{\tau}\big)\, d\tau ,\, \e\! \int_0^{1} \sin\!\big(\sqrt{\e}\,\o_{\tau}\big)\, d\tau\Big)\Big]  $$ 
$$ = \int_{\sR} e^{-w^2/2\e}\; \E_0^{w/\sqrt{\e}}\Big[ \f \Big( \sqrt{\e}\,\o_1\Big) \Phi\Big( \e\! \int_0^{1} \cos\!\big(\sqrt{\e}\,\o_{\tau}\big)\, d\tau ,\, \e\! \int_0^{1} \sin\!\big(\sqrt{\e}\,\o_{\tau}\big)\, d\tau\Big)\Big]  \frac{dw}{\sqrt{2\pi\,\e}}  
$$ 
\begin{equation*}  \label{f.conditpont}
= \int_{\sR}\, \frac{e^{-w^2/2\e}}{\sqrt{2\pi\,\e}}\; \E_0^{0}\Big[ \Phi\Big( \e\! \int_0^{1} \cos\!\big(\tau w+\sqrt{\e}\,\o_{\tau}\big)\, d\tau ,\, \e\! \int_0^{1} \sin\!\big(\tau w+\sqrt{\e}\,\o_{\tau}\big)\, d\tau\Big)\Big] \, \f(w)\, dw \, , 
\end{equation*}   
so that 
$$ \int_{\sR^2}\! \Phi(u,v)\, p_\e\big(0;(w,u,v)\big) du dv $$
$$ = \frac{e^{-w^2\!/2\e}}{\sqrt{2\pi\,\e}}\; \E_0^{0}\Big[ \Phi\Big( \e\! \int_0^{1} \cos\!\big(\tau w+\sqrt{\e}\,\o_{\tau}\big)\, d\tau ,\, \e\! \int_0^{1} \sin\!\big(\tau w+\sqrt{\e}\,\o_{\tau}\big)\, d\tau\Big)\Big]. \;\;\diamond  $$ 

\subsection{Fourier-Plancherel Formula} \label{sec.FPGF} \indf 
   We here use the Fourier transform and the Plancherel inversion formula, as ([BA1], (3.9) p.\,319), but without potential function. Note that  the smooth density $(w,u,v)\mapsto p_\e\big(0\,;(w,u,v)\big)$ is bounded and integrable, hence square integrable. Then we dilate the integration variable according to $(\xi',\xi)\mapsto \big(\e^{-3/2}\,\xi', \e^{-3/2}\,\xi\big)$. Thus we deduce from Lemma \ref{lem.Bbdes} that we have\,:   
$$ 4\pi^2\e^{3}\sqrt{2\pi\,\e}\times e^{w^2/(2 \e)}\times\, p_\e\big(0\, ; (w,\e\, y\,, \e\, z)\big)  \hskip1mm \vspace{-1mm}  $$
$$ =\, \e^{3}\, \sqrt{2\pi\e}\,\, e^{w^2/(2 \e)} \int_{\sR^2}\! \bigg[\!\int_{\sR^2}\! e^{\rt1\! (\xi'u+\xi v) }\, p_\e\big(0\,; (w,u,v)\big)\, du\, dv\bigg] e^{-\rt1\! \left(\xi' \e \, y+\xi\, \e \,z\right)}\, d\xi'd\xi $$ 
$$ = \int_{\sR^2}\! \bigg[\!\int_{\sR^2}\!  e^{\rt1\! \big(\frac{\xi'u}{\e^{3/2}}+ \frac{\xi v}{\e^{3/2}}\!\big) } \sqrt{2\pi \e}\,e^{\frac{w^2}{2\e}}\, p_\e\big(0\,; (w,u,v)\big)\, du dv \bigg] e^{-\rt1\! \big(\!\frac{\xi' y}{\e^{1/2}}+ \frac{\xi\, z}{\e^{1/2}}\!\big)} d\xi'd\xi \, $$ 
{ \small 
$$ =  \int_{\sR^2} \E_0^{0}\bigg[\! \exp\!\bigg(\!{\ts\frac{\rt1}{\sqrt{\e}}}\! \bigg({\xi'} \left[ \int_0^{1}\! \cos\!\big(ws +\sqrt{\e}\,\o_{s}\big) ds - y\right] 
+ {\xi} \left[\int_0^{1}\! \sin\!\big(ws +\sqrt{\e}\,\o_{s}\big) ds - z\right]\bigg) \!\bigg) \bigg] d\xi'd\xi\, .  
$$ }

   We re-state this as follows. 
\blem \label{lem.calcavpontFb} \   For any $\,\e>0$ and $(w,y,z)\in\R^3$, we have\,:  
\begin{equation}  \label{f.expr7} 
p_\e\big(0\, ; (w,\e\, y\,, \e\, z)\big) =\, \frac{e^{ -w^2/(2 \e)}}{4\pi^2\,\e^{3}\,\sqrt{2\pi\,\e}}\,  \int_{\sR^2} P_\e(\xi',\xi)\, d\xi'd\xi\,,  \vspace{-1mm}
\end{equation}  
with \vspace{-0.5mm}
{\small 
$$ P_\e(\xi',\xi) \, := \, \E_0^{0}\bigg[\! \exp\!\bigg(\!{\ts\frac{\rt1}{\sqrt{\e}}}\! \bigg({\xi'} \left[ \int_0^{1}\! \cos\!\big(ws +\sqrt{\e}\,\o_{s}\big) ds - y\right] + {\xi} \left[\int_0^{1}\! \sin\!\big(ws +\sqrt{\e}\,\o_{s}\big) ds - z\right]\bigg) \!\bigg) \bigg] .  \vspace{-1mm}  $$ } 
\elem  
This lemma incites to compare $\,P_\e(\xi',\xi) $ with 
$$ \overline{P}_\e(\xi',\xi) := \, \exp\!\left[{ {\frac{\ss\rt1}{\sqrt{\e}} \left({\xi'} \left[ \int_0^{1}\! \cos(ws) ds - y\right] + {\xi} \left[ \int_0^{1}\! \sin(ws) ds - z\right] \right)}}\right]  \times \vspace{-2mm}\quad  $$
$$ \qquad \times \, \E_0^{0}\bigg[ \exp\!\bigg(\!{\rt1}\!\! \int_0^{1} \big(\xi\, \cos(ws) - \xi' \sin(ws)\big)\,  \o_s\, ds \bigg)\bigg] . $$ 
However this comparison needs careful estimates, which must overcome integration against $\,d\xi'd\xi\,$, and constitute the content of the following section \ref{sec.changbsuit1scomp}. Though the replacement of $\,P_\e(\xi',\xi) $ by $\,\overline{P}_\e(\xi',\xi)$ seems a priori reasonable, we have indeed to be very careful, owing to the high sensibility of oscillatory integrals\,: see for example Section 4 of [I-M]. 

\section{Analysis and domination of $\,P_\e(\xi',\xi)$}  \label{sec.changbsuit1scomp} \indf 
   By Lemma \ref{lem.calcavpontFb} and the order 1 Taylor formula we have  
{
\begin{equation}  \label{f.Pe()} 
P_\e(\xi',\xi) = \exp\!\left[{ {\frac{\ss\rt1}{\sqrt{\e}}\! \left({\xi'} \left[ \int_0^{1}\! \cos(ws) ds - y\right] + {\xi} \left[ \int_0^{1} \!\sin(ws) ds - z\right] \right)}}\right]\! \times \EE_\e(\xi',\xi) ,\vspace{-0mm}
\end{equation} }
with \vspace{-2mm} 
\begin{equation}  \label{df.Eeps} 
\EE_\e(\xi',\xi) :=\, \E_0^{0}\Big[ \exp\!\left({{\rt1}\! A_\e(\xi',\xi)}\right) \Big]  
\end{equation}
and 
\begin{equation} \label{dfAe}
A_\e(\xi',\xi) := \int_0^{1}\!\int_0^{1}\! \big[ \xi\, \cos\!\big(ws+\!\sqrt{\e}\,\tau\, \o_s\big) - {\xi'} \sin\!\big(ws+\!\sqrt{\e}\,\tau\, \o_s\big) \big]\, \o_s\, ds\, d\tau\,  . 
\end{equation}
\parsn 
Let us decompose the phase according to \  $A_\e(\xi',\xi) = \xi'\,U_1^\e + \xi\, U_2^\e\,$, \  where for $\,\e\ge 0\,$: \vspace{-0mm} 
\begin{equation} \label{fde.Ve} 
U^\e(\o) := \left(-\int_{[0,1]^2}\! \sin\!\big(ws +\tau \sqrt{\e}\,\o_{s}\big)\,d\tau\; \o_s  \, ds \, , \int_{[0,1]^2}\! \cos\!\big(ws +\tau \sqrt{\e}\,\o_{s}\big)\,d\tau\; \o_s  \, ds  \right) \! .   \vspace{1mm}
\end{equation}   

   Thus we have to analyse the behaviour of $\,\EE_\e(\xi',\xi)$ as $\e\sea 0\,$.\,  As the term $\,A_\e(\xi',\xi)$ is clearly almost surely continuous with respect to $\e\,$, this amounts to dominating $\,\EE_\e(\xi',\xi)$ with respect to the Lebesgue measure $\,d\xi'd\xi\,$ on $\R^2$. \pars 
   To proceed, we shall use Theorem  2.1 of [T2], viewing $\,\EE_\e(\xi',\xi)$ as an infinite-dimen\-sional oscillatory integral with phase $A_\e(\xi',\xi)$. This requires a careful analysis of the (non-deterministic) Malliavin matrix $DU^\e$ of $U^\e$. 
 
\subsection{Analysis of the phase vector $\;{\ds U^\e}$ {\rm (defined by (\ref{fde.Ve}))}} \label{handlerrort}  \indf
    In order to control the error term by a $d\xi'd\xi$-integrable asymptotically vanishing term, we resort to Malliavin's  calculus, as in [BA1], to take advantage of the oscillatory nature of the integrand, which has phase $A_\e(\xi',\xi)$. Namely we shall perform the analogue of Lemma (3.48) in [BA1], in order to obtain a  $\,d\xi'd\xi$-integrable control on $\, P_{\e}(\xi',\xi)$. However the non-deterministic nature of the Malliavin matrix $DU^{\e}$ we have to face here causes difficulties that did not arise in [BA1]. We shall thus use Malliavin's Calculus in another way than [BA1], namely with the help of [T2] in order to fully take advantage of the oscillatory nature of the error term (which crude $L^p$ estimates would ruin).    \par 

\subsubsection{Integration by parts with respect to $\P_0^0$} \label{MalliavCalc} \indf  
   Let us use integration by parts with respect to the pinned Wiener measure $\P_0^0\,$. \parn
For an accessible reference on Malliavin's Calculus, see for example the book [Fa] or [M]. \parm  

   Classically denoting by $\WW_0$ the pinned Wiener space of continuous real maps $\o$ on $[0,1]$ such that $\,\o_0=\o_1=0\,$ and by $H_0$ the associated Cameron-Martin space, we successively have \big($F,G$ denoting real-valued elements of $D^2_1(\WW_0)$, and $Z$ any element of $C(\WW_0\,,H_0)$\big): \parn
the gradient  $\nabla$,  
defined for $\,h\in H_0\,, \o\in \WW_0\,$ by $\langle \nabla F(\o),h\rangle_{H_0} \equiv \nabla_h F(\o) := \frac{d_o}{dt} F(\o+th)\,$; \parn 
the divergence operator or Skorohod integral $\delta\,$, given by $\,{\ds \delta(Z)(\o) := \int_0^1\! \dot Z_s(\o) D\o_s}\,$, and in particular\,: $\,{\ds \delta(Fh)(\o) = F(\o)\! \int_0^1\! \dot h(s)\, d\o_s -\nabla_hF }\,$, $\,{\ds \delta\Big(\int_0^\cdot\! f(t,\o_t)\,dt\Big) = \int_0^1\! f(t,\o_t)\, d\o_t}$\,; \parn
the integration by parts formula\,: 
${\ds \E_0^0\big[ \langle \nabla F,\nabla G\rangle_{H_0} \big] = \E_0^0\big[ \delta (\nabla F) \times G\big] }$\,; \parn the Ornstein-Uhlenbeck operator $\LL\,$, given by $\,\LL F := \, \delta (\nabla F)$\,; \ 
the energy estimate (see for example [Fa] page 57) \quad {${\ds  \E_0^0\big[ (\LL F)^2\big] \le 4\, \E_0^0\big[ \|\nabla^2 F\|^2_{H_0\otimes H_0}+ \|\nabla F\|_{H_0}^2 \big]}$\,.} \par 

\subsubsection{Gradients and Malliavin matrix of $\,U^\e$} \label{MalliavCalcVe} \indf  
   According to (\ref{fde.Ve}), for any $\,h\in H_0\,$ we have\,:   \vspace{-2mm}   
$$ \big\langle \nabla U_2^{\e}(\o), h\big\rangle_{H_0} = \nabla\!_{h} U_2^{\e}(\o) =  \int_{[0,1]^2}\!\big[\cos(ws\! +\!\sqrt{\e}\,\tau\,\o_{s}) - \!\sqrt{\e}\,\tau\, \o_s\, \sin(ws\! +\!\sqrt{\e}\,\tau\,\o_{s})\big] h_s\, ds \, d\tau\,   \vspace{-2mm} $$
$$ =   \int_0^1\!  \cos(ws\! +\!\sqrt{\e}\,\o_{s}) \, h_s\, ds = \int_0^1\left[\int_t^1\! \cos(ws\!+\!\sqrt{\e}\,\o_{s}) \, ds \right] \dot h_t\, dt \, ,  \vspace{-2mm} $$
i.e., \vspace{-2mm} 
$$ \frac{d}{dt} \nabla U_2^{\e}(\o)(t) = \int_t^1\!\cos\!\big(ws\! +\!\sqrt{\e}\,\o_{s}\big) ds - \int_0^1\! du\!\int_u^1\!\cos\!\big(ws\! +\!\sqrt{\e}\,\o_{s}\big) ds \, =\, \frac{d}{dt} \nabla U_2^1(\!\sqrt{\e}\,\o)(t) \; ; $$ 
similarly, \vspace{-2mm} 
$$  \frac{d}{dt} \nabla U_1^{\e}(\o)(t) = \int_0^1\! du\!\int_u^1\!\sin\!\big(ws\! +\!\sqrt{\e}\,\o_{s}\big) ds - \int_t^1\!\sin\!\big(ws\! +\!\sqrt{\e}\,\o_{s}\big) ds \, =\, \frac{d}{dt} \nabla U_1^1(\!\sqrt{\e}\,\o)(t) \; ; \vspace{-1mm}  $$ 
and \vspace{-2mm} 
{\small 
$$ \frac{d}{dt} \nabla U^0(\o)(t) = \left( { \frac{\sin w}{w^2}  - \frac{\cos(wt)}{w}}\; \raise1.4pt\hbox{,}\; { \frac{1-\cos w}{w^2}  - \frac{\sin(wt)}{w}}\right) = \frac{d}{dt} \nabla U^0(0)(t) = \frac{d}{dt} \nabla U^\e(0)(t)\, .  $$ }
Then the Malliavin covariance matrix is\,: 
\begin{equation} \label{fde.MallMatrix0}�
DU^{\e}(\o) :=  \begin{pmatrix} \langle \nabla U^{\e}_1, \nabla U^{\e}_1\rangle_{H_0} & \langle \nabla U^{\e}_1, \nabla U^{\e}_2\rangle_{H_0} \cr  \langle \nabla U^{\e}_1, \nabla U^{\e}_2\rangle_{H_0} & \langle \nabla U^{\e}_2, \nabla U^{\e}_2\rangle_{H_0} \end{pmatrix}\! (\o) = DU^1(\!\sqrt{\e}\,\o) \, , \vspace{-2mm} 
\end{equation} 
with \vspace{-2mm} 
$$ \langle \nabla U^{\e}_1,\nabla U^{\e}_1\rangle_{H_0} = \int_0^1\! \left(\int_t^1\!\sin\!\big(ws\! +\!\sqrt{\e}\,\o_{s}\big) ds\right)^{\! 2}\! dt - \left[\int_0^1\! \left( \int_t^1\!\sin\!\big(ws\! +\!\sqrt{\e}\,\o_{s}\big) ds\right) dt \right]^{2} \le \frac{1}{3}\, \raise1.5pt\hbox{,} $$ 
$$ \langle \nabla U^{\e}_2,\nabla U^{\e}_2\rangle_{H_0} = \int_0^1\!\left(\int_t^1\!\cos\!\big(ws\! +\!\sqrt{\e}\,\o_{s}\big) ds\right)^{\! 2}\! dt - \left[\int_0^1\!\left( \int_t^1\!\cos\!\big(ws\! +\!\sqrt{\e}\,\o_{s}\big) ds\right) dt \right]^{2} \le \frac{1}{3}\, \raise1.5pt\hbox{,} $$ 
so that \quad  $\,{\Tr}(DU^{\e}) \ge 2\sqrt{\det(DU^{\e})}\, \ge 6\, {\det(DU^{\e})} \,$,  \quad and \quad 
$$ \langle \nabla U^{\e}_1, \nabla U^{\e}_2\rangle_{H_0}\, = \int_0^1\! \left[\int_t^1\! \cos(ws\! +\!\sqrt{\e}\,\o_s)\, ds\right]\! dt\times \int_0^1\! \left[\int_t^1\! \sin(ws\! +\!\sqrt{\e}\,\o_s)\, ds \right]\! dt\hskip10mm \vspace{-2mm} $$
$$ \hskip30mm - \int_0^1\! \left[\int_t^1\! \cos(ws\! +\!\sqrt{\e}\,\o_s)\, ds \right]\!\times \! \left[\int_t^1\! \sin(ws\! +\!\sqrt{\e}\,\o_s)\, ds \right]\! dt \, . $$  
\brem \label{rem.detDV0} \  {\rm   In particular, at $\,\e=0\,$ we have the following deterministic limiting Malliavin covariance matrix\,: 
\begin{equation*}
DU^{0} :=  \begin{pmatrix} \frac{1}{2w^2}\!\left(1+\frac{\sin(2w)}{2w}\right) - \frac{\sin^2\!w}{w^4} & \frac{\sin^2\!w}{2w^3} - \frac{(1-\cos w)\sin w}{w^4} \cr  \frac{\sin^2\!w}{2w^3} - \frac{(1-\cos w)\sin w}{w^4} & \frac{1}{2w^2}\!\left(1-\frac{\sin(2w)}{2w}\right) - \frac{(1-\cos w)^2}{w^4} \end{pmatrix}\!  , \vspace{-0mm} 
\end{equation*} 
and then 
$$ \det(D U^0) = \det(D U^0)(0) = \det(DU^\e)(0) $$
$${\ts  = \left[\frac{1}{2w^2}\!\left(1-\frac{\sin(2w)}{2w}\right) - \frac{(1-\cos w)^2}{w^4}\right]\!\times\!\left[ \frac{1}{2w^2}\!\left(1+\frac{\sin(2w)}{2w}\right) - \frac{\sin^2\!w}{w^4}\right] - \left[\frac{\sin^2\!w}{2w^3} - \frac{(1-\cos w)\sin w}{w^4}\right]^2 } $$ 
$$ = \frac{1}{4w^4} - \frac{4(1-\cos w)+ \sin^2\!w}{4\,w^6} + \frac{(1-\cos w)\sin w}{w^7}  = \frac{w^2}{8640} + \O(w^4)\, \hbox{ (for small $|w|$)} . $$  
This  determinant $\,\det(D U^0)$ would vanish \iff both centred functions $\,{ \frac{\sin w}{w^2}  - \frac{\cos(wt)}{w}}$ and $\,{ \frac{1-\cos w}{w^2}  - \frac{\sin(wt)}{w}}\,$ were proportional\,; but this would imply proportionality between $\,t\mapsto \cos(wt)\,$ and $\,t\mapsto \sin(wt)$, which does not hold provided $\,w\not=0$. \  Hence $\,\det(D U^0) >0\,$ for any $\,w\not=0$. For the same reason, and since $\P_0^0$-almost surely $\,\o \not= [s\mapsto -ws/\sqrt{\e}\,]$, we almost surely have $\,\det(D U^\e) >0$, for any $\,\e\ge 0$.
}\erem 

   Then we have \vspace{-2mm} 
\begin{equation} \label{f.LL ue}
\LL U^{\e}_2(\o) = \int_0^1 \frac{d}{dt} \nabla U_2^{\e}(\o)(t)\, d\o_t\, = \int_0^1\!\cos\!\big(ws +\!\sqrt{\e}\,\o_{s}\big) \,\o_s \, ds\, , \vspace{-1mm} 
\end{equation}
and similarly \quad ${\ds \LL U^{\e}_1(\o) =\, - \int_0^1\!\sin\!\big(ws +\!\sqrt{\e}\,\o_{s}\big) \,\o_s \, ds\,}$. \parsn 
Thence \vspace{-2mm} 
$$  \frac{d}{dt} \nabla \LL U_1^{\e}(\o)(t) = \int_0^t\! \big[\sqrt{\e}\,\omega_s \cos(ws\! + \!\sqrt{\e}\,\o_{s}) + \sin(ws\! +\!\sqrt{\e}\,\o_{s})\big] ds - \int_0^1\!\int_0^t\hbox{idem} \,  ; \vspace{-1mm}  $$ 
$$  \frac{d}{dt} \nabla \LL U_2^{\e}(\o)(t) = \int_0^t\! \big[\sqrt{\e}\,\omega_s \sin(ws\! +\!\sqrt{\e}\,\o_{s}) -  \cos(ws\! +\!\sqrt{\e}\,\o_{s})\big] ds - \int_0^1\!\int_0^t\hbox{idem} \,  ; \vspace{-1mm}  $$ 
$$ \big\| \nabla \LL U_1^{\e} \big\|^2_{H_0}(\o) = {\rm Var}_{[0,1]}\left[\int_0^{\bullet}\! \big[\sqrt{\e}\,\omega_s \cos(ws\! + \!\sqrt{\e}\,\o_{s}) + \sin(ws\! +\!\sqrt{\e}\,\o_{s})\big] ds \right] \le 1+ \e\! \int_0^1\!\omega^2\, ; $$ 
$$ \big\| \nabla \LL U_2^{\e} \big\|^2_{H_0}(\o) = {\rm Var}_{[0,1]}\left[\int_0^{\bullet}\! \big[\sqrt{\e}\,\omega_s \sin(ws\! +\!\sqrt{\e}\,\o_{s}) -  \cos(ws\! +\!\sqrt{\e}\,\o_{s})\big] ds  \right] \le 1+ \e\! \int_0^1\!\omega^2 . $$ 

\subsubsection{Estimates for the successive gradients of $\,U^\e$ and $\,\LL U^\e$} \label{successgrads} \indf  
We begin with an explicit formula for $\,\nabla^kU_2^\e\,$. \vspace{-1mm}
\blem \label{lem.estimUe} \  For $\,k\ge 1\,$, setting \vspace{-1.5mm} 
$$ g_k^\o(t) := (\sqrt{\e}\,)^{k-1}\! \int_t^1 \cos^{(k-1)}\!\big(wu +\!\sqrt{\e}\,\o_u\big)\, du\, , \vspace{-2mm}  $$ 
we almost surely have (recall the definition (\ref{fde.Ve} ) of $U^\e$)\,: \vspace{-1.5mm}  
 $$ \frac{d^k}{ds_1...ds_k} \nabla^kU_2^\e(\o)\big(s_1,\ldots,s_k\big) =\, \sum_{j=0}^k (-1)^j\! \sum_{1\le i_1<...<i_j\le k}  \int_{[0,1]^j} g_k^\o(\max\{s_1,...,s_k\})\, ds_{i_1}...ds_{i_j}\, . \vspace{-1.5mm}  $$  
\elem  
\ub{Proof} \quad  We proceed by induction on $k\,$. We already saw the case $k=1$ in the beginning of Section \ref{MalliavCalcVe}. Then set $\,\delta_k := \frac{d^k}{ds_1...ds_k} \nabla^kU_2^\e(\o)$, denote by $\,\Delta_k\,$ the right hand side of the above formula, and assume that $\,\delta_k\big(s_1,\ldots,s_k\big) = \Delta_k\,$. Then for any $\,h\in H_0^1\,$ and any $\,t\in [0,1]$ we have \vspace{-2mm}  
$$ \nabla_h\, g_k^\o(t)\, =\, \e^{k/2}\! \int_t^1 \cos^{(k)}\!\big(wu +\!\sqrt{\e}\,\o_u\big)\, h(u)\, du\, \vspace{-1.5mm}  $$
$$ =\, \e^{k/2} \! \int_0^1 \int_s^11_{[t,1]}(u) \cos^{(k)}\!\big(wu +\!\sqrt{\e}\,\o_u\big)\, du\; \dot h(s)\, ds\,  = \int_0^1 g_{k+1}^\o(\max\{s,t\}) \, \dot h(s)\, ds\, . $$ 
Hence, 
$$ \int_0^1 \delta_{k+1}\big(s_1,\ldots,s_{k+1}\big)\, \dot h(s_{k+1})\, ds_{k+1}\, = \nabla_h\big[\delta_k\big(s_1,\ldots,s_k\big)\big] = \nabla_h\,\Delta_k $$ 
$$ = \sum_{j=0}^k (-1)^j\! \sum_{1\le i_1<...<i_j\le k}  \int_{[0,1]^j}\! \nabla_h\, g_k^\o(\max\{s_1,...,s_k\})\, ds_{i_1}...ds_{i_j}\, $$
$$ = \int_0^1 \sum_{j=0}^k (-1)^j\! \sum_{1\le i_1<...<i_j\le k}  \int_{[0,1]^j} g_{k+1}^\o(\max\{s_1,...,s_{k+1}\})\, ds_{i_1}...ds_{i_j}\, \dot h(s_{k+1})\, ds_{k+1} \, , $$ 
so that \vspace{-1.5mm}  
$$ \delta_{k+1}\big(s_1,\ldots,s_{k+1}\big) = \sum_{j=0}^k (-1)^j\! \sum_{1\le i_1<...<i_j\le k}  \int_{[0,1]^j} g_{k+1}^\o(\max\{s_1,...,s_{k+1}\})\, ds_{i_1}...ds_{i_j}\, \hskip 14mm \vspace{-2mm}  $$ 
$$ \hskip 30mm -  \int_0^1 \left[ \sum_{j=0}^k (-1)^j\! \sum_{1\le i_1<...<i_j\le k}  \int_{[0,1]^j} g_{k+1}^\o(\max\{s_1,...,s_{k+1}\})\, ds_{i_1}...ds_{i_j}\right]\! ds_{k+1} $$ 
$$ =\, \sum_{j=0}^k (-1)^j\! \sum_{1\le i_1<...<i_j\le k} \left[ \int_{[0,1]^j}\! g_{k+1}^\o(\max\{s_1,...,s_{k+1}\})\, ds_{i_1}...ds_{i_j} \right.  \hskip 35mm \vspace{-2.5mm}  $$
$$ \hskip 65mm \left.  -  \int_{[0,1]^{j+1}}\! g_{k+1}^\o(\max\{s_1,...,s_{k+1}\})\, ds_{i_1}...ds_{i_j}\, ds_{k+1} \right] $$ 
$$ =\, \Delta_{k+1}\, . \;\;\diamond $$ 
Of course, the same holds with `$\,\sin\,$' instead of `$\,\cos\,$'. This entails the following estimate. 
\blem \label{lem.estimA} \  For $\,k\ge 1\,$ and $\,j=0,1\,$ we almost surely have\,: \vspace{-1.5mm}  
$$ \|\nabla^kU_j^\e\|_{H_0^{\otimes k}}\, \le\, 2\sqrt{2}\, k\1\, (2\sqrt{\e}\,)^{k-1}\, , \vspace{-2mm}   $$    
and (recall (\ref{dfAe}))\,: 
$$ \|\nabla^kA_\e(\xi',\xi)\|_{H_0^{\otimes k}}\, \le \,4 \sqrt{\xi^2+\xi'\,\!^2}\; k^{-1}\, (2\sqrt{\e}\,)^{k-1}\, . \vspace{-1mm} $$    
\elem  
\ub{Proof} \quad  The latter follows at once from the former, owing to (\ref{dfAe}).  Now by Lemma \ref{lem.estimA} for any $\,k\ge 1\,$ we have \  $\big| g_k^\o(t)\big| \le \,\e^{\frac{k-1}{2}} \, (1-t)\,$, \  and 
$$ \|\nabla^kU_j^\e\|_{H_0^{\otimes k}}^2\, = \int_{[0,1]^k}\! \delta_k\big(s_1,\ldots,s_k\big)^2 \, ds_1 \ldots ds_k\, =\int_{[0,1]^k}\! \Delta_k^2\, ds_1 \ldots ds_k\, $$
$$ \le\, 4^k \int_{[0,1]^k} g_k^\o(\max\{s_1,...,s_k\})^2 \, ds_1 \ldots ds_k\,  \le\, 4^k \e^{{k-1}} \!\int_{[0,1]^k} (1-\max\{s_1,...,s_k\})^2 \, ds_1 \ldots ds_k  $$ 
$$  =\, 4^k \e^{{k-1}}\times \frac{2}{(k+1)(k+2)}\,<\, 8\,k\2 \,(4\,{\e}\,)^{{k-1}} \hbox{.} \;\;\diamond $$

   We can proceed analogously to handle  $\,\nabla^k\LL U_2^\e\,$. \vspace{-1mm} 
\blem \label{lem.estimLUe} \  For $\,k\ge 1\,$, setting \vspace{-1.5mm} 
$$ \tilde g_k^\o(t) := (\sqrt{\e}\,)^{k-1}\! \int_t^1\!\Big[k \,\cos^{(k-1)} +\,\sqrt{\e}\,\o_u \cos^{(k)}\!\Big] (wu + \sqrt{\e}\,\o_u)\, du\, , \vspace{-2mm}  $$ 
we almost surely  have\,: \vspace{-1mm}  
 $$ \frac{d^k}{ds_1...ds_k} \nabla^k\LL U_2^\e(\o)\big(s_1,\ldots,s_k\big) =\, \sum_{j=0}^k (-1)^j\! \sum_{1\le i_1<...<i_j\le k}  \int_{[0,1]^j} \tilde g_k^\o(\max\{s_1,...,s_k\})\, ds_{i_1}...ds_{i_j}\, . \vspace{-2mm}  $$  
\elem  
\ub{Proof} \quad  We proceed as for the above lemma \ref{lem.estimUe}. The case $k=1$ reduces to (\ref{f.LL ue}) in Section \ref{MalliavCalcVe}. Then the induction on $k\,$ is as in the proof of Lemma \ref{lem.estimUe}, except that the computation of $ \nabla_h\, g_k^\o(t)$ must be replaced by the following\,: 
$$ \nabla_h\, \tilde g_k^\o(t)\, =\, \e^{k/2}\! \int_t^1\!\Big[(k+1) \,\cos^{(k)} +\,\sqrt{\e}\,\o_u \cos^{(k+1)}\!\Big] \!\big(wu +\!\sqrt{\e}\,\o_u\big)\, h(u)\, du\, \vspace{-1.5mm}  $$
$$ =\, \e^{k/2} \! \int_0^1 \int_s^11_{[t,1]}(u) \Big[(k+1) \,\cos^{(k)} +\,\sqrt{\e}\,\o_u \cos^{(k+1)}\!\Big]\!\big(wu +\!\sqrt{\e}\,\o_u\big)\,du\;  \dot h(s)\, ds\,  $$
$$ = \int_0^1 \tilde g_{k+1}^\o(\max\{s,t\}) \, \dot h(s)\, ds\, . \;\;\diamond $$ 
Of course, the same holds with `$\,\sin\,$' instead of `$\,\cos\,$'. This entails the following estimate. 
\blem \label{lem.estimLLUe} \  For $\,k\ge 1\,$ and $\,j=0,1\,$ we almost surely have\,: \vspace{-1.5mm}  
$$ \|\nabla^k\LL U_j^\e\|_{H_0^{\otimes k}}\, \le\, 2\sqrt{2}\;  (2\sqrt{\e}\,)^{k-1}\times \left[ 2 + \e\!\int_0^1\!\o^2 \right]^{1/2} . \vspace{-2mm}   $$    
\elem  
\ub{Proof} \quad More or less as for Lemma \ref{lem.estimA}, but owing now to Lemma \ref{lem.estimLUe}, on the one we have \parn 
\centerline{${\ds  \big|\tilde g_k^\o(t)\big| \le (\sqrt{\e}\,)^{k-1} \int_t^1\sqrt{k^2+\e\,\o^2_u}\,du\, \le (\sqrt{\e}\,)^{k-1} k\, (1-t) + (\sqrt{\e}\,)^{k} \int_0^1 |\o_u|\, du\, }$,} \parsn
and on the other hand, 
$$ \|\nabla^k \LL U_j^\e\|_{H_0^{\otimes k}}^2\,  \le\, 4^k \int_{[0,1]^k} \tilde g_k^\o(\max\{s_1,...,s_k\})^2 \, ds_1 \ldots ds_k\, $$
$$ \le\, 4^k \e^{{k-1}}\times 2 \int_{[0,1]^k} \left( k^2\, (1-\max\{s_1,...,s_k\})^2 + \e \int_0^1 \o^2_\cdot \right)   \, ds_1 \ldots ds_k  $$ 
$$ <\, (4\,\e)^{{k-1}}\times 8\left( 2 + \e \int_0^1 \o^2_\cdot \right) . \;\;\diamond  $$ 

   As a consequence, for the random variables  
\begin{equation}  \label{df.M12} 
M_1[\f,r] := \sum_{n\ge 1}\, \frac{r^n}{n!}\, \big\| \nabla^n \f \big\|^2_{H_0^{\otimes n}}\quad \hbox{and} \quad M_2[\f,r] := \sum_{n\ge 2}\, \frac{r^n}{n!}\,\big\| \nabla^n \f \big\|^2_{H_0^{\otimes n}}\,
\end{equation}
which will intervene in the following section \ref{sec.IPPT2}, we directly deduce the following. 
\blem \label{lem.estimMu} \  For any positive $\,r,\e\,$ we almost surely have\,: \vspace{-0mm}  
$$  M_1[U_1^{\e},r] +  M_1[U_2^{\e},r] < 16\, r\, e^{4r\,\e} \quad \hbox{and} \quad M_2[U_1^{\e},r] +  M_2[U_2^{\e},r] < 16\,r^2\e\, (e^{4r\,\e}-1) \, ; $$ 
$$  M_1[\LL U_1^{\e},r] +  M_1[\LL U_2^{\e},r] < 16\,r\, e^{4r\,\e}\!\left[ 2 + \e\!\int_0^1\!\o^2 \right] ; $$
$$ M_2[\LL U_1^{\e},r] +  M_2[\LL U_2^{\e},r] < 16\,r^2\e\, (e^{4r\,\e}-1)\!\left[ 2 + \e\!\int_0^1\!\o^2 \right] . \vspace{-2mm}  $$ 
\elem  
\vspace{-5mm}  

\subsection{Applying Theorem 2.1 of [T2]} \label{sec.IPPT2} \indf 
   We use the treatment of infinite-dimensional oscillatory integrals as developed in  [T2], to estimate $\,\EE_s(\xi',\xi)$, recall (\ref{df.Eeps}), (\ref{dfAe}).  The correspondence between the notation used by Taniguchi S.  in [T2] and ours is given by\,:  \parsn
\centerline{$N=1$, $\psi=1$,  $\, q_0=U_1^\e\,,\, q_1 = U_2^\e\,$, $\lambda = \sqrt{\xi^2+\xi'\!\,^2}\,$, $\,z^0= \frac{\xi}{\lambda}\,$\raise1pt\hbox{,} $\,z^1= \frac{\xi'}{\lambda}\,$\raise1pt\hbox{,} $\, A_\e(\xi',\xi) = \lambda\, q^{(z)}$.} \parsn 
Thus \quad  $\big\| \nabla q^{(z)} \big\|^2_{H_0} = \big\| z^1\,\nabla U_1^{\e} + z^0\,\nabla U_2^{\e} \big\|^2_{H_0} = (z^1,z^0)(DU^{\e})\,^t(z^1,z^0)$, \  so that 
$$ \min\!\Big\{ \big\| \nabla q^{(z)} \big\|^2_{H_0} \Big|\, |z|=1\Big\}\, =\, \frac{\ss \| \nabla U_1^{\e}\|^2_{H_0} + \| \nabla U_2^{\e}\|^2_{H_0} - \sqrt{\left(\| \nabla U_1^{\e}\|^2_{H_0} - \| \nabla U_2^{\e} \|^2_{H_0}\right)^2 + 4 \left\langle \nabla U^{\e}_1, \nabla U^{\e}_2\right\rangle^2_{H_0}} }{2} $$
\begin{equation} \label{f.estgradq2} 
= \frac{2\, \det(DU^{\e})}{{\ss  \| \nabla U_1^{\e}\|^2_{H_0} + \| \nabla U_2^{\e}\|^2_{H_0} + \sqrt{\left(\| \nabla U_1^{\e}\|^2_{H_0} - \| \nabla U_2^{\e} \|^2_{H_0}\right)^2 + 4 \left\langle \nabla U^{\e}_1, \nabla U^{\e}_2\right\rangle^2_{H_0}} } }  \ge \,  \frac{\det(DU^{\e})}{  \Tr(DU^{\e}) } \,\ge \, \det(DU^{\e})\,\raise0pt\hbox{.} 
\end{equation} 
\pars  

   Then the estimate (2.1) in Theorem 2.1 of [T2] provides the following\,: \parsn
There exists $\,C>0\,$ such that for any $\,r\in\,]1,e[\,$ and $\xi^2+\xi'\!\,^2 >1\,$, we have 
$$ \big|\EE_\e(\xi',\xi)\big|^2 \le\, C\; \E_0^{0}\bigg[\exp\!\bigg(\! \frac{L_2[U^{\e}\!,r]}{L_1[U^{\e}\!,1]^3}\!\bigg) \bigg] \times \E_0^{0}\bigg[\!\exp\!\bigg(\! -  \frac{\min\{\frac{\det(DU^{\e})}{\Tr(DU^{\e})},1\}}{147\,C\,L_1[U^{\e}\!,1]} \sqrt{\xi^2+\xi'\!\,^2}\,\bigg)\bigg]  $$ 
hence according to Section \ref{MalliavCalcVe} \big(recall that $\,{\Tr}(DU^{\e}) \ge 2\sqrt{\det(DU^{\e})}\, \ge 6\, {\det(DU^{\e})}$\big)\,:  \vspace{1mm}  
\begin{equation} \label{f.estT2} 
\big|\EE_\e(\xi',\xi)\big|^2 \le\, C\; \E_0^{0}\bigg[\exp\!\bigg(\! \frac{L_2[U^{\e}\!,r]}{L_1[U^{\e}\!,1]^3}\!\bigg) \bigg] \times \E_0^{0}\bigg[\!\exp\!\bigg(\! -  \frac{\det(DU^{\e})/\Tr(DU^{\e})}{147\,C\,L_1[U^{\e}\!,1]} \sqrt{\xi^2+\xi'\!\,^2}\,\bigg)\bigg]  , 
\end{equation}  
where \ (for $\,i = 1,2$)  \quad  \vspace{-2mm}  
$$ L_i[U^{\e}\!,r] := \sum_{j=1}^2 \big(M_i[U_j^{\e},r] + M_i[\LL U_j^{\e},r]\big) + \bigg[1+ \sum_{j=1}^2\, (\LL U_j^{\e})^2\bigg]^{1/2} , \vspace{-2mm}  $$ 
with $\,M_1,\, M_2\,$ given by (\ref{df.M12}). \parn 
Now by the above lemma \ref{lem.estimMu} and (\ref{f.LL ue}), for $\,r \in [1,2]\,$ and $\,{\e}\,< 1\,$ we have\,: 
$$ 1\,\le\,  L_2[U^{\e}\!,r] \, \le\, L_1[U^{\e}\!,r]\, \le\, 16\,r\, e^{4r\,\e}\!\left[ 3 + \e\!\int_0^1\!\o^2 \right] + \left[ 1 + 2\!\int_0^1\!\o^2 \right]^{1/2}\,  \vspace{-1mm} $$ 
\begin{equation} \label{f.exprL1Ve1*}
\le \, \O(1) + \left[  2\!\int_0^1\!\o^2 \right]^{1/2}\! + \O(\e)\! \int_0^1\!\o^2\, \le\, \O(1) + \o^* + \O(\e) (\o^*)^2 =:\Omega_\e(\o)\, , 
\end{equation} 
where \  $\o^* := \max_{[0,1]}\limits  |\o|\,$.  \parn 
Now \big(see for example ([R-Y], I.(3.10)) or [B-O]\big) for any $\,t>0\,$ we have
\begin{equation} \label{f.estYor}
\P_0^0\big[\o^* >\sqrt{t}\,\big] = \P_0\bigg[\sup_{s\ge 0}\limits\, \big\{|\beta_s| - s\big\} > t\bigg] < 2\, e^{-2\,t}\, . 
\end{equation}
As a consequence, for some large enough constant $C$ we have\,: \vspace{-1mm} 
$$ \E_0^{0}\bigg[\exp\!\bigg(\! \frac{L_2[U^{\e}\!,r]}{L_1[U^{\e}\!,1]^3}\!\bigg) \bigg] \le\, \E_0^{0}\bigg[e^{\log \sqrt{C} +  \o^* + \e \log C\,  (\o^*)^2 }  \bigg] $$
$$ = \sqrt{C} \, \int_0^\infty  \P_0^0\Big[\o^* + \e \log C\,  (\o^*)^2 > \log u\Big] \, du\,= \sqrt{C} \, \int_{\sR} \P_0^0\bigg[\o^*> \frac{\sqrt{1+4s\,\e\log C } -1}{2\e\log C} \bigg] \, e^s\,ds\, $$ 
$$ < \sqrt{C} + 2 \sqrt{C} \, \int_{0}^\infty e^{s -\frac{\left(\sqrt{1+4s\,\e\log C} -1\right)^2}{2\e^2\log^2C}} \, ds\, = \sqrt{C} + 2 \sqrt{C} \, \int_{0}^\infty e^{\frac{\sqrt{1+4s\,\e\log C}}{\e^2\log^2C} - \left(\frac{2}{\e\log C}-1 \right) s}\, ds $$ 
$$ \le \, C  \quad \hbox{ for small enough $\,\e\,$.} $$ 
Then by (\ref{f.estT2}), for large enough positive constant $C$ and small enough positive $\,\e$ we have\,: 
\begin{equation} \label{f.estL2/L1} 
\big|\EE_\e(\xi',\xi)\big|^2 \le\, C^2\; \E_0^{0}\!\left[\exp\!\bigg(\! -  \frac{\det(DU^{\e})/\Tr(DU^{\e})}{147\,C\,L_1[U^{\e}\!,1]} \sqrt{\xi^2+\xi'\!\,^2}\,\bigg)\right] . 
\end{equation}

\subsection{Domination of the decreasing integral of (\ref{f.estL2/L1})}  \label{sec.changbsuit1sapple[T2]} \indf 
    We here establish the key domination property for the $\e$-dependent decreasing integral arising from \big([T2],Th.2.1(2.1)\big) and Section \ref{sec.IPPT2} above, in the non-singular case $\,w\not= 0\,$. Thus the aim of this section is to prove the following crucial estimate. \vspace{-1mm} 
\bpro \label{pro.estintdecre} \  For any fixed positive $C$ and $\,w\in\R^*$, we have 
{
$$ \sup_{0\le \e\le 1}\, \E_0^{0}\bigg[\exp\!\bigg(\! -{\frac{\det(DU^{\e})/\Tr(DU^{\e})}{147\,C\,L_1[U^\e,1]}} \sqrt{\xi^2+\xi'\!\,^2}\bigg)\!\bigg]  \in L^{1/2}(\R^2, d\xi'd\xi) \hbox{.} \vspace{-2mm} $$ } 
\epro 
\vspace{-2mm} 

\subsubsection{Beginning of the proof of Proposition \ref{pro.estintdecre}}  \indf 
   Set \quad $\RR_\e := {\ds \frac{\det(DU^{\e})}{\Tr(DU^{\e}) L_1[U^\e,1]}}\,$ \   and denote by $\,\lambda_\e\,$ the lowest eigenvalue of $\,DU^\e$, which is almost surely positive.    \   
By (\ref{f.exprL1Ve1*}) we have \quad  ${\ds \RR_\e \ge \,\frac{\det(DU^{\e})}{\Tr(DU^{\e})\,\Omega_\e } \, \ge \,\frac{\lambda_{\e}}{2\,\Omega_\e}\, \raise1.9pt\hbox{.}}$ \parn  
Now according to Subsection \ref {MalliavCalcVe}, for any $v=(v_1,v_2)\in \S^1\subset\R^2$ we have  \vspace{-2mm} 
$${ \langle v, DU^\e v\rangle =  {\rm Var}_{[0,1]}\!\left[\int_\bullet^1\!\big[ v_1 \sin\!\big(ws\! +\!\sqrt{\e}\,\o_{s}\big)- v_2 \cos\!\big(ws\! +\!\sqrt{\e}\,\o_{s}\big)\big] ds\right] }\,, \vspace{-2mm}  $$
so that 
$$ \lambda_{\e}\, =\, \inf_\theta\,  {\rm Var}_{[0,1]}\!\left[{\int_\bullet^1} \sin\!\big(ws\! +\! \theta\! +\!\sqrt{\e}\,\o_s\big)ds\right] .  $$
Hence, for any deterministic positive $R\,$ we have\,:  
\begin{equation}  \label{f.estimReVar} 
\E_0^{0}\!\left[ e^{-\RR_\e\times R}\,\right] \le\, \E_0^{0}\!\left[ \exp\!\left(-\, \frac{\lambda_{\e} R}{2\,\Omega_\e}\right)\right] = \,  \E_0^{0}\!\left[ e^{\frac{-R}{2\,\Omega_\e}\,\inf_\theta\limits\, {\rm Var}_{[0,1]}\!\left[{\ts \int_\bullet^1} \sin\left(ws+\theta+\sqrt{\e}\,\o_s\right)ds\right] }\right] . 
\end{equation}

   Then, for some positive constant $C$ and any $\,R>3/2\,, \, 0\le \e<1\,$,  we have\,: 
$$ \E_0^{0}\!\left[ e^{\frac{-R}{2\,\Omega_\e}\,\inf_\theta\limits\, {\rm Var}_{[0,1]}\!\left[{\ts \int_\bullet^1} \sin\left(ws+\theta+\sqrt{\e}\,\o_s\right)ds\right] }\right] \le\,  \E_0^{0}\!\left[ e^{ \frac{-\,R/C}{1+ \omega^*+ \e (\omega^*)^2}\,\inf_\theta\limits\, {\rm Var}_{[0,1]}\!\left[{\ts\int_\bullet^1} \sin\left(ws+\theta+\sqrt{\e}\,\o_s\right) ds\right] }\right] $$ 
$$ \le\, \P_0^{0}\!\left[ \omega^*> {R}^{1/3} \right] + \E_0^{0}\!\left[ e^{ \frac{-\,R/C}{1+ R^{1/3}+ \e\, R^{2/3}}\,\inf_\theta\limits\, {\rm Var}_{[0,1]}\!\left[{\ts \int_\bullet^1} \sin\left(ws+\theta+\sqrt{\e}\,\o_s\right) ds\right]} 1_{\{\omega^* \le R^{1/3}\}}  \right] $$ 
$$ \le \, 2\, e^{-2\, R^{2/3}} + \E_0^{0}\!\left[ e^{ \frac{-\, R^{1/3}}{\,2C }\,\inf_\theta\limits\, {\rm Var}_{[0,1]}\!\left[{\ts \int_\bullet^1} \sin\left(ws+\theta+\sqrt{\e}\,\o_s\right) ds\right] } \right]  . $$ 
Note that we applied (\ref{f.estYor}) to get the last line. \  Hence, according to (\ref{f.estimReVar}), so far we obtain\,: 
\begin{equation} \label{f.MajorIntOsc} 
\sqrt{\E_0^{0}\big[ e^{-\RR_\e\times R}\,\big]}\,  \le\, 2\, e^{-\, R^{2/3}} + \sqrt{\E_0^{0}\!\left[ e^{ -\,\frac{R^{1/3}}{2C}\,\inf_\theta\limits\, {\rm Var}_{[0,1]}\!\left[{\ts \int_\bullet^1} \sin\left(ws+\theta+\sqrt{\e}\,\o_s\right) ds\right] } \right]} . \vspace{-2mm} 
\end{equation} 
\vspace{-3mm} 

\subsubsection{Estimating the Variance from below}  \label{sec.MinVar} \indf 
   We must now estimate the last term in the above (\ref{f.MajorIntOsc}). We first write a tractable expression of the variance\,: \  setting  \vspace{-2mm} 
\begin{equation}  \label{df.phit} 
\phi_t^\theta := {\ds \int_t^1} \sin\!\left(ws+\theta+\sqrt{\e}\,\o_s\right) ds\, ,  
\end{equation}
we have \quad ${\ds {\rm Var}_{[0,1]}(\phi^\theta) = \,\5 \int_{[0,1]^2} (\phi^\theta_t-\phi^\theta_u)^2\, dt\, du\, = \int_{0\le t < u\le 1} (\phi^\theta_t-\phi^\theta_u)^2\, dt\, du\,}$, \  and then \vspace{-0mm} 
\begin{equation}  \label{f.VarVar} 
{\rm Var}_{[0,1]}\!\left[{ \int_\bullet^1} \sin\!\left(ws\!+\!\theta\!+\!\sqrt{\e}\,\o_s\right) ds\right] = \int_{0\le t < u\le 1}\! \left[\int_t^u \sin\!\left(ws\!+\!\theta\!+\!\sqrt{\e}\,\o_s\right) ds\right]^2 dt\, du\, . 
\end{equation}   

   We shall strongly use this expression (\ref{f.VarVar}) of the variance, and the following remark. 
\brem   \label{rem.modcontpontB} \  {\rm For $\,0\le s<t\le 1\,$, for any $p\ge 2\,$ and for some constant $c_p\,$ we have \parn   
$\E_0^0\Big[\big|\omega_t-\omega_s\big|^{2p}\Big] = \E_0\Big[\big|B_t-B_s - (t-s)B_1\big|^{2p}\Big] = \E\Big[\big|\NN\big(0\,; (t-s)(1-t+s)\big)\big|^{2p}\Big] \le c_p\,|t-s|^p . $   \parsn   
   According to \big([R-Y], I.Theorem (2.1)\big), this entails that setting  \parn   
${\ds S_p(\omega) := \sup_{0\le t<u\le 1}\frac{|\omega_t-\omega_u|}{|t-u|^{\frac{1}{2}-\frac{1}{2p}-\frac{1}{2p^2}}} =\, S_p(-\omega)\, }$ \  defines  a random variable $\,S_p \in L^{2p}\big(\WW_0\,,\P_0^0\big)$, \  such that \pars   
\centerline{$\sup_{t\le s\le u}\limits \big|\omega_s-\omega_t\big| \le\, S_p(\omega)\times |u-t|^{\frac{1}{2}-\frac{1}{2p}-\frac{1}{2p^2}}\,$ \   for any $\,0\le t<u\le 1\,$.} \parn  
}\erem 

   The following crucial estimate amounts in the present context to control the determinant of the Malliavin matrix from below, provided $\,w\not= 0\,$. 
\bpro \label{pro.minorVarianc} \  For any non-null real $\,w\in\R^*$ and any $\,p\ge 2$, using the random variable $\,S_p \in L^{2p}\big(\WW_0\,,\P_0^0\big)$ of Remark \ref{rem.modcontpontB}, \  we have\,:  \vspace{-1mm} 
$$ \inf_{0\le \e\le 1} \inf_{\theta}\,
{\rm Var}_{[0,1]}\!\left[{ \int_\bullet^1} \sin\!\left(ws + \theta + \sqrt{\e}\,\o_s\right) ds \right] \ge \, { \frac{4}{5}} \left[{ \frac{\min\{\pi^2, w^2\}}{2\pi^2}}\right]^{\frac{10p-2-2/p}{p-1- {1}/{p}}}\! \big(|w|+ S_p\big)^{\frac{-8\,p}{p-1- {1}/{p}}} . \vspace{-1mm} $$ 
\epro 
\ub{Proof} \quad  A delicate point here is to obtain the uniformity with respect to $\theta$. \parn 
   By periodicity and symmetry, we can of course restrict to $\,0\le \theta < \pi$. \parn  
Set  \  $\,M_\e^+(\o):= \max_{0\le s\le1}\limits \!\big\{ws+\sqrt{\e}\,\o_s \big\}$,  $\, M_\e^-(\o):= -\min_{0\le s\le1}\limits \!\big\{ws+\sqrt{\e}\,\o_s \big\}\,$  and $\,M_\e:= M_\e^++M_\e^-$. \parn  
Note that \  $M_\e\ge |w|\,$ $\,\P^0_0$-almost surely, for any real $w$.  \parn 
   Thus $\,\P_0^0$-almost surely, the image of $\,[0,1]\ni s\mapsto \big|\sin\!\left(ws + \theta + \sqrt{\e}\,\o_s\right)\!\big|\,$ is the random interval \   
$J_\e^\theta := \big|\sin\!\big([\theta-M_\e^-, \theta+M_\e^+]\big)\big|$, image under $|\sin(\cdot)|$  of the interval $[\theta-M_\e^-, \theta+M_\e^+]$. \parn 
Now, on the one hand \  $\inf_\theta\limits \max J_\e^\theta\,$ is attained for $\,\theta\in \frac{U_\e^--U_\e^+}{2}+\pi\Z\,$ and equal to $\,\tilde M_\e :=\sin\!\big(\frac{\min\{U_\e,\pi\}}{2}\big) \ge \min\{1,M_\e/\pi\} >0\,$, \  and on the other hand, \   the minimal amplitude of $\,J_\e^\theta\,$ is attained for $\,\theta\in \frac{U_\e^--U_\e^++\pi}{2}+\pi\Z\,$ and equal to $\,\inf_\theta\limits {\rm range}(J_\e^\theta) = 2\sin^2\!\big(\frac{\min\{U_\e,\pi\}}{4}\big) \ge \min\{1,M_\e^2/2\pi^2\} >0\,$  \big(we let $\,{\rm range}(J_\e^\theta):= \max J_\e^\theta - \min J_\e^\theta\,$\big).  \parn 
Consider also  \  $  \widetilde{J_\e^\theta} := \max\{\min J_\e^\theta\,, \tilde M_\e/2\} \in J_\e^\theta\,$, \ and 
the random times\,: \parsn 
$(i)$ if $\,\widetilde{J_\e^\theta(\o)}\,$ is hit before $\,\max J_\e^\theta(\o)\,$: \parn 
\centerline{$T_1(\o) := \min\!\Big\{s\in [0,1]\,\Big|\, \big|\sin(ws+\theta+\sqrt{\e}\,\o_s)\big| = \widetilde{J_\e^\theta}(\o)\Big\}\; $;}  \parsn 
\centerline{$T_2(\o) := \min\!\Big\{s\in [T_1(\o),1]\,\Big|\, \big|\sin(ws+\theta+\sqrt{\e}\,\o_s)\big| = \max J_\e^\theta(\o)\Big\}\; $;}  \parn 
$(i)$ if not\,: \parn 
\centerline{$T_1(\o) := \min\!\Big\{s\in [0,1]\,\Big|\, \big|\sin(ws+\theta+\sqrt{\e}\,\o_s)\big| = \max J_\e^\theta(\o)\Big\}\; $;}  \parsn 
\centerline{$T_2(\o) := \min\!\Big\{s\in [T_1(\o),1]\,\Big|\, \big|\sin(ws+\theta+\sqrt{\e}\,\o_s)\big| = \widetilde{J_\e^\theta}(\o)\Big\}\,$.}  \parsn 
Then on the one hand \  we have \quad  $\big|\sin\!\big(ws+\theta+\sqrt{\e}\,\o_s\big)\big| \ge \tilde M_\e/2$ \   for $\,T_1\le s\le T_2\,$, \parn  
and on the other hand, by Remark \ref{rem.modcontpontB} we have\,:   
$$ \5\,\min\{1, M_\e^2/\pi^2\}\, \le \, \min\!\big\{\tilde M_\e/2\, , \inf_\theta\limits {\rm range}(J_\e^\theta)\big\} $$
$$ \le\, \big| \sin\!\big(w\,T_1(\o)+\theta+\sqrt{\e}\,\o_{T_1}\big) - \sin\!\big(w\,T_2(\o)+\theta+\sqrt{\e}\,\o_{T_2}\big)\big| $$
$$ \le\, |w|\, |T_2(\omega)-T_1(\omega)| + \sqrt{\e}\, S_p(\omega)\times \big|T_2(\omega)-T_1(\omega)\big|^{\frac{1}{2}-\frac{1}{2p}-\frac{1}{2p^2}} $$
$$ \le \, \big(|w|+ \sqrt{\e}\, S_p(\omega)\big)\times  \big|T_2(\omega)-T_1(\omega)\big|^{\frac{1}{2}-\frac{1}{2p}-\frac{1}{2p^2}}\, , $$ 
whence $\P_0^0$-almost surely,
$$  T_2(\o) -T_1(\o)\, \ge\, \left(\frac{\min\{1, M_\e^2/\pi^2\}/2}{|w|+ \sqrt{\e}\, S_p(\o)}\right)^{\!\left(\frac{1}{2}-\frac{1}{2p}-\frac{1}{2p^2}\right)\1} =\, \left(\frac{\min\{1, M_\e^2/\pi^2\}/2}{|w|+ \sqrt{\e}\, S_p(\o)}\right)^{\frac{2\,p}{p-1- {1}/{p}}} \, . $$  
Thence we deduce\,:  
$$ \int_{0\le t < u\le 1}\! \left[\int_t^u \sin\!\left(ws\!+\!\theta\!+\!\sqrt{\e}\,\o_s\right) ds\right]^2\! dt\, du\, \ge \int_{T_1\le t < u\le T_2}\! \left|\int_t^u \sin\!\left(ws\!+\!\theta\!+\!\sqrt{\e}\,\o_s\right) ds\right|^2\! dt\, du  $$  
$$ \ge \int_{T_1\le t < u\le T_2}\! \big[(u-t) \,\tilde M_\e/2\big]^2\, dt\, du = \,(\tilde M_\e)^{2} (T_2-T_1)^4\big/ {48} $$
$$ \ge \, {\frac{(\tilde M_\e)^{2}}{48}}\! \left[{ \frac{\min\{\pi^2,\, M_\e^2\}}{2\pi^2}}\right]^{\frac{8\,p}{p-1- \frac{1}{p}}} \big(|w|+ \sqrt{\e}\, S_p\big)^{\frac{-8\,p}{p-1- \frac{1}{p}}} $$
$$ \ge \, { \frac{4}{5}} \left[{ \frac{\min\{\pi^2, M_\e^2\}}{2\pi^2}}\right]^{\frac{10p-2-2/p}{p-1- {1}/{p}}} \big(|w|+ \sqrt{\e}\, S_p\big)^{\frac{-8\,p}{p-1- {1}/{p}}} . $$ 
By (\ref{f.VarVar}), the claim is now clear, using that \  $M_\e\ge |w|>0\,$ $\,\P^0_0$-almost surely.  
\   $\;\diamond$ \pars  

\subsubsection{End of the proof of Proposition \ref{pro.estintdecre}}  \indf 
   The above proposition \ref{pro.minorVarianc} and   (\ref{f.MajorIntOsc}) together provide\,:  \  for any $\,w\in\R^*$,
$$ \sqrt{\sup_{0\le \e\le 1} \E_0^{0}\big[ e^{-\RR_\e\times R}\,\big]}\, \le\, 2 e^{- R^{2/3}}\! +  \sqrt{\E_0^{0}\!\left[ \exp\!\left(\!{ -\,{\ts\frac{2 R^{1/3}}{5C}}\left[{\ts\frac{\min\{\pi^2,\,w^2\}}{2\pi^2}}\right]^{\!\frac{10p-2-\frac{2}{p}}{p-1- \frac{1}{p}}}\! \left(|w|+ S_p\right)^{\frac{- 8\,p}{p-1- \frac{1}{p}}}}\right)\! \right]} .  $$ 
Now, for any positive random variable $V$ and any deterministic $\,Y>0, \, q\in\N^*\,$ we have\,:
$$ \E\big[e^{-Y/V}\big] = \int_0^1 \P\big[e^{-Y/V}>x\big] dx = \int_1^\infty \P\big[V>Y/\log x\big] x\2\, dx =\, Y \int_0^\infty \P\big[tV>1\big] e^{-Yt}\, dt $$
$$ \le \, Y\, \E(V^q) \int_0^\infty t^q\, e^{-Yt}\, dt\, = \, q!\; \E(V^q)\, Y^{-q} . $$ 
In particular we thus have 
$$ \E_0^0\!\left[e^{-Y \left(|w|+ S_p\right)^{\frac{- 8\,p}{p-1- \frac{1}{p}}}}\right] \le\, q!\;  \E_0^0\!\left[ \left(|w|+ S_p\right)^{\frac{8\,p\,q}{p-1- \frac{1}{p}}}\right]\, Y^{-q} , $$ 
and then by the above, setting \   $C(p,q,w):= \left({\ts\frac{5C}{2}}\right)^{q/2}\left[{\ts\frac{2\pi^2}{\min\{\pi^2,\,w^2\}}}\right]^{\!\frac{(5p-1-\frac{1}{p})q}{p-1- \frac{1}{p}}}$ we have\,: 
$$ \sqrt{\sup_{0\le \e\le 1} \E_0^{0}\big[ e^{-\RR_\e\times R}\,\big]}\, \le\, 2 e^{- R^{2/3}}\! + C(p,q,w) \sqrt{q!\;  \E_0^0\!\left[ \left(|w|+ S_p\right)^{\frac{8\,p\,q}{p-1- \frac{1}{p}}}\right]} \, R^{-q/6} .  $$ 
As a consequence, taking for example $\,q=13$ \big(in order to have $\,1-q/6<-1$\big) and $\,p=54$ \Big(in order to have $\,\frac{8\,p\,q}{p-1- \frac{1}{p}} \le 2p\,$; recall Remark \ref{rem.modcontpontB}\,: $S_p \in L^{2p}\big(\WW_0\,,\P_0^0\big)$\Big), we obtain\,: 
$$ \int_1^\infty \sqrt{\sup_{0\le \e\le 1} \E_0^{0}\big[ e^{-\RR_\e\times R}\,\big]}\,  R\,dR\, \le\, 6 + 6\, C(54,13,w) \sqrt{13!\;  \E_0^0\big[ \left(|w|+ S_p\right)^{2p}\big]} \,   < \infty\, .  $$ 
This completes the proof of Proposition \ref{pro.estintdecre}. $\;\diamond$ 
\parb  

  As a consequence, we have the following key approximation result. \vspace{-0mm} 
\bpro  \label{pro.equivIntPe} \  As $\,\e\sea 0\,$, for any $\,w\in\R^*$ and uniformly with respect to $(y,z)\in\R^2$, we have 
$$  \int_{\sR^2} P_\e(\xi',\xi)\, d\xi'd\xi\, = (1+o(1))\!  \int_{\sR^2} \exp\!\left[{\ts {\frac{\ss\rt1}{\sqrt{\e}}\! \left({\xi'} \left[ \frac{\sin w}{w} - y\right] + {\xi} \left[ \frac{1- \cos w}{w} - z\right] \right)}}\right] \EE_0(\xi',\xi)\, d\xi'd\xi\, . \vspace{-1mm} $$ 
\epro 
\ub{Proof} \quad  By (\ref{f.Pe()}) we have  \vspace{-1mm} 
$$ P_\e(\xi',\xi) =\, \exp\!\left[{\ts {\frac{\ss\rt1}{\sqrt{\e}}\! \left({\xi'} \left[ {\ts\frac{\sin w}{w}} - y\right] + {\xi} \left[ {\ts\frac{1-\cos w}{w}} - z\right] \right)}}\right]\! \times \EE_\e(\xi',\xi) \vspace{-0mm}\, . $$ 
By (\ref{f.estL2/L1}) and Proposition \ref{pro.estintdecre} we have  $\,\sup_{0\le \e\le \e_0}\limits\! \big|\EE_\e(\xi',\xi)\big|\! \in\! L^1(\R^2\!, d\xi'd\xi)$, for some $\e_0>0$. \parn  This provides the wanted domination of  $\, P_\e(\xi',\xi) $, which allows  to apply the Lebesgue theorem with respect to $d\xi'd\xi\,$, using the continuity (recall (\ref{df.Eeps}),(\ref{dfAe})) of \mbox{$\, \EE_\e(\xi',\xi)$ at $\e=0\,$.   $\,\diamond$}

\section{Small-time equivalent for $\,p_\e\,$, provided $\,w\not=0$}  \label{sec.asymptpe} \indf 
   The above proposition \ref{pro.equivIntPe} allows to substitute $\,\overline{P}_\e(\xi',\xi) \,$
for $\,P_\e(\xi',\xi)$ in Lemma \ref{lem.calcavpontFb}. \  Doing this and owing to (\ref{df.Eeps}) and (\ref{dfAe}), we directly obtain the following. 
\bPro \label{pro.calcavpontFb0} \   For any $\,w\in\R^*$ and uniformly with respect to $(y,z)\in\R^2$,  as $\,\e\sea 0\,$ we have 
$$ p_\e\big(0\, ; (w, \e\, y\,, \e\, z)\big) =\, \frac{\big(1+o(1)\big)\, e^{- {w^2}/{2 \e}}}{4\pi^2\,\e^{3}\,\sqrt{2\pi\,\e}}  \int_{\sR^2} \overline{P}_\e(\xi',\xi)\, d\xi'd\xi \, , \vspace{-2mm} $$ 
with  \vspace{-1mm} 
$$ \overline{P}_\e(\xi',\xi) = \,  \exp\!\left[{ {\frac{\ss\rt1}{\sqrt{\e}} {\ts \Big( \xi' \left[ \frac{\sin w}{w} - y\right] +  \xi \left[ \frac{1- \cos w}{w} - z\right] \Big)}} }\right]  \times \vspace{-2mm} \hskip30mm  \vspace{-1mm}  $$
$$ \hskip30mm  \times\,\E_0^{0}\bigg[ \exp\!\bigg(\!{\rt1}\! \int_0^{1} \big(\xi\,\cos(ws) - \xi'\sin(ws)\big)\, \o_s\,ds \bigg)\bigg] . $$ 
\ePro

   Now we have the following easy lemma, which merely amounts to a Gaussian computation, or can as well be seen as the definition of the Gaussian (pinned Wiener) measure $\P_0^0$ on the Wiener space $(\WW_0,H_0)$ of the standard Brownian bridge. 
\bLem  \label{lem.Gauss1'} \  For any complex deterministic continuous function $\,u_\cdot$ on $[0,1]$, we have 
$$ \E_0^{0}\!\left[ e^{\int_0^1 u_s\, \omega_s\, ds}\right] = \E_0^{0}\!\left[ e^{\int_0^1 \left(\int_s^1u_\cdot\right)\, d\omega_s}\right] =\, \exp\!\left({\frac{1}{2}} \left[\int_0^1\! \Big(\int_s^1u\Big)^{\!2} ds - \Big(\int_0^1\! \Big(\int_s^1u\Big) ds\Big)^{\!2} \right]\right) . $$
\eLem       

   Therefore we have \vspace{-0mm} 
$$ \E_0^{0}\bigg[ \exp\!\bigg(\!{\rt1}\! \int_0^{1} \big(\xi\,\cos(ws) - \xi'\sin(ws)\big)\, \o_s\,ds \bigg)\bigg] =\, \exp\!\left[-\frac{1}{2}\,(\xi,\xi') (DU^0) \,^t\!(\xi,\xi')\right] . $$ 
Then setting \quad  ${\ds \lambda := \frac{1- \cos w}{w} - z\,}$ \  and  \   ${\ds \mu := \frac{\sin w}{w} - y\,}$, \quad  we have\,: 
$$  \int_{\sR^2} \overline{P}_\e(\xi',\xi)\, d\xi'd\xi \, = \int_{\sR^2} \exp\!\left[-\frac{1}{2}\,(\xi,\xi') (DU^0) \,^t\!(\xi,\xi') +  {\frac{\rt1}{\sqrt{\e}}}\, \big(\lambda\, \xi + \mu\, \xi' \big)\right] d\xi'd\xi \,  $$ 
$$ =\, \frac{2\pi}{\sqrt{\det(U^0)}}\, \exp\!\left( \frac{-1}{2\e}\,(\lambda,\mu) (DU^0)\1 \,^t\!(\lambda,\mu) \right) = \,\frac{2\pi}{\sqrt{\det(U^0)}}\, \exp\!\left(-\,\frac{ \psi(w,y,z) }{4w^2\,\e \,\det(U^0) } \right) ,  $$ 
with
$${\ts  \psi(w,y,z) :=\,2w^2 \left(y-\frac{\sin w}{w}\, \raise1.3pt\hbox{,}\, \frac{1- \cos w}{w} - z \right) \times DU^0\times\, ^t\!\left(y-\frac{\sin w}{w}\, \raise1.3pt\hbox{,}\, \frac{1- \cos w}{w} - z \right) . } \vspace{1mm} $$
Recall also from Remark \ref{rem.detDV0} that 
$$ 0<  \det(U^0) = \frac{1}{4w^4} - \frac{4(1-\cos w)+ \sin^2\!w}{4\,w^6} + \frac{(1-\cos w)\sin w}{w^7}  = \frac{w^2}{8640} + \O(w^4)\,$$  
\hbox{for small $|w|$}. \   With Proposition \ref{pro.calcavpontFb0} and Remark \ref{rem.detDV0}, this yields the the case $\,w\not=0\,$ of Theorem \ref{th.mainres}, and then (by means of banal computations) Remark \ref{rem.coeffwnot=0} as well. Note that the uniformity with respect to $(y,z)\in\R^2$ in Proposition \ref{pro.calcavpontFb0} (already in Proposition \ref{pro.equivIntPe}) indeed allows to replace $(y,z)$ eventually by $(y,z)/\e\,$. 

\section{The singular case $\,w=0$  \if{(1er essai)}\fi }  \label{sec.casew=0}  \indf 
    Let us start again from Lemma \ref{lem.calcavpontFb}\,:  \vspace{-1mm} 
\begin{equation*}
p_\e\big(0\, ; (0, \e\, y\,, \e \,z)\big) =\, \frac{1}{4\pi^2\,\e^{3}\,\sqrt{2\pi\,\e}}\,  \int_{\sR^2} P^0_\e(\xi',\xi)\, d\xi'd\xi\,,  \vspace{-1mm}
\end{equation*}  
with \quad \vspace{-2mm}
$$ P^0_\e(\xi',\xi) \, :=  \,\E_0^{0}\bigg[\! \exp\!\bigg(\!{\ts\frac{\rt1}{\sqrt{\e}}}\! \bigg({\xi'} \left[ \int_0^{1}\! \cos\!\big(\sqrt{\e}\,\o_{s}\big) ds - y\right] + {\xi} \left[\int_0^{1}\! \sin\!\big(\sqrt{\e}\,\o_{s}\big) ds - z\right]\!\bigg) \!\bigg) \bigg] = \vspace{-0mm}  $$  
{\small 
$$ \exp\!\left[{ {\frac{\ss\rt1}{\sqrt{\e}}\! \left({\xi'} \big[1 - y\right] - {\xi} \, z\big)}}\right]\! 
\times\, \E_0^{0}\bigg[\! \exp\!\bigg(\rt1\!\! \int_0^{1}\!\int_0^{1}\! \big[ \xi\, \cos\!\big(\sqrt{\e}\,\tau\, \o_s\big) - {\xi'} \sin\!\big(\sqrt{\e}\,\tau\, \o_s\big) \big]\, \o_s\, ds\, d\tau \!\bigg) \bigg] . \vspace{1mm}  $$  } 
Let us change $\,\xi'$  into $\,\xi'/\sqrt{\e}\,$ in (\ref{f.expr7}), so that we instead  have\,: \vspace{1mm} 
$$ 4\pi^2\sqrt{2\pi}\, \e^{4}\times \, p_\e\big(0\, ; (0,\e\, y\,, \e\, z)\big)  \hskip1mm \vspace{-1mm}  $$ 
{\small 
$$ =  \int_{\sR^2} \E_0^{0}\bigg[\! \exp\!\bigg(\!{\frac{\rt1}{\sqrt{\e}}} \bigg({\frac{\xi'}{\sqrt{\e}}} \left[ \int_0^{1}\! \cos\!\big(\sqrt{\e}\,\o_{s}\big) ds - y\right] + {\xi} \left[\int_0^{1}\! \sin\!\big(\sqrt{\e}\,\o_{s}\big) ds - z\right]\! \bigg) \!\bigg) \bigg] d\xi'd\xi\,   $$  }
\begin{equation}  \label{f.expr70'} 
= \int_{\sR^2} \e^{ \rt1\! \left({\frac{\xi'}{\e}} [1 - y] - {\frac{\xi}{\sqrt{\e}}}\, z \right)} \times\, \E_0^{0}\bigg[\! \exp\!\bigg(\rt1\!\! \int_0^{1}\! \left[ \xi\, {\ts\frac{\sin\left(\sqrt{\e}\, \o_s\right)}{\sqrt{\e}}} - {\xi'}\, {\ts\frac{1- \cos\left(\sqrt{\e}\, \o_s\right)}{{\e}}} \right] ds \!\bigg) \bigg] d\xi'd\xi\,  . \vspace{1mm}  
\end{equation}   
The phase is now \  $\,{\ds \tilde A_\e(\xi',\xi) := \xi'\,\tilde U_1^\e + \xi\, \tilde U_2^\e }\,$, \quad 
with the phase vector\,: 
\begin{equation} \label{fde.Ueo} 
\tilde U^\e(\o) := \left(-\int_{0}^1 \frac{1- \cos\left(\sqrt{\e}\, \o_s\right)}{{\e}}\, ds \, , \int_{0}^1 \frac{\sin\!\big(\sqrt{\e}\,\o_{s}\big)}{\sqrt{\e}}\, ds  \right) \! .  \vspace{-1mm}
\end{equation}   
The above replaces (\ref{f.Pe()}), (\ref{dfAe}) and (\ref{fde.Ve}). \   We have again to dominate the last term of (\ref{f.expr70'})\,:  
\begin{equation} \label{df.Etilde} 
\tilde \EE_\e(\xi',\xi):= \E_0^0\!\left[e^{\rt1 (\xi'\,\tilde U_1^\e + \xi\, \tilde U_2^\e)}\right]
\end{equation}
(analogue  to (\ref{df.Eeps})) by an integrable $\e$-independent expression. \  

\subsection{Analysis of the phase vector $\,\tilde U^\e$ and second use of [T2]} \label{sec.Analw=0} 
   We here perform the analogue of Sections \ref{handlerrort} and \ref{sec.IPPT2} for the singular case $w=0\,$.  We successively have\,:  
\begin{equation} \label{fde.Ueoo} 
\tilde U^0(\o) = \left(-{\frac{1}{2}}\int_{0}^1\! \o_{s}^2 \, ds \, , \int_{0}^1\! \o_s  \, ds  \right) \! ,  \vspace{1mm}
\end{equation}   
$$ \frac{d}{dt} \nabla \tilde U_2^{\e}(\o)(t) = \int_t^1\!\cos\!\big(\sqrt{\e}\,\o_{s}\big) ds - \int_0^1\! du\!\int_u^1\!\cos\!\big(\sqrt{\e}\,\o_{s}\big) ds \; ; $$  
$$  \frac{d}{dt} \nabla \tilde U_1^{\e}(\o)(t) = \int_0^1\! du\!\int_u^1\frac{\sin\!\big(\sqrt{\e}\,\o_{s}\big)}{\sqrt{\e}}\, ds - \int_t^1\frac{\sin\!\big(\sqrt{\e}\,\o_{s}\big)}{\sqrt{\e}}\, ds \; ; \vspace{-0mm}  $$ 
$$ \LL \tilde U^{\e}_1(\o) =\, -\, \e^{-1/2} \int_0^1\!\sin\!\big(\sqrt{\e}\,\o_{s}\big) \,\o_s \, ds\, ,  \vspace{-0mm}  $$ 
$$  \frac{d}{dt} \nabla \LL \tilde U_1^{\e}(\o)(t) = \int_0^t\! \bigg[\omega_s \cos(\sqrt{\e}\,\o_{s}) + \frac{\sin(\sqrt{\e}\,\o_{s})}{\sqrt{\e}}\bigg] ds - \int_0^1\!\int_0^t\! \bigg[\omega_s \cos(\sqrt{\e}\,\o_{s}) + \frac{\sin(\sqrt{\e}\,\o_{s})}{\sqrt{\e}}\bigg] ds\,  ; \vspace{-0mm}  $$ 
and \vspace{-2mm} 
$$ D\tilde U^{0}(\o) :=  \begin{pmatrix} \int_0^1\!\left[\int_{t}^1 \o\right]^2\! dt - \left[\int_0^1\! dt\int_{t}^1\! \o\right]^2 & \int_0^1 (t-\5)\!\left[\int_{t}^1\! \o\right] \!dt \cr \int_0^1 (t-\5)\!\left[\int_{t}^1\! \o\right] \!dt &  {1}/{12} \end{pmatrix}\! . \vspace{0mm} $$ 
Thus the estimates of Lemmas \ref{lem.estimA} and \ref{lem.estimLLUe} relating to $\|\nabla^k\tilde U_2^\e\|_{H_0^{\otimes k}}\,$  and $ \|\nabla^k\LL \tilde U_2^\e\|_{H_0^{\otimes k}}\,$ remain the same as before (regarding $\,U_2^\e$). Whereas for $\|\nabla^k\tilde U_1^\e\|_{H_0^{\otimes k}}\,$  and $ \|\nabla^k\LL \tilde U_1^\e\|_{H_0^{\otimes k}}\,$, a straightforward adaptation provides the following. \vspace{-1mm} 
\blem \label{lem.estimLLUetilde} \  We almost surely have\,: \vspace{-1.5mm}  
$$ \|\nabla^k\tilde U_1^\e\|_{H_0^{\otimes k}}\, \le\, 4\sqrt{2}\, k\1\, (2\sqrt{\e}\,)^{k-2}  \  ; \quad \|\nabla^k\LL \tilde U_1^\e\|_{H_0^{\otimes k}}\, \le\, 4\sqrt{2}\;  (2\sqrt{\e}\,)^{k-2}\times \left[ 2 + \e\!\int_0^1\!\o^2 \right]^{1/2} , \vspace{-2mm}    $$    
for $\,k\ge 2\,$, \  and for $\,k=1\,$:  \vspace{-2mm}
$$ \|\nabla\tilde U_1^\e\|_{H_0}\, \le \int_0^1\!|\o|  \  ; \quad \|\nabla\LL \tilde U_1^\e\|_{H_0}\, \le\, \sqrt{2}\;  \left[ \int_0^1\!\o^2 \right]^{1/2} . \vspace{-2mm}   $$   
Thus \vspace{-1mm} 
$$ \|\nabla^k\tilde A_\e(\xi',\xi)\|_{H_0^{\otimes k}}\, \le \,6\, \sqrt{\xi^2+\xi'\,\!^2}\; (2\sqrt{\e}\,)^{k-2}\times \left[ 2 + \e\!\int_0^1\!\o^2 \right]^{1/2} ,  \vspace{-1mm} $$     
for $\,k\ge 2\,$, \  and for $\,k=1\,$:  \vspace{-2mm}
$$ \|\nabla\tilde A_\e(\xi',\xi)\|_{H_0}\, \le \,\sqrt{3}\,\sqrt{\xi^2+\xi'\,\!^2}\;  \left[ \int_0^1\!\o^2 \right]^{1/2}\, . \vspace{-2.5mm} $$     
\elem  
  As a consequence, we directly deduce the following analogue of Lemma \ref{lem.estimMu}. 
\blem \label{lem.estimMu0} \  For any positive $\,r,\e\,$ we almost surely have\,: \vspace{-0mm}  
$$  M_1[\tilde U_1^{\e},r] +  M_1[\tilde U_2^{\e},r] < 3+ 8\,r^2(1+\e)\, e^{4r\,\e}+ {\ts\int_0^1\!\o^2} \ \hbox{ ; } \  M_2[\tilde U_1^{\e},r] +  M_2[\tilde U_2^{\e},r] < 8\,r^2(1+\e)\, e^{4r\,\e}\, ; $$ 
$$  M_1[\LL \tilde U_1^{\e},r] +  M_1[\LL \tilde U_2^{\e},r] < 8\,r\,(2r+1)\, e^{4r\,\e}\!\left[ 2 + \e\!\int_0^1\!\o^2 \right] + 2\int_0^1\!\o^2 \;  ; $$
$$ M_2[\LL \tilde U_1^{\e},r] +  M_2[\LL \tilde U_2^{\e},r] < (16+ 8\e)\,r^2\, e^{4r\,\e}\left[ 2 + \e\!\int_0^1\!\o^2 \right] . \vspace{-2mm} $$ 
\elem  
Then this entails\,: \vspace{-1mm} 
$$  L_2[\tilde U^{\e},r] \,\le\, 8\,r^2(1+\e)\, e^{4r\,\e}+ (16+ 8\e)\,r^2\, e^{4r\,\e}\left[ 2 + \e\!\int_0^1\!\o^2 \right] + 1 + \int_0^1\!\o^2 , \vspace{-2mm} $$ 
and (to replace (\ref{f.exprL1Ve1*}))\,: \vspace{-2mm} 
$$ 1\le  L_1[\tilde U^{\e},r]\, \le\, 4+ 8\,r^2(1+\e)\, e^{4r\,\e} + 8\,r\,(2r+1)\, e^{4r\,\e}\!\left[ 2 + \e\!\int_0^1\!\o^2 \right] + 4\int_0^1\!\o^2 \vspace{-1mm}  $$ 
\begin{equation}  \label{f.estOmegt} 
\le\, \O(1)+ \big(4+\O(\e)\big) (\o^*)^2 =: \tilde \Omega_\e\,. \vspace{1mm}
\end{equation}

   Then we proceed as in Section \ref{sec.IPPT2}, in order to apply (2.1) in Theorem 2.1 of [T2] again, of course with $\tilde U^{\e}\equiv \tilde q\, $ instead of $U^{\e}$. As in (\ref{f.estgradq2}) and for the same reason, we again have  \vspace{-1mm} 
$$ \min\!\Big\{ \big\| \nabla \tilde q^{(z)} \big\|^2_{H_0} \Big|\, |z|=1\Big\}  \ge \,  \frac{\det(D\tilde U^{\e})}{  \Tr(D\tilde U^{\e}) } \,\raise1.9pt\hbox{,} \vspace{-1mm} $$ 
so that for some $\,C>0\,$, any $\,r\in\,]1,e[\,$ and any $\,\xi^2+\xi'\!\,^2 >1\,$, we again have (\ref{f.estT2})\,:
{\small 
\begin{equation*} \label{f.estT2o} 
\left|\E_0^0\!\left[e^{\rt1 (\xi'\,\tilde U_1^\e + \xi\, \tilde U_2^\e)}\right]\right|^{2} \le C\; \E_0^{0}\bigg[\exp\!\bigg(\! \frac{L_2[\tilde U^{\e}\!,r]}{L_1[\tilde U^{\e}\!,1]^3}\!\bigg) \bigg] \times \E_0^{0}\bigg[\!\exp\!\bigg(\! -  \frac{\det(D\tilde U^{\e})/\Tr(D\tilde U^{\e})}{147\,C\,L_1[\tilde U^{\e}\!,1]} \sqrt{\xi^2+\xi'\!\,^2}\,\bigg)\bigg]  , \vspace{-1mm} 
\end{equation*}  } \noindent 
\hskip-1mm and then \vspace{-1mm} 
$$ \E_0^{0}\bigg[\exp\!\bigg(\! \frac{L_2[\tilde U^{\e}\!,r]}{L_1[\tilde U^{\e}\!,1]^3}\!\bigg) \bigg] \le \E_0^{0}\Big[e^{L_2[\tilde U^{\e}\!,r]} \Big] \le \sqrt{C}\; \E_0^{0}\Big[e^{\left(1+\e\log C\right) \int_0^1\!\o^2} \Big] = 
\sqrt{ \ts\frac{ \sqrt{2\left(1+\e\log C\right)}  } { {\ts \sin}\left(\!\sqrt{2\left(1+\e\log C\right)}\right)} } \le\, 2 $$ 
for small enough $\,\e\,$. \  See Corollary \ref {cor.calculQuadrF} below for the Laplace transform of the law of $\int_0^1\o^2$ under $\P^0_0\,$. \  Therefore the analogue of (\ref{f.estL2/L1}) holds\,: 
$$  \left|\E_0^0\!\left[e^{\rt1 (\xi'\,\tilde U_1^\e + \xi\, \tilde U_2^\e)}\right]\right|^{2}  \le\, 2C\; \E_0^{0}\!\left[\exp\!\bigg(\! -  \frac{\det(D\tilde U^{\e})/\Tr(D\tilde U^{\e})}{147\,C\,L_1[\tilde U^{\e}\!,1]} \sqrt{\xi^2+\xi'\!\,^2}\,\bigg)\right] , $$ 
and then, owing to (\ref{df.Etilde}) and (\ref{f.estOmegt})\,: \vspace{-2mm} 
\begin{equation} \label{f.applT20}
\left|\tilde \EE_\e(\xi',\xi)\right|^{2} \le\, 2C\; \E_0^{0}\!\left[\exp\!\bigg(\! -  \frac{\det(D\tilde U^{\e})}{\Tr(D\tilde U^{\e})\, \tilde \Omega^{\e}}\times  \frac{\sqrt{\xi^2+\xi'\!\,^2}}{147\,C}\,\bigg)\right] . \vspace{-2mm} 
\end{equation}

\subsection{Domination of the decreasing integral in the singular case} \label{sec.Analogw=0}   \indf 
Here we perform the analogue of Section \ref{sec.changbsuit1sapple[T2]} for the singular case $w=0\,$. 
    We have to estimate the new key variance (determinant of the Malliavin matrix) from below, and to dominate the related expected value. Proceeding as for Proposition \ref{pro.estintdecre}, let us denote by $\,\tilde \lambda_{\e}\,$ the (almost surely positive) lowest eigenvalue of $\,D\tilde U^\e$. Then on the one hand we have  
the following analogue of (\ref{f.estimReVar})(\ref{f.MajorIntOsc}), for any deterministic positive $R\,$:  
$$ \E_0^{0}\!\left[ e^{\ts -\frac{\det(D\tilde U^{\e})}{\Tr(D\tilde U^{\e})\, \tilde \Omega_\e}\times R}\right] \le \E_0^{0}\!\left[ \exp\!\left(-\, \frac{\tilde \lambda_{\e} R}{2\left[\O(1)+ \big(4+\O(\e)\big) (\o^*)^2\right]}\right)\right] $$ 
$$ \le \, \P_0^0\big(\o^*> R^{1/3}/C\big) + \E_0^{0}\!\left[ \exp\!\left(\frac{-\, \tilde \lambda_{\e}\, R}{\O(1)+ (4+\O(\e)) C\2 R^{2/3}}\right) 1_{\{\o^*\le R^{1/3}/C\}} \right] $$
\begin{equation}  \label{f.estimw=0} 
\le \, 2\, e^{-2C\2 R^{2/3}} +\, \E_0^{0}\!\left[ e^{-\, \tilde \lambda_{\e}\, R^{1/3}}\,1_{\{\o^*\le R^{1/3}\}} \right] \quad \hbox{(for $0\le \e\le 1$ and large enough $R$)}. 
\end{equation}
On the other hand, we have\,:     \vspace{-2mm} 
$$ \langle v, D\tilde U^\e v\rangle =  {\rm Var}_{[0,1]}\!\left[\int_\bullet^1\!\Big[ v_1\,{\ts\frac{\sin\left(\sqrt{\e}\,\o_{s}\right)}{\sqrt{\e}}} +v_2 \cos\!\big(\sqrt{\e}\,\o_{s}\big)\Big] ds\right] \quad \hbox{for any }\, v\in\S^1 ,  \vspace{-2mm}  $$
so that, setting \quad 
\begin{equation}  \label{df.Vartheps} 
V^\theta_\e(\o) := 
{\rm Var}_{[0,1]}\!\left[ \sin\theta  \int_\bullet^1\! \cos\!\left(\sqrt{\e}\,\o_\cdot\right) +  \cos\theta \int_\bullet^1 \frac{\sin\!\left(\sqrt{\e}\,\o_\cdot\right)}{\sqrt{\e}} \right]  
\end{equation} 
$$ = \int_{0<t<u<1}\left[ \sin\theta  \int_t^u\! \cos\!\left(\sqrt{\e}\,\o_\cdot\right) +  \cos\theta \int_t^u \frac{\sin\!\left(\sqrt{\e}\,\o_\cdot\right)}{\sqrt{\e}} \right]^2 dt\, du \qquad \hbox{according to (\ref{f.VarVar}),} \vspace{1mm}  $$ 
we have\,:  \vspace{-1mm} 
\begin{equation}  \label{f.VarVar0} 
\tilde \lambda_{\e}\, =\, \inf_\theta\,  V^\theta_\e\, . 
\end{equation}   

   According to (\ref{f.applT20}), our aim is to establish that    
\begin{equation}  \label{f.butdomw=0} 
\int_1^\infty \sup_{0\le \e\le 1}\limits\sqrt{ \E_0^{0}\!\left[ e^{-\,\tilde \lambda_{\e}\, R^{1/3}}\,1_{\{\o^*\le R^{1/3}\}} \right]}\, R\, dR\, <\infty\,. \vspace{-3mm} 
\end{equation} 

\subsubsection{Getting the extremum over $\theta$ out of the expectation $\E_0^{0}$}  \label{sec.AstH} \indf 
    We must now obtain the analogue of Proposition \ref{pro.minorVarianc} relating to Variance $\,V^\theta_\e\,$. \parn
Note that by symmetry it is enough to consider $\theta\in [0,\pi[\,$, and that using (\ref{f.VarVar}) we easily have\,: for any $\theta, \theta', \e, \o$,  
\begin{equation}  \label{f.Varth-Varth'} 
\big| V^\theta_\e(\o) - V^{\theta'}_\e(\o)\big|\, \le |\theta-\theta'|\times \frac{1+(\o^*)^2}{6}\, \raise1.5pt\hbox{.} 
\end{equation}   

    We now use a nice idea of ([H], Lemma 4.7), which consists in getting an extremum out of a probability, by means of a simple counting argument. Namely, for any positive $t\,$ consider the angles $\,k\,t\,$, for $\N\ni k<\pi/t\,$. Since by (\ref{f.Varth-Varth'}), restricting on the event $\{\o^*\le R^{1/3}\}$ we have \  $\big| V^\theta_\e - V^{kt}_\e\big|\, \le |\theta-kt|\,\frac{1+R^{2/3}}{6} \le |\theta-kt|\,\frac{R^{2/3}}{4}$ (for $R\ge3$)\raise0pt\hbox{,} \  we obtain   \parn   
$\inf_\theta\limits\, V^\theta_\e \ge \inf_k\limits\, V^{kt}_\e - \,t\,\frac{R^{2/3}}{4}\,$\raise0.6pt\hbox{,} \   and then, on the event $\{\o^*\le R^{1/3}\}\,$:  \vspace{-2mm} 
$$ \left\{\inf_\theta\limits\, V^\theta_\e < R^{-1/3} s\right\}\, \subset\, \left\{\inf_k\limits\, V^{kt}_\e < R^{-1/3} s+ {\ts\frac{R^{2/3}}{4}}\,t \right\}\, \subset \bigcup_{0\le k <\pi/t}\, \left\{V^{kt}_\e < R^{-1/3} s + {\ts\frac{R^{2/3}}{4}} \,t\right\} . \vspace{-2mm} $$ 
Choosing \  $t= \pi \,s/R\,$, we thus obtain\,: 
$$ \P_0^0\!\left(\tilde \lambda_{\e} < R^{-1/3}\, s\, ,\,\o^*\le R^{1/3}\right)  \le\, \frac{R}{s}\, \sup_\theta\limits\, \P_0^0\!\left(V^{\theta}_\e < \, 2\,R^{-1/3}\, s \right) . $$ 
Hence, 
$$ \E_0^{0}\!\left[ e^{-\, \tilde \lambda_{\e}\, R^{1/3}}\,1_{\{\o^*\le R^{1/3}\}} \right] = \int_0^\infty  \P_0^0\!\left(\tilde \lambda_{\e} < R^{-1/3}\, s\, ,\,\o^*\le R^{1/3}\right)  e^{-s}\, ds $$
$$ \le\, R \int_0^\infty  \sup_\theta\limits\, \P_0^0\!\left(V^{\theta}_\e < \, 2\,R^{-1/3}\, s \right)  e^{-s}\, \frac{ds}{s}\, =\, R \int_0^\infty  \sup_\theta\limits\, \P_0^0\!\left(V^{\theta}_\e < \, 2\,s \right)  e^{-R^{1/3}\, s}\, \frac{ds}{s}\,\raise1.9pt\hbox{,}  $$ 
whence  \vspace{-1mm} 
$$ \E_0^{0}\!\left[ e^{-\, \tilde \lambda_{\e}\, R^{1/3}}\,1_{\{\o^*\le R^{1/3}\}} \right]  \le\, R^{2/3}\, e^{-R^{1/3}} + R \int_0^1 \sup_\theta\limits\, \P_0^0\!\left(V^{\theta}_\e < \, 2\,s \right)  e^{-R^{1/3}\, s}\,\frac{ds}{s}\,  
\raise1.9pt\hbox{.} $$ 
As a consequence, the estimate (\ref{f.butdomw=0}) will follow now from 
$$\int_1^\infty \sup_{0\le \e\le 1}\limits\sqrt{\int_0^1 \sup_\theta\limits\, \P_0^0\!\left(V^{\theta}_\e < \, 2\,s \right)  e^{-R^{1/3}\, s}\,\frac{ds}{s}}  \; R^{3/2}\, dR\, <\infty\,, $$
or equivalently, 
\begin{equation}  \label{f.estimlambdt} 
\int_1^\infty \sup_{0\le \e\le 1}\limits\sqrt{\int_0^1 \sup_\theta\limits\, \P_0^0\!\left(V^{\theta}_\e < \, 2\,s \right)  e^{-R\, s}\,\frac{ds}{s}}  \; R^{13/2}\, dR\, <\infty\,. \vspace{-2mm} 
\end{equation}

\subsubsection{Estimating the Variance from below for $w=0$}  \label{sec.MinVar0} \indf 
   We here establish the following substitute for Proposition \ref{pro.minorVarianc}, with the big difference that in the present degenerate case $\,w=0$, an uniform estimate of the variance from below can no longer be sufficient. Thus a somewhat more subtle proof is here necessary. 
\bpro \label{pro.MinVar-w=0} \  The following estimate holds\,: \  for any $\,x\in [0,1]$ and $\,p\ge 2\,$, there exists an explicit positive finite constant $C_p$ such that we have\,:  \vspace{-1mm} 
$$ \sup_{\theta}\, \P_0^0\Big(\inf_{0\le \e\le 1}\limits V^\theta_\e <x\Big) \le  C_p\, \E_0^0\big(S_{p}^{2p}\big)\, x^{(2p-3- 3/p)/4} + \min\!\left\{2\,\raise1pt\hbox{,}\, \frac{4\sqrt{\pi}\, x^{-1/8}\, e^{-\, \frac{\pi}{196}\, x^{-1/4}}}{7 \left(1-e^{-\,\frac{2\pi}{49}\, x^{-1/4}}\right)}\right\} \raise1.6pt\hbox{.} \vspace{-2mm} $$ 
\epro
\ub{Proof} \quad It is enough to consider $\,0\le \theta<\pi$. \  Again let \  $\o^* := \max_{0\le s\le 1}\limits |\o_s|\,$.
\parsn 
$(i)$ \  Consider first the stopping time  \  $T:= \min\big\{1,\inf\!\big\{ s>0\,\big| \, |\o_s| = \5 |\sin\theta|\big\}\big\}$.   \  Then on the one hand for $\,0\le s\le T\,$ we have \vspace{-1mm} 
$$ \sin\theta  \cos\!\left(\sqrt{\e}\,\o_s\right) + \cos\theta\, \frac{\sin\!\left(\sqrt{\e}\,\o_s\right)}{\sqrt{\e}}\, \ge\, \5 |\sin\theta|\, \cos\!\left(\5\sqrt{\e}\,|\sin\theta|\right) - \5 \cos\theta\,|\sin\theta|\, \ge\, \5 |\sin\theta|\,, \vspace{-1mm}  $$ 
so that by (\ref{df.Vartheps}) we have\,: \vspace{-2mm} 
$$ V^\theta_\e(\o) \ge  \int_{0<t<u<T}\left[ \int_t^u\! \5 |\sin\theta|\, ds\right]^2 dt\, du\, =\, \frac{T^2\, \sin^2\!\theta}{48}\, \raise1.9pt\hbox{,} $$ 
and on the other hand, by Remark \ref{rem.modcontpontB} we have\,:   \  either $\,T=1$ (on $\{\o^*<\5 |\sin\theta\,|\}$) \   or 
$$ \5 |\sin\theta\,|\,=  | \o_T-\o_0| \le\, S_p(\omega)\times \big|T(\omega)\big|^{\frac{1}{2}-\frac{1}{2p}-\frac{1}{2p^2}}\, , $$ 
so that $\,\P_0^0$-almost surely \vspace{-3mm} 
$$ V^\theta_\e\,\ge\, \frac{|\sin\theta\,|^\frac{4p-2-2/p}{p-1- {1}/{p}}}{48}\, (2\,S_p)^{\frac{-2\,p}{p-1- {1}/{p}}} \, . $$ 
$(ii)$ \  Suppose then that $\,\sin\theta < {1}/{9}\,$\raise0pt\hbox{.} We distinguish two sub-cases, as follows. \parsn
$(ii1)$ \  Consider first the case \   $\o^*> \sqrt{1/\e}\,$, \  and thereon the stopping times\,: \parn
$$ T_1:= \inf\!\big\{ s>0\,\big| \, |\sin\!\left(\sqrt{\e}\,\o_s\right)\!| = \sin1\big\}  \quad \hbox{and} \quad T_2:= \inf\!\big\{ s>T_1\,\big| \, |\sin\!\left(\sqrt{\e}\,\o_s\right)\!| = \5\,\sin1\big\}\, . $$ 
Then for any $\,s\in[T_1,T_2]\,$ we have \  
${\ds \left|\frac{\sin\!\left(\sqrt{\e}\,\o_s\right)}{\sqrt{\e}}\right| \ge\, \frac{\sin1}{2\sqrt{\e}} \,\raise1.6pt\hbox{,} }$ \quad  which entails 
$$ \left|\, \cos\theta\, \frac{\sin\!\left(\sqrt{\e}\,\o_s\right)}{\sqrt{\e}} \pm \sin\theta  \cos\!\left(\sqrt{\e}\,\o_s\right) \right| \ge\, \sqrt{80/81}\,\, \frac{\sin1}{2\sqrt{\e}}- \frac{1}{9} \,>\, \frac{2}{7\sqrt{\e}}\,\raise1.5pt\hbox{,} $$ 
so that 
$$ V^\theta_\e\,\ge \int_{T_1<t<u<T_2}\big[2(u-t)/7\big]^2 dt\, du\, =\, \frac{(T_2-T_1)^2}{147}\, \raise1.9pt\hbox{.} $$ 
Since moreover (as previously) \vspace{-1.5mm} 
$$ {\ts\frac{1}{2\sqrt{2}}} < \5\,\sin1\, = \big| \o_{T_2} - \o_{T_1}\big| \le \, S_p(\omega)\times \big|T_2(\omega)-T_1(\omega)\big|^{\frac{1}{2}-\frac{1}{2p}-\frac{1}{2p^2}} \, , $$ 
we obtain\,: $\,\P_0^0$-almost surely, for any $\,p\ge 2\,$: \vspace{-2.5mm} 
$$ V^\theta_\e\, \ge\, \frac{1}{150} \left(2\sqrt{2}\, S_p \right)^{\frac{-2\,p}{p-1- {1}/{p}}}\, 1_{\{9\,|\sin\theta\,| < 1<\sqrt{\e}\,\o^*\}} \ge\, \frac{1}{48} \left(5\, S_p \right)^{\frac{-2\,p}{p-1- {1}/{p}}}\, 1_{\{9\,|\sin\theta\,| < 1<\sqrt{\e}\,\o^*\}}\, . $$ 
$(ii2)$ \  Consider then the case \   $7\,|\sin\theta\,| \le \o^*\le \sqrt{1/\e}\,$, \  and thereon the stopping times\,: \parn
$$ T_1:= \inf\!\big\{ s>0\,\big| \, |\o_s| = \o^*\big\}  \quad \hbox{and} \quad T_2:= \inf\!\big\{ s>T_1\,\big| \, |\o_s| = 4\,\o^*/5\big\}\, . $$ 
Then for any $\,s\in[T_1,T_2]\,$ we have\,: 
$$ \left|\frac{\sin\!\left(\sqrt{\e}\,\o_s\right)}{\sqrt{\e}}\right| \ge\, |\o_s|\times \left(1-\frac{\e\, (\o^*)^2}{6}\right) \ge\, \frac{2}{3}\,\o^*\, \raise0pt\hbox{,} \vspace{-2mm}  $$ 
whence \vspace{-1mm} 
$$ \left|\, \cos\theta\, \frac{\sin\!\left(\sqrt{\e}\,\o_s\right)}{\sqrt{\e}} \pm \sin\theta  \cos\!\left(\sqrt{\e}\,\o_s\right) \right| \ge\, \sqrt{80/81}\, \frac{2}{3}\, \o^* - |\sin\theta\,|\, > \,\o^*/2\, , $$ 
so that 
$$ V^\theta_\e\, \ge \int_{T_1<t<u<T_2}\big[(u-t)\, \o^*/2\big]^2 dt\, du\, =\, \frac{(T_2-T_1)^2}{48}\,(\o^*)^2\,  \raise0pt\hbox{.}  $$ 
Since moreover (as previously) \  
${\ds \frac{\o^*}{5} = \big| \o_{T_2} - \o_{T_1}\big| \le \, S_p(\omega)\!\times\! \big|T_2(\omega)-T_1(\omega)\big|^{\frac{1}{2}-\frac{1}{2p}-\frac{1}{2p^2}} ,}$ \  
we obtain\,: $\,\P_0^0$-almost surely, \vspace{-2mm} 
$$ V^\theta_\e\,\ge\, \frac{(\o^*)^\frac{4p-2-2/p}{p-1- {1}/{p}}}{48}\, (5\,S_p)^{\frac{-2\,p}{p-1- {1}/{p}}}\, 1_{\big\{|\sin\theta\,| < 1/9\,,\,7\,|\sin\theta\,| < \o^*\le \sqrt{1/\e}\,\big\}} \, . $$ 
$(iii)$ \  Summing up the above cases $(i), (ii1)$ and $(ii2)$, so far we have\,:  \  for any $\,p\ge 2$, 
$$ \inf_{0\le \e\le 1}\limits\, V^\theta_\e\,\ge\, \frac{1}{48}\,\Big(\max\!\big\{|\sin\theta\,|\,,\, \min\{\o^*,1\}\!\times\! 1_{\{7\,|\sin\theta\,| \le \o^*\}}\big\}\Big)^\frac{4p-2-2/p}{p-1- {1}/{p}}\times  (5\,S_p)^{\frac{-2\,p}{p-1- {1}/{p}}} \, . $$
This entails\,: \  for any $\,x\in [0,1]$ and any $\,p\ge 2$, 
$$ \P_0^0\Big(\inf_{0\le \e\le 1}\limits\, V^\theta_\e <x\Big) \le\,  \P_0^0\Big(\o^* < 7\, |\sin\theta|\,, (\sin^2\!\theta)^{2p-1-1/p} < (5\,S_p)^{2 p}\, (48\, x)^{p-1- {1}/{p}}\Big) \hskip12mm \vspace{-2mm}  $$ 
$$ \hskip32mm + \,  \P_0^0\Big(\o^* \ge 7\,|\sin\theta|\,, \big(\min\{\o^*,1\}\big)^{4p-2-2/p} < (5\,S_p)^{2 p}\, (48\, x)^{p-1- {1}/{p}}\Big)   $$ 
$$ \le 1_{\{|\sin\theta|> x^{1/8}\}}\, \P_0^0\Big(x^{(2p-1-1/p)/4} < (5\,S_p)^{2 p}\, (48\, x)^{p-1- {1}/{p}}\Big) + 2\times  1_{\{|\sin\theta|\le x^{1/8}\}}\, \P_0^0\big(\o^* < 7\,x^{1/8}\big) $$ 
$$ +\, \P_0^0\Big( \big(x^{1/8}\big)^{4p-2-2/p} < (5\,S_p)^{2 p}\, (48\, x)^{p-1- {1}/{p}}\Big) , $$ 
so that for any $\,p\ge 2\,$ we obtain\,: 
\begin{equation}  \label{f.estimsupPV<x0} 
\P_0^0\Big(\inf_{0\le \e\le 1}\limits\, V^\theta_\e <x\Big)\le \, 2\, (48)^{p-1- {1}/{p}}\, \E_0^0\big((5\,S_p)^{2 p}\big)\, x^{(2p-3- 3/p)/4} + 2\, \P_0^0\big(\o^* < 7\,x^{1/8}\big) \, . 
\end{equation} 
$(iv)$ \  Now, according to \big([B-O] (4.12)\big) and to the Poisson formula \big(see for example (38) in [B-P-Y]\big), for any positive $y\,$ we have  \vspace{-2mm} 
$$ \P_0^0(\o^* < y) = \sum_{n\in\sZ} (-1)^n\, e^{-2n^2 y^2} = \frac{2\sqrt{\pi}}{y}\, \sum_{n\in \sN}\,  \exp\!\left({-\,\frac{(2n+1)^2\,\pi^2}{4\, y^2}}\right) ,  \vspace{-2mm} $$ 
or equivalently \vspace{-2mm} 
\begin{equation}  \label{f.estimloiw*} 
\P_0^0\big(\o^* < \sqrt{\pi}\, y\big) = \, \frac{2}{y}\, e^{-\, y\2/4}\, \sum_{n\in \sN}\, e^{-n(n+1)/y^2} 
\, <\, \frac{2\, e^{-\, y\2/4}}{x\, (1-e^{-\,2/ y^2})}\, \raise1.9pt\hbox{.} 
\end{equation} 
By (\ref{f.estimsupPV<x0}) and (\ref{f.estimloiw*}), setting $\,C_p := 2\times 5^{2 p}\times (48)^{p-1- {1}/{p}}$,  
we finally  obtain\,: \  for $\,0\le x\le 1$, 
$$ \sup_{\theta}\, \P_0^0\Big(\inf_{0\le \e\le 1}\limits V^\theta_\e <x\Big) \le  C_p\, \E_0^0\big(S_{p}^{2p}\big)\, x^{(2p-3- 3/p)/4} + \min\!\left\{2\,\raise1pt\hbox{,}\, \frac{4\sqrt{\pi}\, x^{-1/8}\, e^{-\, \frac{\pi}{196}\, x^{-1/4}}}{7 \left(1-e^{-\,\frac{2\pi}{49}\, x^{-1/4}}\right)}\right\} \raise1.6pt\hbox{.} \;\;\diamond $$ 
 
\subsubsection{Key domination relating to the singular case $w=0$}  \label{sec.KDomin0} \indf 
    The following crucial estimate is the singular analogue of Proposition \ref{pro.estintdecre}. 
\bpro \label{pro.estintdecre0} \  For any fixed positive $C$, we have 
{
$$ \sup_{0\le \e\le 1}\, \E_0^{0}\!\left[\exp\!\bigg(\! -  \frac{\det(D\tilde U^{\e})}{\Tr(D\tilde U^{\e})\, \tilde \Omega^{\e}}\times  \frac{\sqrt{\xi^2+\xi'\!\,^2}}{147\,C}\,\bigg)\right] \in L^{1/2}(\R^2, d\xi'd\xi) \hbox{.} \vspace{-2mm} $$ } 
\epro 
\ub{Proof} \quad  Recall that according to (\ref{f.estimw=0}), we have to show that (\ref{f.butdomw=0}) holds, which by Section \ref{sec.AstH} amounts to showing that (\ref{f.estimlambdt}) holds. Now, applying Proposition \ref{pro.MinVar-w=0} with $\,p\ge 34\,$ we obtain\,:  
$$ \left( \int_1^\infty \sup_{0\le \e\le 1}\limits\sqrt{\int_0^1 \sup_\theta\limits\, \P_0^0\!\left(V^{\theta}_\e < \, 2\,s \right)  e^{-R\, s}\,\frac{ds}{s}}  \; R^{13/2}\, dR\right)^{\!2}  $$ 
$$ \le \int_1^\infty \sup_{0\le \e\le 1}\limits \left[\int_0^1 \sup_\theta\limits\, \P_0^0\!\left(V^{\theta}_\e < 2 s \right)  e^{-R\, s}\,\frac{ds}{s}\right] R^{15}\, dR\, \times \int_1^\infty R^{-2}\, dR \vspace{1mm} $$ 
$$ \le\, \int_1^\infty \left[\int_0^1C_p\, \E_0^0\big(S_{p}^{2p}\big)\, (2s)^{(2p-3- 3/p)/4} \,e^{-R\, s}\,\frac{ds}{s}\right] R^{15}\, dR \qquad $$ 
$$ \qquad + \int_1^\infty \left[\int_0^1\min\!\left\{2\,\raise1pt\hbox{,}\, \frac{4\sqrt{\pi}\, (2s)^{-1/8}\, e^{\ts-\, \frac{\pi}{196} (2s)^{-1/4}}}{7 \left(1-e^{-\,\frac{2\pi}{49} (2s)^{-1/4}}\right)}\right\} e^{-R\, s}\,\frac{ds}{s}\right] R^{15}\, dR $$ 
$$ =\, C_p\, \E_0^0\big(S_{p}^{2p}\big) 2^{(2p-3- 3/p)/4} \int_1^\infty \left[\int_0^R s^{(2p-3- 3/p)/4-1} \,e^{- s}\, ds\right] R^{15-(2p-3- 3/p)/4}\, dR \qquad $$ 
$$ \qquad + \int_0^1 \left[ \int_1^\infty e^{-R\, s}\,R^{15}\, dR \right] \min\!\left\{2\,\raise1pt\hbox{,}\, \frac{4\sqrt{\pi}\, (2s)^{-1/8}\, e^{\ts -\, \frac{\pi}{196} (2s)^{-1/4}}}{7 \left(1-e^{-\,\frac{2\pi}{49} (2s)^{-1/4}}\right)}\right\} \frac{ds}{s} $$ 
$$ =\, C_p\, \E_0^0\big(S_{p}^{2p}\big) 2^{(2p-3- 3/p)/4} \int_0^\infty \left[\int_s^\infty  R^{15-(2p-3- 3/p)/4}\, dR\right] s^{(2p-3- 3/p)/4-1} \,e^{- s}\, ds \qquad \vspace{-1mm} $$ 
$$ \qquad + \int_0^1 \left[ \int_s^\infty e^{-R}\,R^{15}\, dR \right] \min\!\left\{2\,\raise1pt\hbox{,}\, \frac{4\sqrt{\pi}\, (2s)^{-1/8}\, e^{\ts-\, \frac{\pi}{196} (2s)^{-1/4}}}{7 \left(1-e^{\ts -\,\frac{2\pi}{49} (2s)^{-1/4}}\right)}\right\} s^{-17}\, ds $$ 
{
$$ \le\, C_p\, \E_0^0\big(S_{p}^{2p}\big) \frac{2^{(2p+5- 3/p)/4}}{2p-67-3/p} \int_0^\infty\! e^{- s} s^{15} ds + 2^{17} (15)!\! \int_0^2\! \min\!\left\{1\raise1pt\hbox{,}\, \frac{2\sqrt{\pi}\, s^{-1/2}\, e^{\ts\frac{-\pi}{196}\, s^{-1}}}{7 \left(1-e^{-\,\frac{2\pi}{49}\, s^{-1}}\right)}\right\} s^{-65}\, ds $$ } 
which (owing to Remark \ref{rem.modcontpontB}) is clearly finite, for example for $\,p=34\,$. $\;\diamond$

\subsection{Intermediate small-time equivalent for $\,p_\e\,$ in the case $\,w=0$}  \label{sec.KApprox0} \indf 
    The preceding estimate entails the following wanted approximation result, which is the singular analogue of Propositions \ref{pro.equivIntPe} and \ref{pro.calcavpontFb0}. 
\bpro  \label{pro.equivIntPe0} \   Uniformly with respect to $(y,z)\in \R^2\,$, as $\,\e\sea 0\,$ we have
$$ p_\e\big(0\, ; (0,y, z)\big) \, =\, \frac{1+o(1)}{4\pi^2\,\e^{4}\,\sqrt{2\pi}}\times\II_\e(y,z)\,, \vspace{-2mm}  $$
with \vspace{-2mm}  
$$ \II_\e(y,z) := \int_{\sR^2} \exp\!\left[{ \rt1\! \left({\frac{\xi'}{\e^2}} (\e - y) - {\frac{\xi}{\e^{3/2}}}\, z\right)}\right] \tilde\EE_0(\xi',\xi)\, d\xi'd\xi\, ,\vspace{-1mm} $$ 
and \qquad ${\ds  \tilde\EE_0(\xi',\xi) =\, \E_0^{0}\bigg[\! \exp\!\bigg(\rt1\! \int_0^{1}\! \left[ \xi\, \o_s - {\xi'}\, \o_s^2/2 \right] ds \!\bigg) \bigg] \in L^1(\R^2, d\xi'd\xi)}$. 
\epro 
\ub{Proof} \quad  By (\ref{f.expr70'}), (\ref{fde.Ueo}) and (\ref{df.Etilde}) (in the beginning of Section \ref{sec.casew=0}), we have 
$$ p_\e\big(0\, ; (0, \e\, y\,, \e\, z)\big) =\, \frac{1}{4\pi^2\,\e^{4}\,\sqrt{2\pi}}\,  \int_{\sR^2} P^0_\e(\xi',\xi)\, d\xi'd\xi\,,  \vspace{-1mm} $$ 
with 
$$ P^0_\e(\xi',\xi) = \exp\!\left[{ \rt1\! \left({\frac{\xi'}{\e}} (1 - y) - {\frac{\xi}{\sqrt{\e}}}\, z\right)}\right] \! \times \tilde \EE_\e(\xi',\xi) \vspace{-0mm}\,  . $$ 
By (\ref{f.applT20}) and Proposition \ref{pro.estintdecre0} we have \  $\sup_{0\le \e\le1}\limits \big|\tilde \EE_\e(\xi',\xi)\big| \in L^1(\R^2, d\xi'd\xi)$. \  This provides the wanted domination of  $\, P_\e^0(\xi',\xi) $, which allows  to apply the Lebesgue theorem with respect to $d\xi'd\xi\,$, using the continuity of $\,\tilde\EE_\e(\xi',\xi)$ at $\,\e=0\,$ and then (\ref{fde.Ueoo}). As $\,\tilde\EE_\e(\xi',\xi)$ does not depend on $(y,z)$, this convergence is uniform with respect to $(y,z)\in \R^2$, which we can finally replace by $(y,z)/\e\,$. $\;\diamond$  

\section{Computation of a quadratic Laplace transform} \label{sec.calculFonctQuad} \indf
    Proposition \ref{pro.equivIntPe0} yields the wanted equivalent for $ p_\e\big(0\, ; (0,y, z)\big)$ in terms of the  Laplace transform ${\ds  \tilde\EE_0(\xi',\xi)}$ of a quadratic functional of the Brownian bridge. 
    
    We here perform the computation of a slightly more general Brownian bridge Laplace transform. This was already needed in Section \ref{sec.Analw=0}, and will be crucial in the forthcoming section \ref{sec.CompLoi(om,om^2)}. The principle of this type of computation goes back to Yor [Y]. \vspace{-2mm} 
\bPro \label{lem.resultcalcGpont} \   We have \quad ${\ds \E_0^{0}\bigg[ \exp\!\bigg(\! \int_0^{1}\! \big[ \alpha_s \, \o_{s} +  \gamma_s\, \o_s^2\big] ds\! \bigg)\bigg] =  }$ \vspace{-2mm} 
$$ = \, \frac{\exp\!\bigg[ \5\int_0^1\!\big(\!\int_s^1 e^{\int_\tau^s g}\, \alpha_\tau\,d\tau\! \big)^{\!2} ds -\,  { {\big(\int_0^1\big(\!\int_s^1e^{\int_\tau^s g}\, \alpha_\tau\, d\tau\! \big) e^{\int_1^s g}\, ds\big)^{\!2} }{\big(2\int_0^1e^{2\int_1^s g}\, ds}\big)^{\!\1} } - \5\int_0^1 g \bigg]}{{\sqrt{\int_0^1e^{2\int_1^s g}\, ds}}}\, \raise1.5pt\hbox{,} \vspace{-2mm}  $$ 
where $\,\alpha_s\,,\,\gamma_s\,$ are real deterministic, $\,\gamma_s\le 0\,$, and $\,g$ solves the Riccati equation \  $g'= g^2 + 2 \gamma$ \big(equivalent to the linear equation \  
$ \frac{d^2}{ds^2}\,\exp\!\big(-\int^s_0 g\big) = -2 \gamma_s\, \exp\!\big(-\int^s_0 g\big)$\big) a.e.  on $[0,1]$.
\ePro 
\ub{Proof} \quad  Set \quad ${\ds J^v:= \E_0^{v}\bigg[ \exp\!\bigg(\! \int_0^{1}\! \big[ \alpha_s \, \o_{s} +  \gamma_s\, \o_s^2\big] ds\! \bigg)\bigg]}$, \   so that \parn 
\begin{equation}  \label{f.defJvY}  
 \int_{\sR} \frac{e^{-(\theta+1)v^2/2}}{\sqrt{2\pi}}\, J^v\, dv\, =\, \E_0\!\bigg[ \exp\!\bigg(\! \int_0^{1}\! \big[ \alpha_s \, \o_{s} +  \gamma_s\, \o_s^2\big] ds - \5\, \theta\, \o_1^2 \bigg)\bigg] =:\, Y_\theta\,. 
\end{equation} 
Consider the exponential martingale defined by 
$$ M_s^{g} := \exp\!\bigg(\! - \int_0^{s}\! g_\tau\,\o_\tau\, d\o_\tau - \5\!\int_0^{s}\! g_\tau^2\,\o_\tau^2\, d\tau\!\bigg)   
= \exp\!\bigg(\5 \int_0^{s} g\, - \5\,g_s\,\o_s^2 + \5\!\int_0^{s}\big(g'_\tau - g_\tau^2\big)\,\o_\tau^2\, d\tau\bigg) .  $$ 
Denoting by $\,\P^{g}$ the new probability law having $M_s^{g}$ as density on $\FF_s\,$ with respect to $\P_0\,$,  we have\,: \vspace{-1mm} 
$$ Y_\theta\; e^{\int_0^1f/2} =\, \E^g\bigg[ \exp\!\bigg(\! \5 (g_1-\theta)\o_1^2 + \int_0^{1}\! \Big[ \alpha_s\,\o_{s} +  \big(\gamma_s -\5 (g'_s-g^2_s) \big)\,\o_s^2\Big] ds\! \bigg)\bigg] \vspace{-1mm}  $$ 
\begin{equation}  \label{f.calcintermY} 
=\, \E^g\bigg[ \exp\!\bigg(\! \5 (g_1-\theta)\o_1^2 + \int_0^{1} \alpha_s \,\o_{s} \,ds \bigg)\bigg] , 
\end{equation}
by taking $\,g\,$ almost everywhere solving the Ricatti equation \  $g'= g^2 + 2 \gamma$ \big(equivalent to the linear equation \  $\exp\!\big[-\int g\big] '' = -2 \gamma \exp\!\big[-\int g\big]$\big). \parn 
On the other hand, the Girsanov formula provides a $(\P^g,\FF_s)$ Brownian motion $B$ such that \  ${\ds \o_s= B_s - \int_0^s g_\tau\,\o_\tau \, d\tau}\,$, and then \  ${\ds \o_s=  \int_0^s \exp\!\Big(\int_s^\tau g\Big) dB_\tau}\,$. \   Hence, for any real $\,r\,$: 
$$ \E^g\bigg[\Big(\! r\,\o_1 + \int_0^{1} \alpha_s \,\o_{s} \,ds \Big)^2\bigg] = \E^g\bigg[\Big( \int_0^{1} \Big[r +\int_s^1\alpha\Big] d\o_s \Big)^2\bigg] $$
$$ = \E^g\bigg[\bigg( \int_0^{1} \Big[r +\int_s^1\alpha\Big] \Big[ dB_s - g_s \Big(\int_0^s e^{\int_s^\tau g}\, dB_\tau \Big) ds  \Big]  \bigg)^2\bigg]  $$
$$ = \E^g\bigg[\bigg( \int_0^{1}  \Big[r +\int_s^1\alpha - \int_s^1 \Big(r +\int_\tau^1\alpha\Big) g_\tau\, e^{\int_\tau^s g} d\tau\Big] dB_s \bigg)^2\bigg]  $$ 
$$ = \int_0^{1}  \Big[r \, e^{\int_1^s g} +\int_s^1 e^{\int_\tau^s g}\, \alpha_\tau\, d\tau \Big]^2 ds \, .    $$
This yields the covariance matrix of the $\,\P^g$-Gaussian variable ${\ds \Big( \o_1 \, , \int_0^{1}\alpha_s \,\o_{s} \,ds\Big)}$, namely  
$$ K = \begin{pmatrix} \int_0^1 e^{2\int_1^s g} ds\; & \int_0^1\big(\!\int_s^1e^{\int_\tau^s g}\, \alpha_\tau\, d\tau\! \big) e^{\int_1^s g} ds  \cr 
\int_0^1 \big(\!\int_s^1e^{\int_\tau^s g}\, \alpha_\tau\, d\tau\! \big) e^{\int_1^s g} ds\; & \int_0^1 \big(\!\int_s^1e^{\int_\tau^s g}\, \alpha_\tau\, d\tau\! \big)^2 ds
\end{pmatrix}\!  ,  $$
whence the density of ${\ds \Big(\o_1 \, , \int_0^{1}\alpha_s \,\o_{s} \,ds\Big)}$ with respect to $\,\P^g$.  \ Thus 
$$ \E^g\bigg[ \exp\!\bigg(\! \int_0^{1}\! \alpha_s\, \o_{s} \,ds - r\,\o_1^2 \bigg)\bigg] =  
\int_{\sR^2} e^{v- r\, u^2} \exp\!\bigg[  \frac{-1}{2\,\det K} \bigg( u^2\! \int_0^1\! \Big(\!\int_s^1e^{\int_\tau^s g}\, \alpha_\tau\, d\tau\! \Big)^2 ds \hskip 1mm $$
$$ \hskip 10mm - 2\,uv \int_0^1\big(\!\int_s^1e^{\int_\tau^s g}\, \alpha_\tau\, d\tau\! \big) e^{\int_1^s g} ds +  v^2\! \int_0^1 e^{2\int_1^s g} ds\bigg) \bigg] \frac{du\, dv}{2\pi\sqrt{\det K}}  \, = $$
$$ = \int_{\sR^2}\! \exp\!\bigg[{\ts\frac{-1}{2}} \bigg( \bigg[{\ts \frac{\int_0^1 \left(\int_s^1e^{\int_\tau^s\! g}\, \alpha_\tau\, d\tau\! \right)^2 ds}{\det K}} +2r \bigg] u^2  + \bigg[{\ts \frac{\int_0^1 e^{2\int_1^s\! g} ds}{\det K}} \bigg] v^2 -2v \hskip 19mm $$
$$ \hskip 40mm- 2\bigg[{\ts \frac{\int_0^1\left(\int_s^1e^{\int_\tau^s g}\, \alpha_\tau\, d\tau\! \right) e^{\int_1^s g} ds}{\det K}}\bigg] uv  \bigg) \bigg] \frac{du\, dv}{2\pi\sqrt{\det K}}  $$ 
\begin{equation} \label{f.resultcalcG} 
= \left(1+2r\! \int_0^1 e^{2\int_1^s\! g} ds \right)^{\!-1/2}\times \exp\!\left(\frac{\int_0^1 \big(\!\int_s^1e^{\int_\tau^s g}\, \alpha_\tau\, d\tau\! \big)^2 ds  + 2r\, \det K }{2 \left(1+2r\! \int_0^1 e^{2\int_1^s\! g} ds \right)} \right) 
 \end{equation}  
(by a classical Gaussian computation). \  Finally by (\ref{f.defJvY}), (\ref{f.calcintermY}) and (\ref{f.resultcalcG}) we have 
$$ \E_0^{0}\bigg[ \exp\!\bigg(\! \int_0^{1}\! \big[ \alpha_s \,\o_{s} +  \gamma_s\, \o_s^2\big] ds\! \bigg)\bigg]= J^0 = \lim_{\theta\to\infty} \sqrt{\theta+1} \int_{\sR} \frac{e^{-(\theta+1)\,v^2/2}}{\sqrt{2\pi}}\, J^v\, dv $$
$$ =\, \lim_{\theta\to\infty} \sqrt{\theta+1}\,\, Y_\theta\, = \lim_{r\to\infty} \sqrt{2r}\times \E^g\bigg[ \exp\!\bigg(\! \int_0^{1}\! \alpha_s\,\o_{s} \,ds - r\,\o_1^2 \bigg)\bigg] \times e^{-\int_0^1g/2}$$
$$ = \, \left(\int_0^1 e^{2\int_1^s g}\, ds \right)^{\!-1/2}\! \times \exp\!\left(\frac{ \det K }{2 \int_0^1 e^{2\int_1^s\! g}\, ds } -\5 \int_0^1g\right) $$ 
$$ =\, \frac{\exp\!\bigg(\frac{ \left(\int_0^1 e^{2\int_1^s g} ds\right)\! \left(\int_0^1 \big(\!\int_s^1e^{\int_\tau^s g}\, \alpha_\tau\, d\tau\! \big)^2 ds \right) - \left( \int_0^1\big(\!\int_s^1e^{\int_\tau^s g}\, \alpha_\tau\, d\tau\! \big) e^{\int_1^s g} ds \right)^{\!2} }{2 \int_0^1 e^{2\int_1^s g}\, ds} -\5 \int_0^1g\bigg)}{\sqrt{\int_0^1 e^{2\int_1^s\! g}\, ds }}\,\raise2pt\hbox{.} \;\; \diamond  $$ 
   
\bCor  \label{cor.calculQuadrF}  \  For any real constants $\,\alpha\,$ and $\,\gamma\,$, \  we have 
$$ \E_0^{0}\bigg[ \exp\!\bigg(\! \int_0^{1}\! \big[ \alpha\, \o_{s} -\5 \gamma^2\, \o_s^2\big] ds\! \bigg)\bigg] = \sqrt{\frac{\gamma}{\sh \gamma}}\, \exp\!\bigg[ \frac{\alpha^2}{\gamma^3}\! \left( {\frac{\gamma}{2}} -   {\th\!\big(\frac{\gamma}{2}}\big) \right) \!\bigg] \vspace{-2mm} . $$ 
\eCor 
\ub{Proof} \quad We apply Proposition \ref{lem.resultcalcGpont} with $\, \gamma_s\equiv -\gamma^2/2\,$, so that we can take $\,g \equiv \gamma\,$, yielding\,: \  
$$ \E_0^{0}\bigg[ \exp\!\bigg(\! \int_0^{1}\! \big[ \alpha\, \o_{s} -\5 \gamma^2\, \o_s^2\big] ds\! \bigg)\bigg] $$ 
$$ = \, \frac{\exp\!\bigg[ \frac{\alpha^2}{2}\int_0^1\!\big(\frac{1-e^{\gamma (s-1)}}{\gamma}\big)^{\!2} ds -\,  { {\big(\int_0^1\big(\frac{1-e^{\gamma(s-1)}}{\gamma}\big) e^{\gamma (s-1)}\, ds\big)^{\!2} } \frac{\alpha^2 \gamma}{1-e^{-2\gamma}}  } - \5\,\gamma \bigg]} {{\sqrt{\frac{1-e^{-2\gamma}}{2\gamma} }}} $$ 
$$ =  \sqrt{\frac{2\gamma}{1-e^{-2\gamma}} }\times \exp\!\bigg[ \frac{\alpha^2}{2\gamma^2} \int_0^1\!\big(1-e^{\gamma (s-1)}\big)^{\!2} ds - { {\left(\int_0^1\!\big[e^{\gamma (s-1)}-e^{2\gamma (s-1)}\big] ds\right)^{\!2} } \frac{\alpha^2/\gamma}{1-e^{-2\gamma}} } - \5\,\gamma \bigg] $$ 
$$ = \sqrt{\frac{2\gamma\, e^{-\gamma}}{1-e^{-2\gamma}} }\times \exp\!\bigg[ \frac{\alpha^2}{2\gamma^2}\times \left(1-\frac{3-4\, e^{-\gamma} + e^{-2\gamma}}{2\gamma}\right)  - \frac{(1-e^{-\gamma} )^4}{4\gamma^2}\times \frac{\alpha^2/\gamma}{1-e^{-2\gamma}} \bigg] $$ 
$$ = \sqrt{\frac{\gamma}{\sh \gamma}}\, \exp\!\bigg[ \frac{\alpha^2}{4\gamma^3} \big(2\gamma - (1-e^{-\gamma})(3-e^{-\gamma})\big) - \frac{\alpha^2}{4\gamma^3}\,\frac{(1-e^{-\gamma})^3}{1+e^{-\gamma}}\bigg] $$ 
$$ = \sqrt{\frac{\gamma}{\sh \gamma}}\, \exp\!\bigg[ \frac{\alpha^2}{\gamma^3}\! \left( {\frac{\gamma}{2}} -   {\th\!\big(\frac{\gamma}{2}}\big) \right) \!\bigg]  .\;\;\diamond \vspace{-1mm} $$ 

\bRem  \label{remRiccati} \  {\rm The above corollary \ref{cor.calculQuadrF} contains the following well known particular case\,: for any real constant $\gamma\,$, we have \vspace{-1mm}
$$ \E_0^{0}\bigg[ \exp\!\bigg(-\!{\frac{\gamma^2}{2}} \int_0^{1}\! \omega_s^2\, ds\! \bigg)\bigg] = \sqrt{\frac{\gamma}{\sh\gamma}}\, \raise1.5pt\hbox{.} \vspace{-2mm} $$ 
Now for $\,0<\beta < 2\,$, by (\ref{f.estYor}) we have 
$$ \E_0^{0}\bigg[ \exp\!\bigg(\beta\!  \int_0^{1}\! \omega_s^2\, ds\! \bigg)\bigg] \le \E_0^{0}\Big[ e^{\beta\,  (\omega^*)^2}\Big] = \beta \int_0^\infty \P_0^0\big[\omega^*>\sqrt{t}\,\big] e^{\beta\, t} dt  \, \le\, 2\beta \int_0^\infty\! e^{(\beta-2)\, t} dt  < \infty\, \raise0pt\hbox{.} \vspace{-1mm} $$ 
By analytical continuation, this implies that \vspace{-1mm}
$$ \E_0^{0}\bigg[ \exp\!\bigg({\frac{\gamma^2}{2}} \int_0^{1}\! \omega_s^2\, ds\! \bigg)\bigg] = \sqrt{\frac{\gamma}{\sin\gamma}} $$ 
is analytical on $\{0\le \gamma < 2\}$, and then on $\{0\le \gamma < \pi\}$ as well by monotone convergence (and it explodes at $\,\gamma = \pi$). Similarly, for any positive $\,q\,$ we have\,: 
$$ \E_0^{0}\bigg[ \exp\!\bigg(q\beta\!\int_0^{1}\! \omega\, + \beta\! \int_0^{1}\! \omega^2\bigg)\bigg] \le \E_0^{0}\Big[ e^{q\beta\, \omega^*+ \beta\,  (\omega^*)^2}\Big] = \beta \int_0^\infty \P_0^0\big[(\omega^*)^2+q\omega^*> {t}\,\big] e^{\beta\, t} dt $$
$$ = \beta \int_0^\infty \P_0^0\big[2\,\omega^*> \sqrt{q^2+4t}-q\big] e^{\beta\, t} dt\, \le\, 2\beta \int_0^\infty\! e^{(\beta-2)\, t +q \sqrt{q^2+4t} -q^2} dt  < \infty\, \raise0pt\hbox{.} \vspace{-1mm} $$ 
This entails that the function $\,{\ds (\alpha,\beta)\lmt  \E_0^{0}\bigg[ \exp\!\bigg(\! \int_0^{1}\! \big[ \alpha\, \o_{s} + \beta\, \o_s^2\big] ds\! \bigg)\bigg] }\,$  is analytical  on $\,\C \times \{\Re(\beta) <2\}$, and even on $\,\C \times \{\Re(\beta) <\pi^2/2\}$\,: $\,\pi^2/2\,$ is the abscissa of convergence. 
}\eRem  

   Now the expression found in Corollary \ref{cor.calculQuadrF}, valid for $\, -\infty < \beta<2$ (and not for $\beta = \pi^2/2$), is not straightforwardly analytically continued in the whole half plane $\{\Re(\beta) <2\}$, whereas we need an explicit analytic continuation. The delicate point is of course  the square root coming in the expression of Corollary \ref{cor.calculQuadrF}. We shall use the following lemma. \vspace{-0mm} 
\bLem \label{lem.ProlAnalEspCompl} \  $(i)$ \  For any positive real $\,a,\,b$, we have 
$$  \sh\!\big[a+\rt1b\big]\, =\, \sqrt{\frac{\ch\!(2a)-\cos(2b)}{2}}\times \exp\!\left[\rt1\!\! \int_0^b  \frac{\sh\!(2a)\, d\beta }{\ch\!(2a)-\cos(2\beta)}\right] \raise1pt\hbox{.} $$ 
$(ii)$ \  For any real $\,\chi\,$ and positive $\,x$, we have the following continuous lift (i.e., determination) of the square root\,: 
$$ \sqrt{\frac{\sqrt{\chi+\rt1\! x}}{{\sh}\!\sqrt{\chi+\rt1\! x}}}\, = \left[{\frac{2\,\sqrt{\chi^2+x^2}}{\ch\!\sqrt{ 2\sqrt{\chi^2+x^2}+2\chi}  -\cos\!\sqrt{ 2\sqrt{\chi^2+x^2}-2\chi}}}\right]^{1/4} \times e^{\rt1 \f(\chi,x)} \raise1pt\hbox{,}  \vspace{-2mm} $$
with \vspace{-2mm} 
$$ \f(\chi,x) = \frac{1}{2}\, \arctg\sqrt{\frac{\sqrt{\chi^2+x^2}-\chi}{\sqrt{\chi^2+x^2}+\chi}}\, - \frac{1}{2} \int_0^{\sqrt{\frac{\sqrt{\chi^2+x^2}-\chi}{2}}}\! \frac{\sh\!\sqrt{ 2\sqrt{\chi^2+x^2}+2\chi}\; d\beta }{\ch\!\sqrt{ 2\sqrt{\chi^2+x^2}+2\chi}-\cos(2\beta)}\, \raise1.9pt\hbox{} $$ 
$${\ts  = \frac{1}{2}\, \arctg\sqrt{\frac{\sqrt{\chi^2+x^2}-\chi}{\sqrt{\chi^2+x^2}+\chi}}\, - \frac{k\pi}{2} - \frac{1}{2} \,\arctg\!\left(\tg\!\!\left[\sqrt{\frac{\sqrt{\chi^2+x^2}-\chi}{2}}-k\pi\right]\times \coth\!\!\left[\sqrt{\frac{\sqrt{\chi^2+x^2}+\chi}{2}}\,\right]\right) } $$
for $\,k\in\N\,$ and $\,\left|\sqrt{\frac{\sqrt{\chi^2+x^2}-\chi}{2}}-k\pi\right| \le {\ds\frac{\pi}{2}}\,$\raise1.4pt\hbox{.}  
\eLem
\ub{Proof} \quad  $(i)$ \  The modulus of \  $\sh\!\big[a+\rt1b\big]\, =\,\sh a\,\cos b + \rt1 \ch a \, \sin b$ \   is \parn
$\big| \sh\!\big[a+\rt1b\big]\big|  = \sqrt{\sh\!^2a + \sin^2\!b} =\! \sqrt{\frac{\ch\!(2a)-\cos(2b)}{2}}>0$, and setting 
 \  $\frac{\sh\!\left[a+\rt1b\right]}{\sqrt{\sh\!^2a + \sin^2\!b}} = e^{\rt1\!\f_a(b)}\,$  with $\,\f_a(0)=0\,$, we have   \  $\f'_a(b) = \frac{\sh\!(2a)}{\ch\!(2a)-\cos(2b)}\, $\raise1.9pt\hbox{,} \  since 
$$ \frac{\sh a \, \cos b\times\f'_a(b)}{\sqrt{\sh\!^2a + \sin^2\!b}} = \frac{\partial}{\partial b} \sin\f_a(b) = \frac{\partial}{\partial b}\frac{\ch a \, \sin b}{\sqrt{\sh\!^2a + \sin^2\!b}}  = \frac{\ch a\, \cos b\times \sh\!^2a}{(\sh\!^2a + \sin^2\!b)^{3/2}}\, \raise1.9pt\hbox{.} \;\; $$ 
Finally
$$ \f_a(b) =  \int_0^b  \frac{\sh\!(2a)\, d\beta }{\ch\!(2a)-\cos(2\beta)}\,, \   \hbox{ whence }\;\f_a(b)  =\, \arctg\big[\coth a\times \tg b\big]  \  \hbox{ for } \; 0\le b\le \pi/2\,, $$ 
and then \  $\, \f_a(b) = k\pi + \arctg\big[\coth a\times \tg (b-k\pi)\big] $ \  for $\; k\pi-\frac{\pi}{2} \le b\le k\pi+\frac{\pi}{2}\,\raise1pt\hbox{,}\; k\in\N\,$. \parn
In particular, we have $\;\f_a(k\pi/2) = k\pi/2\;$ for any $\,k\in\N\,$.  \parsn 
$(ii)$ \  Taking the usual determination of the square root on $\C\moins \rt1\R_-\,$, we first have 
\begin{equation}  \label{f.racinexkhi} 
\sqrt{\chi+\rt1 x} = a+\rt1 b\; \hbox{ with }\; a:= \sqrt{\frac{\sqrt{\chi^2+x^2}+\chi}{2}}\,, \; b:= \sqrt{\frac{\sqrt{\chi^2+x^2}-\chi}{2}}\,\raise1.9pt\hbox{,} 
\end{equation}
so that applying $(i)$ above we have 
{\small 
$$ \frac{\sqrt{\chi+\rt1\! x}}{{\sh}\!\sqrt{\chi+\rt1\! x}}\, = \sqrt{\ts\frac{2\,\sqrt{\chi^2+x^2}}{\ch\!\sqrt{ 2\sqrt{\chi^2+x^2}+2\chi}  -\cos\!\sqrt{ 2\sqrt{\chi^2+x^2}-2\chi}}}\,\times \exp\!\left[\rt1\! \big(\arctg (b/a) - \f_a(b)\big) \right]  . \;\;\diamond$$ }

\bCor \label{cor.ValFourIntomt}  \ For any real $\xi\,$, any $\,\chi > -\pi^2\,$ and $\,x>0\,$,  we have\,: 
$$ \E_0^{0}\bigg[ \exp\!\bigg(\!\int_0^{1}\! \Big[ \rt1\!\xi\, \o_{s} - \frac{\chi+\rt1\! x}{2}\, \o_s^2\,\Big] ds\! \bigg)\bigg] =\, e^{\rt1 \f(\chi,x)}\, \times \hskip22mm \vspace{-2mm}  $$
$$ \hskip12mm \left[{\frac{2\,\sqrt{\chi^2+x^2}}{\ch\!\sqrt{ 2\sqrt{\chi^2+x^2}+2\chi}  -\cos\!\sqrt{ 2\sqrt{\chi^2+x^2}-2\chi}}}\right]^{1/4} \times  \exp\!\left[\frac{-\, \xi^2/2}{\chi+\rt1\!x}\times f(\chi,x) \!\right] \! , \vspace{-0mm} $$ 
with $\,\f(\chi,x)$ as in Lemma \ref{lem.ProlAnalEspCompl}, and setting $\,\sqrt{\chi+\rt1 x} = a+\rt1 b\,$ as in Lemma \ref{lem.ProlAnalEspCompl}\,: 
\begin{equation}  \label{df.fnf} 
f(\chi,x) := \, 1 -   {\frac{\th\!\!\left[\sqrt{\chi+\rt1\!x}\big/2\,\right]}{\sqrt{\chi+\rt1\!x}\big/2}}\, =\, 1- \frac{a\, \sh a + b\, \sin b + \rt1 (a\, \sin b - b\,\sh a)}{(a^2+b^2)\, (\ch a + \cos b)\big/ 2}\,\vspace{-2mm} \raise1.9pt\hbox{.} 
\end{equation} 
\eCor 
\ub{Proof} \quad  Remark \ref{remRiccati} and Lemma \ref{lem.ProlAnalEspCompl} allow to continue Corollary \ref{cor.calculQuadrF} analytically,  to $\,\gamma = \sqrt{\chi+\rt1\!x}\,$, at least for  $\, \chi > -\pi^2$ and $\,x>0$. Lemma  continues the term  $\sqrt{\frac{\gamma}{\sh \gamma}}\, $ of Corollary \ref{cor.calculQuadrF}.  And the remaining term is 
$$ \exp\!\bigg[ \frac{-\,\xi^2}{2\, \gamma^2}\! \left( 1 - {\frac{2}{\gamma}}\times \th\!\!\left[\frac{\gamma}{2}\right] \right) \!\bigg] =   \exp\!\left[ \frac{-\, \xi^2/2}{\chi+\rt1\!x}\, f(\chi,x)\right] ,  $$ 
whence the first expression of the statement for $\,f(\chi,x)$. Using the notation$\,a, b\,$, the second expression follows at once from the observation that  \  
$\th\!\!\left[\frac{\sqrt{\chi+\rt1\!x}}{2}\right] =\, \frac{\sh a + \rt1\! \sin b}{\ch a + \cos b}\, \raise1pt\hbox{.} \;\;\diamond $ 

\section{Evaluation of the oscillatory integral $\,{\ds \II_\e(y,z) }$}  \label{sec.CompLoi(om,om^2)} \indf 
    Proposition \ref{pro.equivIntPe0} yields the wanted equivalent for $ p_\e\big(0\, ; (0,y, z)\big)$, in terms of the sto\-chastic oscillatory integral $\,{\ds \II_\e(y,z) }$, which we have to evaluate carefully.  We hall first use Section \ref{sec.calculFonctQuad}, in order to express $\,{\ds \II_\e(y,z) }$ as a finite-dimensional oscillatory integral. Even so, it will however remain delicate to handle. 
    
\subsection{Reduction to a finite-dimensional oscillatory integral}  \label{sec.UsCalcQuadr} \indf 
   First we have the following.
\blem \label{lem.ExprIepsyz} \  For any $(y,z)\in \R^2\,$ and $\,\e>0\,$ we have \  
$$ \II_\e(y,z) =\, 2\; \Re \Bigg\{\! \int_{\sR}\exp\!\left[\rt1\!{\ts\frac{z}{\e^{3/2}}}\, \xi\right]\!\times \! \left[ \int_{0}^\infty\! \exp\!\left({ {\ts\frac{\e - y}{\e^2}} }\rt1\! x\right) \tilde\EE_0(x,\xi)\, dx\right] d\xi\Bigg\}  \vspace{-2mm} $$ 
$$ =\, 2\; \Re \Bigg\{\! \int_{0}^\infty\! \exp\!\left({ {\ts\frac{\e - y}{\e^2}} }\rt1\! x\right) \left[ \int_{\sR}\exp\!\left[\rt1\!{\ts\frac{z}{\e^{3/2}}}\, \xi\right]\!\times \! \tilde\EE_0(x,\xi)\, d\xi \right] dx \Bigg\} . \vspace{-2mm} $$ 
\elem 
\ub{Proof} \quad  Note first the obvious symmetries\,: \  $\tilde\EE_0(\xi',\xi) = \tilde\EE_0(\xi',-\xi) = \overline{\tilde\EE_0(-\xi',\xi)}\,$, due to the symmetry of $\P_0^0\,$, which permit in particular to restrict below to $\,\xi'\ge0\,$.  
Note that $\,\tilde\EE_0\,$ is a continuous function on $\R^2$, such that $\,\tilde\EE_0(0,\xi) = e^{-\xi^2/24}$ by Lemma \ref{lem.Gauss1'}. \   Since \  $\tilde\EE_0(\xi',\xi) \in L^1(\R^2, d\xi'd\xi)$, we can apply Fubini's Theorem, performing the integration with respect to $\xi'$ first. Using both symmetries of $\,\tilde\EE_0(\cdot,\xi)$, this yields\,: 
$$ \II_\e(y,z) =\, \int_{\sR}\exp\!\left[\rt1\!{\ts\frac{z}{\e^{3/2}}}\, \xi\right]\!\times \!\left[ \int_{-\infty}^\infty\! \exp\!\left({ {\ts\frac{\e - y}{\e^2}} }\rt1\! \xi' \right) \tilde\EE_0(\xi',\xi)\, d\xi'\right] d\xi $$ 
$$ =\, \int_{\sR}\exp\!\left[\rt1\!{\ts\frac{z}{\e^{3/2}}}\, \xi\right]\!\times \! \left[ \int_{0}^\infty\! \exp\!\left({ {\ts\frac{\e - y}{\e^2}} }\rt1\! \xi' \right) \tilde\EE_0(\xi',\xi)\, d\xi'\right] d\xi \,  $$ 
$$  \hskip15mm+\, \int_{\sR}\exp\!\left[\rt1\!{\ts\frac{z}{\e^{3/2}}}\, \xi\right]\!\times \! \left[ \int_{0}^\infty\! \exp\!\left({ -\,{\ts \frac{\e - y}{\e^2}} }\rt1\! \xi' \right) \tilde\EE_0(-\xi',\xi)\, d\xi'\right] d\xi \,  $$ 
$$ =\, \int_{\sR}\exp\!\left[\rt1\!{\ts\frac{z}{\e^{3/2}}}\, \xi\right]\!\times \! \left[ \int_{0}^\infty\! \exp\!\left({ {\ts\frac{1 - y}{\e}} }\rt1\! \xi' \right) \tilde\EE_0(\xi',\xi)\, d\xi'\right] d\xi \,  $$ 
$$ \hskip15mm +\, \int_{\sR}\exp\!\left[-\rt1\!{\ts\frac{z}{\e^{3/2}}}\, \xi\right]\!\times \! \left[ \int_{0}^\infty\! \exp\!\left({ -\,{\ts \frac{\e - y}{\e^2}} }\rt1\! \xi' \right) \overline{\tilde\EE_0(\xi',\xi)}\, d\xi'\right] d\xi \,  $$ 
$$ =\, 2\; \Re \Bigg\{ \int_{\sR}\exp\!\left[\rt1\!{\ts\frac{z}{\e^{3/2}}}\, \xi\right]\!\times \! \left[ \int_{0}^\infty\! \exp\!\left({ {\ts\frac{\e - y}{\e^2}} }\rt1\! \xi' \right) \tilde\EE_0(\xi',\xi)\, d\xi'\right] d\xi\Bigg\} .  $$ 
Finally we can apply Fubini's Theorem again, integrating now with respect to $\xi$ first. $\; \diamond$ \parm 

    Now we use Section \ref{sec.calculFonctQuad}, to deduce the following expression for the integral $ \II_\e(y,z)$ (we have to evaluate), as a finite-dimensional oscillatory integral. 
\bpro \label{pro.exprAvChgCont} \   For any $(y,z)\in \R^2\,$ and $\,\e>0\,$ we have \ 
$$ \II_\e(y,z) =\, 2\; \Re \Bigg\{\! \int_{0}^\infty\! \exp\!\left[{ {\frac{\e - y}{\e^2}} }\rt1\! x - \frac{z^2}{2\e^3}\, \frac{\rt1\! x}{f(0,x)}\right] \Phi(x)\, dx \Bigg\} \, , \vspace{-2mm} $$ 
where $\,\Phi$ is the analytical function given (for any $x\ge 0$) by\,: \vspace{-2mm} 
$$ \Phi(x) := \, e^{\rt1 \f(x)}\times  \left[{\frac{2\,x}{\ch\!\sqrt{2x}  -\cos\!\sqrt{2x}}}\right]^{1/4} \times \sqrt{\frac{2\pi\rt1\! x}{f(x)}}\, \raise1.9pt\hbox{,} $$ 
{\small 
$$ \frac{\rt1\! x}{f(0,x)} = \frac{x\sqrt{x/2}\, \big(\ch\!\sqrt{x/2} + \cos\!\sqrt{x/2}\,\big)}{\sh\!\sqrt{x/2} - \sin\!\sqrt{x/2} + \rt1\!\big(\sh\!\sqrt{x/2} + \sin\!\sqrt{x/2} - (\ch\!\sqrt{x/2} + \cos\!\sqrt{x/2}\,)\sqrt{x/2}\,\big)} \, \raise1.9pt\hbox{,}
$$ }
and \vspace{-2mm} 
$$ \f(0) = \,0\,, \quad  \f'(x) = \, \frac{-1}{2\sqrt{2x}} \times {\frac{\sh\!\sqrt{2x} - \sin\!\sqrt{2x}}{\ch\!\sqrt{2x} - \cos\!\sqrt{2x}}}\, \raise1.9pt\hbox{.}  $$ 
\epro  
\ub{Proof} \quad   By Corollary \ref{cor.ValFourIntomt}, for any real $\,\chi> - \pi^2$ and  $\,x\ge 0\,$  we have \vspace{-1mm} 
\begin{equation}  \label{f.E0tildx-ichi}
\tilde\EE_0\big(x-\rt1\! \chi\,,\,\xi\big) =\, \E_0^{0}\bigg[ \exp\!\bigg(\! \int_0^{1}\! \left[ \rt1\!\xi\, \o_{s} - \frac{\chi + \rt1\! x}{2}\, \o_s^2\right] ds\! \bigg)\bigg] \vspace{-2mm} 
\end{equation}
{\small 
$$ =\, e^{\rt1 \f(\chi,x)} \times  \left[{\frac{2\,\sqrt{\chi^2+x^2}}{\ch\!\sqrt{ 2\sqrt{\chi^2+x^2}+2\chi}  -\cos\!\sqrt{ 2\sqrt{\chi^2+x^2}-2\chi}}}\right]^{1/4} \times  \exp\!\left[\frac{-\, \xi^2/2}{\chi+\rt1\!x}\!\times\! f(\chi,x) \right] \! , \vspace{-2mm} $$  } \par \noindent 
with $\,f(\chi,x)$ given by (\ref{df.fnf}). \quad  In particular, \vspace{-2mm} 
$$ \tilde\EE_0\big(x\,,\,\xi\big) =\, e^{\rt1 \f(x)} \times  \left[{\frac{2\,x}{\ch\!\sqrt{2x}  -\cos\!\sqrt{2x}}}\right]^{1/4} \times  \exp\!\left[\frac{-\, \xi^2\, f(0,x)}{2\rt1\!x} \right] \! , \vspace{-1mm} $$
with \vspace{-1mm} 
$$ \frac{f(0,x)}{\rt1 }\, = \, \frac{\sh\!\sqrt{x/2} - \sin\!\sqrt{x/2} + \rt1\!\big(\sh\!\sqrt{x/2}+ \sin\!\sqrt{x/2}\,\big)}{(\ch\!\sqrt{x/2} + \cos\!\sqrt{x/2}\,)\sqrt{x/2}} -\rt1 , \vspace{-0mm} $$ 
and by Lemma \ref{lem.ProlAnalEspCompl}\,: \quad ${\ds  \f'(x) = \, \frac{-1}{2\sqrt{2x}} \times {\frac{\sh\!\sqrt{2x} - \sin\!\sqrt{2x}}{\ch\!\sqrt{2x} - \cos\!\sqrt{2x}}}\, \raise1.4pt\hbox{,}  }$ \quad and 
$$ \f(x) = \frac{\pi}{8} - \frac{k\pi}{2} - \frac{1}{2} \,\arctg\!\left[\tg\!\!\left(\!\sqrt{\frac{x}{2}}-k\pi\!\right)\times \coth\!\sqrt{\frac{x}{2}}\,\right]  \  
\hbox{ for $\,k\in\N$ and $\left|\sqrt{\frac{x}{2}}-k\pi\right| \le {\frac{\pi}{2}}\,$\raise1.4pt\hbox{.} } $$
Hence we obtain\,:  \vspace{-1mm} 
$$ \int_{\sR}\exp\!\left[\rt1\!{\ts\frac{z}{\e^{3/2}}}\, \xi\right]\times \tilde\EE_0(x,\xi)\, d\xi = e^{\rt1\! \f(x)}\! \left[{\ts\frac{2\,x}{\ch\!\sqrt{2x}  -\cos\!\sqrt{2x}}}\right]^{1/4}\! \sqrt{\ts\frac{2\pi\rt1 x}{f(0,x)}}\times  \exp\left[ -{\ts \frac{z^2}{2\e^3}\!\times\! \frac{\rt1\! x}{f(0,x)} }\right]\! . $$ 
Therefore, using Lemma \ref{lem.ExprIepsyz} we obtain the expression of the statement. \pars  
   Note then that \quad   $\ch\!\sqrt{2x} - \cos\!\sqrt{2x} = 0\, \LRa \,\sqrt{2x}\,\in  (\rt1\pm 1)\pi\Z\, \LRa \, x\in \pm\rt1\!\pi^2\N^2$. \parsn 
Note also that \  $\frac{\rt1\! x}{f(0,x)}\,$ is an even meromorphic function of $\sqrt{x}\,$, locally expandable in a power series of $x\,$. \  Moreover, we have  \vspace{-2mm} 
\begin{equation}  \label{f.partReelleb} 
\frac{\rt1}{f(0,2x^2)}\, =\, \frac{\rt1\!\left(1 - \frac{\sh x+\sin x}{x\,(\ch x + \cos x)}\right) + \frac{\sh x-\sin x}{x\,(\ch x + \cos x)}} {\left(1 - \frac{\sh x+\sin x}{x\,(\ch x + \cos x)}\right)^2 + \left(\frac{\sh x-\sin x}{x\,(\ch x + \cos x)}\right)^2 } \, \raise1.9pt\hbox{.} 
\end{equation}
Since $\,\left(x - \frac{\sh x+\sin x}{\ch x + \cos x}\right)$ has derivative $\,\frac{\sh^2 x-\sin^2 x}{(\ch x + \cos x)^2}\,>0\,$ on $\R_+^*$ and then continuously increases from 0 to infinity, we have $\left(1-\frac{\sh x+\sin x}{x\,(\ch x + \cos x)}\right) >0\, $ on $\R_+^*\,$, as well as $\,x\,\frac{\sh x-\sin x}{\ch x + \cos x}>0\,$\raise1pt\hbox{.} \parn 
In particular the real and imaginary parts in (\ref{f.partReelleb}) are positive continuous on $\,\R_+^*\,$, which shows that in the above the argument of $\,\frac{\rt1\! x}{f(0,x)}\,$ belongs to $[0,\pi/2]$, confirming that in the expression for $\Phi(x)$ we could indeed use the usual continuous determination of the square root.  \  Finally the expression of $\,\Phi\,$ shows indeed that it is analytical on $\R_+\,$, as the product of three analytical terms. $\;\diamond$ \parm   

      To handle the oscillatory integral yielding $ \II_\e(y,z)$ in Proposition \ref{pro.exprAvChgCont}, as is usual we shall change the path of integration. However, for that we need an analytic continuation, which demands more work. Part of it was already done in Section \ref{sec.calculFonctQuad}. The following is another step in this direction. \vspace{-2mm}  
\blem \label{lem.racinf} \  The function  $\,{\ds \frac{\rt1 x}{f(0,x)}\,}$ is meromorphic on $\,\C\moins\R_-\,$, and its poles are the values $\,4\rt1\theta_k^2\,$, $k\ge 1\,$, where $\{\theta_1<\theta_2<\ldots\}$ are the positive roots of $\,\tg \theta_k=\theta_k\,$. \  In particular, $\,\frac{4\pi}{3}< \theta_1< \frac{3\pi}{2}\,$\raise1pt\hbox{.} \   We have $\,f(0,x-\rt1\!\chi) = f(\chi,x)$. Moreover $\,{\ds \frac{\rt1 x}{f(0,x)}\,}$ is real on the segment $\big[0, 4\rt1\!\theta_1^2\big[\,$, positive on $\,\big]0, \rt1\!\pi^2\big[\,$ and negative on $\,\big]\rt1\!\pi^2, 4\rt1\!\theta_1^2\big[\,$.   \vspace{-2mm} 
\elem 
\ub{Proof} \   By (\ref{df.fnf}), we have \  ${ f(0,x) =  1 - {\frac{{\ts\th}\!\!\big[\sqrt{\rt1\!x}\big/2\,\big]}{\sqrt{\rt1\!x}\big/2}}\,}$, hence \  ${\ds f(0,-4\rt1\! x^2) =  1 - {\frac{\th x}{x}}\,}$\raise1.5pt\hbox{.} \parn 
Now the equation  \  $r\,e^{\rt1 \theta} =\th(r\,e^{\rt1 \theta})\,$ has no solution for $\,r>0$ and $\,| \theta| < {\pi}/{2}\,$\raise1.5pt\hbox{.} \parn   
Indeed, this equation is equivalent to both \  $r= \cos\theta\, \th\!(r \cos\theta) + \sin\theta\, \tg(r \sin\theta)$ and  
\  $r\, \th\!(r \cos\theta)\,\tg(r \sin\theta) = \cos\theta\, \tg(r \sin\theta) - \sin\theta\,  \th\!(r \cos\theta) $, which by eliminating $r$ implies \  $\tg\theta = \frac{\sin(2r \sin\theta)}{\sh(2r \cos\theta)}< \tg\theta\;$ for $\,0<\theta<\pi/2\,$, a contradiction. By symmetry, the claim holds for negative $\theta$ as well, and also clearly for $\theta=0$. \  On the contrary, this does not hold for $\,\theta=\pi/2\,$, since the imaginary roots are $\,\pm\rt1\!\theta_k\,$. 
\   Thus \  ${\ds \frac{\rt1\! (-4\rt1\! x^2)}{f(0,-4\rt1\! x^2)} = \frac{4\,x^3}{x-\th x}\;}$ is even holomorphic in $\{\,|\arg(x)|<\pi/2\}$, whereas  \  ${\ds \frac{\rt1\! (4\rt1\! x^2)}{f(0,4\rt1\! x^2)} = \frac{4\,x^3}{\tg x-x}\;}$ has poles. \parn  
By symmetry we can change $\theta$ into $\theta+\pi$,  which entails the analyticity for $\,\pi>\arg(x) >\pi/2\,$ as well. \  Finally, from the expression in Proposition \ref{pro.exprAvChgCont} we get 
$$ \frac{\rt1\! (4\rt1\! x^2)}{f(0,4\rt1\! x^2)} = \frac{4 x^3\, \cos{x}}{\sin{x} - {x}\, \cos{x}} \, \raise1.9pt\hbox{.} \;\;\diamond \vspace{-2mm}  $$  

\brem \label{rem.util2temps} \  {\rm  The dominant part of the integral in Proposition \ref{pro.exprAvChgCont} happens to belong to $\rt1\R$, and then does not contribute to the real projection. This means that the contribution to the heat kernel we have to extract is only a very small part of this integral. Therefore handling it remains delicate. We shall proceed by first changing the contour, in order to suppress the dominant hanger-on contribution. \   A good contour must go through a convenient saddle point of the considered oscillatory integral. However in the present case the saddle points are not easy to find\raise0pt\hbox{.}  
}\erem   \vspace{-2mm}

\subsection{First sub-case\,: $\, y\le0\,,\, z=0$} \label{sec.z=0>y1} \indf 
   We deal here with the sub-case $\,y\le 0\,$, in (recall Proposition \ref{pro.exprAvChgCont}) 
$$ \II_\e(y,0) =\, 2\; \Re\Bigg\{\! \int_{0}^\infty\! \exp\!\left[{ {\frac{\e - y}{\e^2}} }\rt1\! x \right] \Phi(x)\, dx \Bigg\} . $$   

   The dominant saddle point of the above oscillatory integral, happens to be close to $\, 4\rt1\!\! \left(\pi^2-\frac{\e^2}{4(\e-y)}\right)$, as will be confirmed later, by Proposition \ref{pro.ControlIntPhi}. \  The point $\rt1\!\pi^2$ \big(though a root of $\,\ch\!\sqrt{2x} - \cos\!\sqrt{2x}\,$ too, as mentioned in the proof of Proposition \ref{pro.exprAvChgCont}\big) is actually not a true pole, as will be clear in Lemma \ref{lem.LiftSqRtxPsi} below. \par  
   Thus we need to change the path of integration and to use an analytical continuation of the function $\,\Phi\,$ on the closure of a convenient domain. We already have part of this continuation, owing to Lemma \ref{lem.ProlAnalEspCompl}, Corollary \ref{cor.ValFourIntomt} and (\ref{f.E0tildx-ichi})(\ref{df.fnf}). \parn 
For the new contour, we need to change $\,x\,$ into $\,x-\rt1\!\chi\,$, \  with \  $\,x\ge 0\ge\chi \ge 4(\eta-\pi^2)$ \parn 
and $\,\eta = \frac{\e^2}{4(\e-y)}>0\,$\raise0.1pt\hbox{.} 
According to (\ref{f.E0tildx-ichi}) and Proposition \ref{pro.exprAvChgCont}, the function $\Phi$ admits the following continuation at least in  
$\,]-\pi^2, \infty[\,\times\, ]0,\infty[\,$: 
\begin{equation}  \label{df.PhiProlAnal}
\Phi\big(x-\rt1\!\chi\big) \equiv\, \Phi(\chi, x) := \, e^{\rt1 \f(\chi,x)} \times \psi(\chi,x) \times \sqrt{2\pi \,\frac{\chi + \rt1\! x}{f(\chi,x)}}\, \raise1.5pt\hbox{,}  \vspace{-2mm}
\end{equation}
where \vspace{-2mm} 
\begin{equation} \label{df.psiProlAnal}
\psi(\chi,x) := \left[{\frac{2\,\sqrt{\chi^2+x^2}}{\ch\!\sqrt{ 2\sqrt{\chi^2+x^2}+2\chi}  -\cos\!\sqrt{ 2\sqrt{\chi^2+x^2}-2\chi}}}\right]^{1/4}   
\end{equation}
$$ = \left[{\frac{a^2+b^2}{(\ch a-\cos b)(\ch a+\cos b)}}\right]^{1/4} \quad \hbox{with $\,a,b\,$ as in (\ref{f.racinexkhi}) and Lemma \ref{lem.LiftSquaRoot} below,} $$ 
$\f(\chi,x)$ is given by Lemma \ref{lem.ProlAnalEspCompl}, $\,f(\chi,x)$ is given by (\ref{df.fnf}) (in Corollary \ref{cor.ValFourIntomt}). 
   Now, in order that the above continuation be analytic, we need a continuous lift of the square root ${\ds\sqrt{\frac{\chi + \rt1\! x}{f(\chi,x)}}\, }\raise1pt\hbox{,}$ \   for $\,-\pi^2-\eta\le \chi \le 0\le x\,$. It is provided as follows.  
\blem  \label{lem.LiftSquaRoot} \  The square root $\sqrt{\frac{\chi + \rt1\! x}{f(\chi,x)}}\,$ admits the following  analytical lift, for any $\,\chi,x\,$ such that  $\,\chi > -4 \theta^2_1\,$ and  $\,x> 0\,$:  
{
$$ \sqrt{\frac{\chi + \rt1\! x}{f(\chi,x)}}\,= \, \frac{ (a^2+b^2)\,(\ch a + \cos b)^{1/4} \times e^{\rt1 \tilde\f(\chi,x)/2}} {\Big[(a^2+b^2)(\ch a + \cos b) - 4(a\, \sh a + b\, \sin b) + 4(\ch a - \cos b)\Big]^{1/4} }
$$ }
with the notation (\ref{f.racinexkhi})\,: \vspace{-2mm} 
$$ \sqrt{\chi+\rt1 x} = a+\rt1 b\quad \hbox{ where }\quad a:= \sqrt{\frac{\sqrt{\chi^2+x^2}+\chi}{2}}\,, \; b:= \sqrt{\frac{\sqrt{\chi^2+x^2}-\chi}{2}}\,\raise1.6pt\hbox{,} $$ 
and the (bounded for bounded $\chi$) argument\,: 
{\small 
$$ \tilde\f(\chi,x)\, =\, \arg\!\left(\frac{\rt1\! x}{f(0,x)}\right) + \frac{\pi}{2} - \frac{3}{2}\, \arctg\!\bigg[\frac{\chi}{x}\bigg] \,+  $$ 
$$ +  \int_0^\chi { \frac{(b\, \sh a - a\,\sin b)(a\,\sh a+b\,\sin b-\ch a + \cos b\big)}{\sqrt{\chi^2+x^2}\,(\ch a + \cos b)\big[(a^2+b^2)(\ch a + \cos b) - 4(a\, \sh a + b\, \sin b) + 4(\ch a - \cos b)\big]}}\, d\chi\, \raise1.9pt\hbox{.} $$ }  
\elem
\ub{Proof} \quad According to (\ref{df.fnf}), we have 
$$ \frac{1}{f(\chi,x)} = \frac{(a^2+b^2)\, (\ch a + \cos b)} {(a^2+b^2)\, (\ch a + \cos b)- 2a\, \sh a - 2b\, \sin b + 2\rt1 \! (b\,\sh a - a\, \sin b)} = \rr\, e^{\rt1\alpha} \vspace{1mm} \raise0.1pt\hbox{} $$ 
(which is real for $\,x=0\,$), \  with \ $\alpha\in \R\,$, \vspace{-1mm} 
$$ \rr := \,{\ts\frac{(a^2+b^2)\, (\ch a + \cos b)} {\sqrt{\left[(a^2+b^2)\, (\ch a + \cos b)- 2a\, \sh a - 2b\, \sin b \right]^2 + 4(b\,\sh a - a\, \sin b)^2} } }\,= \, \frac{\sqrt{a^2+b^2}\, (\ch a + \cos b)} {\sqrt{B} } \,\vspace{-0mm} \raise1.9pt\hbox{,} $$ 
$$ B := (a^2+b^2)(\ch a + \cos b)^2 - 4(a\, \sh a + b\, \sin b)(\ch a + \cos b) + 4(\sh\!^2 a + \sin^2\!b)\, $$
$$  = (\ch a + \cos b)\times\big[(a^2+b^2)(\ch a + \cos b) - 4(a\, \sh a + b\, \sin b) + 4(\ch a - \cos b)\big] , $$
and 
$$ \cos \alpha = \frac{(a^2+b^2)^{1/2}\, (\ch a + \cos b)- 2(a^2+b^2)^{-1/2}\,(a\, \sh a +b\, \sin b)}{\sqrt{B} } \,\raise1.9pt\hbox{,} $$ 
$$ \sin \alpha = \frac{2(a^2+b^2)^{-1/2}\,(a\,\sin b - b\, \sh a)}{\sqrt{B} } \,\raise1.9pt\hbox{.} \vspace{0mm} $$ 
Note that by Lemma \ref{lem.racinf}, $\,B\not=0\,$ for $\,\chi > - 4 \theta_1^2\,$, hence for $\,\chi > - 7 \pi^2$. \parn   
Thence, we successively have\,: 
$$ \frac{\partial B}{\partial a} =  2 (\ch a + \cos b) \big[ a(\cos b-\ch a) + (a^2+b^2-2) \,\sh a\big] - 4(a\, \sh a + b\, \sin b - 2\,\ch a)\, \sh a\, ;  $$ 
$$ \frac{\partial B}{\partial b} = 2 (\ch a + \cos b) \big[ b(\ch a-\cos b) - (a^2+b^2+2)\sin b \big] + 4(a\, \sh a + b\, \sin b +2\,\cos b) \sin b \,;  $$ 
$$ \frac{1}{2}\! \left[a\,\frac{\partial B}{\partial a}-b\,\frac{\partial B}{\partial b}\right] = \big[a\,\sh a+b\,\sin b-\ch a + \cos b\big] \big[(a^2+b^2)(\ch a + \cos b)-2(a\,\sh a+b\,\sin b)\big]\, ; $$ 
$$ \frac{1}{2} \cos \alpha \,\frac{\partial\alpha}{\partial a} = \frac{\partial\sin \alpha}{2\,\partial a} = \, \frac{\partial}{\partial a}\,  \frac{(a^2+b^2)^{-1/2}\,(a\,\sin b - b\, \sh a)}{\sqrt{B} }\, =\, \frac{A}{B^{3/2}\, (a^2+b^2)^{3/2}}\, \raise0.9pt\hbox{;} $$ 
$$ \frac{1}{2} \cos \alpha \,\frac{\partial\alpha}{\partial b} = \frac{\partial\sin \alpha}{2\,\partial b} = \, \frac{\partial}{\partial b}\,  \frac{(a^2+b^2)^{-1/2}\,(a\,\sin b - b\, \sh a)}{\sqrt{B} }\, =\, \frac{A'}{B^{3/2}\, (a^2+b^2)^{3/2}}\, \raise0.9pt\hbox{;} $$ 
with \vspace{-2mm} 
$$ A := \big(a\,\sh a + b\,\sin b -(a^2+b^2)\,\ch a\big) b B + \frac{1}{2} (a^2+b^2)(b\, \sh a - a\,\sin b)\, \frac{\partial B}{\partial a} \, ; $$ 
$$ A' := \big((a^2+b^2)\,\cos b - a\,\sh a - b\,\sin b\big) a B + \frac{1}{2} (a^2+b^2)(b\, \sh a - a\,\sin b)\, \frac{\partial B}{\partial b} \; ; $$ 
$$ a\,A - b\, A' =  \big[2(a\,\sh a+b\,\sin b) - (a^2+b^2)(\ch a + \cos b)\big]\, ab\,B \vspace{-1mm} $$
$$ \hskip22mm  +\, \frac{1}{2}\, (a^2+b^2)(b\, \sh a - a\,\sin b) \left[a\,\frac{\partial B}{\partial a} - b\,\frac{\partial B}{\partial b}\right]  $$ 
$$ = \, \big[(a^2+b^2)(\ch a + \cos b)-2(a\,\sh a+b\,\sin b)\big] \times \vspace{-1mm} $$
$$ \times  \Big( (a^2+b^2)(b\, \sh a - a\,\sin b)\big[a\,\sh a+b\,\sin b-\ch a + \cos b\big] - ab\,B\Big) . $$
Hence
$$ \frac{\partial}{\partial \chi}\, =\, \frac{\partial b}{\partial \chi}\,\frac{\partial}{\partial b}\,+\frac{\partial a}{\partial \chi}\,\frac{\partial}{\partial a}\, =\, \frac{1}{2(a^2+b^2)}\left(a\,\frac{\partial}{\partial a} - b\,\frac{\partial}{\partial b}\right) ; 
\quad  \cos \alpha \,\frac{\partial\alpha}{\partial \chi}\, =\,  \frac{a\,A - b\, A'}{ B^{3/2}\, (a^2+b^2)^{5/2}}\; ; \raise1.9pt\hbox{}   $$
$$ \frac{\partial\alpha}{\partial \chi}\, =\,  \frac{(b\, \sh a - a\,\sin b)\big[a\,\sh a+b\,\sin b-\ch a + \cos b\big]}{\sqrt{\chi^2+x^2}\,\, B} - \frac{x}{2(\chi^2+x^2)}\, \raise1.9pt\hbox{,} $$ 
since $\,(a^2+b^2) = \sqrt{\chi^2+x^2}\,$ and $\,x=2ab\,$. \    Then 
$$ \frac{\chi + \rt1\! x}{f(\chi,x)} = (\chi + \rt1\! x)\, \rr\, e^{\rt1\alpha} = (a^2+b^2)\, \rr\, e^{\rt1\tilde\f(\chi,x)} \vspace{-1.5mm} $$ 
with \quad ${\ds \tilde\f(\chi,x) = \alpha + \pi\,1_{\{\chi\le 0\}} + \arctg(x/\chi) \, =\, \alpha + \frac{\pi}{2} - \int_0^\chi  \frac{x \; d\chi}{\chi ^2+ x^2}\,}$\raise2pt\hbox{.} \parn 
Finally we obtain \vspace{-2mm} 
$$ \tilde\f(\chi,x)\, =\, \frac{\pi}{2} - \int_0^\chi  \frac{x \; d\chi}{\chi ^2+ x^2} + \arg\!\left[\frac{\rt1\! x}{f(0,x)}\right] + \int_0^\chi  \frac{\partial\alpha}{\partial \chi}\, d\chi\, $$ 
$$ =\,\frac{\pi}{2} - \frac{3}{2}\, \arctg\!\bigg[\frac{\chi}{x}\bigg] + \arg\!\left[\frac{\rt1\! x}{f(0,x)}\right] + \int_0^\chi { \frac{(b\, \sh a - a\,\sin b)(a\,\sh a+b\,\sin b-\ch a + \cos b)}{\sqrt{\chi^2+x^2}\; B}}\, d\chi\, \raise1.9pt\hbox{.} $$  
We already noticed that $\,0\le \arg\!\left[\frac{\rt1\! x}{f(0,x)}\right]\! \le\frac{\pi}{2}\,$\raise1.1pt\hbox{.} \  Moreover, as $\,a,b\,$ grow with $\,x\,$, asymptotically as $\sqrt{x}\,$, the above integral is bounded too, provided $\,\chi\,$ remains bounded. $\;\diamond$  \parmn   

   Lemma \ref{lem.LiftSquaRoot}, (\ref{df.PhiProlAnal}) and the expression following (\ref{df.psiProlAnal}) (i.e., $\,\psi\,$ in terms of $\,a,b$) together easily yield the following. 
\blem  \label{lem.LiftSqRtxPsi} \  For any real $\,\chi\,, x\,$ such that  
$\,x> 0\,$ we have the following analytical lift\,:  
$$ \Phi(\chi, x) =  \frac{ (a^2+b^2)^{5/4}\, \left(\ch a - \cos b\right)^{-1/4} \times  e^{\rt1\! (\f(\chi,x)+\tilde\f(\chi,x)/2)}} {\Big[(a^2+b^2)(\ch a + \cos b) - 4(a\, \sh a + b\, \sin b) + 4(\ch a - \cos b)\Big]^{1/4} } \, \raise1.9pt\hbox{,} $$  
with the notation (\ref{f.racinexkhi}) as in Lemma \ref{lem.LiftSquaRoot}, $\f(\chi,x)$ given by Lemma \ref{lem.ProlAnalEspCompl}, and the same argument $\,\tilde\f(\chi,x)$ as in Lemma \ref{lem.LiftSquaRoot}.  \par 
Moreover the above extends continuously to $\,x=0\,$, except for the singular points where the denominator vanishes, i.e., $\,\chi= -4 \theta_k^2\,$ ($k\ge 1$, recall Lemma \ref{lem.racinf}), and $\,\chi= -4 k^2\pi^2\,$ ($k\ge 1$, namely the zeros of $(\ch a - \cos b)$).
\elem  \par   

\if{  
   We shall also need the following, in order to control the contribution near the singular point $\rt1\!\pi^2$. \vspace{-2mm} 
\blem  \label{lem.ControlPhi(ipi2)}  \   On $\big\{0< \eta\le 1\,,\,|\theta|\le \pi/2\big\}$,  \  $\left| \Phi\big(-\pi^2-\eta\, \sin\theta\,, \eta\,\cos\theta\big)\right|\,$ is bounded. 
\elem
\ub{Proof} \quad   According to (\ref{df.psiProlAnal}) we have \vspace{-1mm} 
$$ \psi\big(-\pi^2-\eta\, \sin\theta\,, \eta\,\cos\theta\big) = \vspace{-2mm} $$ 
{\small %
$$ \left[{\frac{2\,\sqrt{\pi^4 + 2\pi^2\, \eta\,\sin\theta +\eta^2}}{\ch\!\sqrt{ 2\sqrt{\pi^4\! + 2\pi^2 \eta\sin\theta +\eta^2} - 2\pi^2\! - 2\eta \sin\theta} - \cos\!\sqrt{ 2\sqrt{\pi^4\! + 2\pi^2 \eta\sin\theta +\eta^2} + 2\pi^2\! + 2\eta \sin\theta}}}\right]^{1/4} $$ } 
$$ = \left[{\frac{2\pi^2 + 2\eta\,\sin\theta +\pi\2\,\eta^2\cos^2\!\theta +\O(\eta^3)} {\ch\!\sqrt{ \pi\2\,\eta^2\cos^2\!\theta +\O(\eta^3)} - \cos\!\sqrt{4\pi^2 + 4\eta\,\sin\theta +\pi\2\,\eta^2\cos^2\!\theta +\O(\eta^3)}}}\right]^{1/4} $$ 
$$ = \left[{\frac{2\pi^2 + 2\eta\,\sin\theta +\O(\eta^2)} {1+\frac{\eta^2\cos^2\!\theta +\O(\eta^4)}{2\,\pi^2} - \cos\!\left[\frac{\eta \sin\theta }{\pi} + \frac{\eta^2\cos(2\theta) +\O(\eta^3)}{4\,\pi^3}\right] }}\right]^{1/4} ,$$
i.e., 
$$ \psi\big(-\pi^2-\eta\, \sin\theta\,, \eta\,\cos\theta\big) =\, \pi\sqrt{2/\eta} \times \Big(1 + \frac{\eta\,\sin\theta}{4\,\pi^2} +\O(\eta^2)\Big) . $$ 
Then we use Lemma \ref{lem.LiftSquaRoot} and its proof  to estimate the modulus \parn
${\ds \left|\sqrt{\frac{ \rt1\! \eta\,\cos\theta - \pi^2-\eta\, \sin\theta}{f\big(-\pi^2-\eta\, \sin\theta\,, \eta\,\cos\theta\big)}}\,\right| = (\pi+\O(\eta))\sqrt{\rr}\,}$\raise0pt\hbox{,} \quad  where 
$$ \rr  = \,{\frac{(a^2+b^2)\, (\ch a + \cos b)} {\sqrt{\left[(a^2+b^2)\, (\ch a + \cos b)- 2a\, \sh a - 2b\, \sin b \right]^2 + 4(b\,\sh a - a\, \sin b)^2} } }\,\raise1.9pt\hbox{.} $$ 
Now, as for $ \psi\big(-\pi^2-\eta\, \sin\theta\,, \eta\,\cos\theta\big)$ above, we have \  ${\ds a = \,\frac{\eta\,\cos\theta}{2\,\pi} +\O(\eta^2) = \sh a \,}$ \  and 
$$ b = \sqrt{\pi^2 + \eta\,\sin\theta + \frac{\eta^2\cos^2\!\theta +\O(\eta^3)}{4\,\pi^2}}\, = \pi + \frac{\eta \sin\theta }{2\pi} + \frac{\eta^2\cos(2\theta) +\O(\eta^3)}{8\,\pi^3}\,\raise1.6pt\hbox{,} $$ 
$$ \ch a = 1+\frac{\eta^2\cos^2\!\theta +\O(\eta^4)}{8\,\pi^2} \,\raise1.4pt\hbox{,}  \  \cos b = -1 + \frac{\eta^2\sin^2\!\theta +\O(\eta^4)}{8\,\pi^2} \,\raise1.4pt\hbox{,} $$
$$ (a^2+b^2)\, (\ch a + \cos b) = \frac{\eta^2 }{8}+ \frac{\eta^3 \sin\theta}{8\,\pi^2}  +\O(\eta^4)\, , \  \;  \sin b = -\frac{\eta\, \sin\theta }{2\pi} - \frac{\eta^2\cos(2\theta) +\O(\eta^3)}{8\,\pi^3}\,\raise1.6pt\hbox{,} $$
$$ a\, \sh a = \,\frac{\eta^2\cos^2\!\theta}{4\,\pi^2} +\O(\eta^3)\,\raise1.6pt\hbox{,} \quad  b\, \sh a = \,\frac{\eta\,\cos\theta}{2} +\O(\eta^2)\,\raise1.6pt\hbox{,} \quad a\,\sin b = -\frac{\eta^2\cos\theta\, \sin\theta }{2\pi} +\O(\eta^3)\,\hbox{,} $$
$$ b\,\sin b = -\frac{\eta\, \sin\theta }{2} +\O(\eta^2)\,\hbox{,} \quad  \hbox{ whence } $$
$$ \rr = \,{\frac{\frac{\eta^2 }{8}+ \frac{\eta^3 \sin\theta}{8\,\pi^2}  +\O(\eta^4)} {\sqrt{\left[\frac{\eta^2 }{8} + \O(\eta^3)- \frac{\eta^2\cos^2\!\theta}{2\,\pi^2} - \eta\, \sin\theta\right]^2 + 4\left(\frac{\eta\,\cos\theta}{2} +\O(\eta^2)\right)^2} } }\,=\, \frac{\eta}{8}+ \O(\eta^2)\,\raise0pt\hbox{.} $$ 
The claim follows now from  the expression (\ref{df.PhiProlAnal}) defining $\Phi(\chi, x)$. $\diamond$ \parb 
}\fi 
\if{ 
\subsubsection{tentative en 2 changements successifs} \indf 
  Let us consider the contour 
\begin{equation} \label{df.nContourbis} 
\Gamma_\eta^1 := \Big[0\,, \rt1\!\pi^2\Big]\, \, \bigcup\, \Big(\rt1\pi^2 + \R_+\Big) .
\end{equation}
Lemma \ref{lem.LiftSqRtxPsi} allows to change the contour  in \ 
 ${\ds  \Re\Bigg\{\! \int_{0}^\infty\! \exp\!\left[{ {\frac{y-\e}{\e^2}} }\rt1\! x \right] \overline{\Phi}(x)\, dx \Bigg\}}$,  into $\,\Gamma_\eta^1\,$ of (\ref{df.nContourb}), in order to obtain the first step of the reduction announced in Remark \ref{rem.util2temps}, namely the following.  
\bpro \label{pro.ChangContz=0>yb1} \  For any $\;y\le0\,$ and small $\,\e >0\,$ we have 
$$ \II_\e(y,0) = \, \exp\!\left[{ {\frac{\e - y}{\e^2}} }\,\pi^2\right]\times 2\; \Re\Bigg\{\! \int_{0}^\infty\! \exp\!\left[{ {\frac{\e - y}{\e^2}} }\rt1\! x \right]  \Phi\big(-\pi^2, x\big)\, dx \Bigg\} \, ,  $$ 
where \  $\Phi(\chi, x)  \equiv \Phi\big(x-\rt1\!\chi\big)\,$ was defined in (\ref{df.PhiProlAnal}). 
\epro   
\ub{Proof} \quad The change of contour (\ref{df.nContourbis}) is justified by the two following observations. \parn
On the one hand, restricting to the range $\,-\pi^2\le \chi\le 0\,$, i.e., to $\,0\le b\le \pi\,$, 

and on the other hand, by 
$$ \limsup_{R\to\infty} \left|\int_{0}^{\pi^2}\! \exp\!\left[{ {\frac{\e-y}{\e^2}} }\,(x -\rt1\!R)\right] \overline{\Phi}\big(-x,R\big)\, dx\,\right|$$
$$ \le\,  \limsup_{R\to\infty}\, \sup_{-\pi^2\le \chi\le 0}\, \left| {\Phi}\big(\chi,R\big)\right| \times \int^{\pi^2}_{0} \exp\!\left[{ {\frac{\e-y}{\e^2}} }\,x \right] dx\, =\, 0\,, $$ 
since Lemma \ref{lem.LiftSqRtxPsi} ensures that $\,\left| \Phi\big(\chi,R\big)\right| = \O\big(R\,e^{-\sqrt{R/8}}\,\big)$, uniformly for bounded $\chi\,$. \parn 
Thus performing this change yields\,: \  {\bf NON !!  $\,\overline{\Phi}\,$ n'est pas analytique !! }
$$ \int_{0}^\infty \exp\!\left[{ {\frac{y-\e}{\e^2}} }\rt1\! x \right] \overline{\Phi}(x)\, dx \, =  \int_{\Gamma_\eta^1} \exp\!\left[{ {\frac{y-\e}{\e^2}} }\rt1\! x \right] \overline{\Phi}(x)\, dx $$
$$ = \rt1\!\! \int_{0}^{\pi^2} \exp\!\left[{ {\frac{\e-y}{\e^2}} }\,x \right] \overline{\Phi}(-x,0)\, dx\, +\, e^{{ {\frac{\e-y}{\e^2}} }\,\pi^2}\! \int_{0}^\infty \exp\!\left[{ {\frac{y-\e}{\e^2}} }\rt1\! x \right]  \overline{\Phi}\big(- \pi^2, x\big)\, dx\, .  $$ \parn 
Then by Lemma \ref{lem.ProlAnalEspCompl}, for $\,\chi\le 0\,$, letting $x\sea 0\,$, on the one hand we find\,: \parsn 
\centerline{$\f(\chi,0)=0\,$ separately for $\,b=\sqrt{-\chi}\, \in [0,\pi/2]\,$ and for $\,b=\sqrt{-\chi}\,\in\,]\pi/2, \pi]$,} \parsn  
and on the other hand, we saw in Lemma \ref{lem.racinf} that $\,\frac{-x}{f(-x,0)}\in\R_+\,$ for $\,0\le x\le \pi^2$. \quad Hence, according to (\ref{df.PhiProlAnal}), for all   $\,x\in[0,\pi^2]\,$ we have \vspace{-2mm} 
$$ \Phi(-x,0) = e^{\rt1 \f(-x,0)} \,\psi(-x,0)\sqrt{\ts\frac{- 2\pi \,x}{f(-x,0)}}\,\in\R\, \raise0pt\hbox{,} \vspace{-2mm} $$ 
so that \vspace{-2mm} 
$$ \rt1\!\! \int_{0}^{\pi^2} \exp\!\left[{\ts {\frac{y-\e}{\e^2}} }\,x \right] \Phi(-x,0)\, dx\, \in\rt1 \R\, $$ 
and then disappears when taking the real part.  \  This yields 
$$ \II_\e(y,0) =\, 2\; \Re\Bigg\{\! \int_{0}^\infty\! \exp\!\left[{ {\frac{\e - y}{\e^2}} }\rt1\! x \right] \Phi(x)\, dx\! \Bigg\}  =\, 2\; \Re\Bigg\{\! \int_{0}^\infty\! \exp\!\left[{ {\frac{y-\e}{\e^2}} }\rt1\! x \right] \overline{\Phi}(x)\, dx\! \Bigg\}  $$
$$ =\,  \exp\!\left[{ {\frac{\e - y}{\e^2}} }\,\pi^2\right]\times 2\; \Re\Bigg\{\! \int_{0}^\infty\! \exp\!\left[{ {\frac{y-\e}{\e^2}} }\rt1\! x \right]  \overline{\Phi}\big(- \pi^2, x\big)\, dx\! \Bigg\}  $$
$$ =\,  \exp\!\left[{ {\frac{\e - y}{\e^2}} }\,\pi^2\right]\times 2\; \Re\Bigg\{\! \int_{0}^\infty\! \exp\!\left[{ {\frac{\e-y}{\e^2}} }\rt1\! x \right] {\Phi}\big(- \pi^2, x\big)\, dx\! \Bigg\} . \;\; \diamond $$
}\fi 

   Let us consider the contour 
\begin{equation} \label{df.nContourb} 
\Gamma_\eta := \Big[0\,, 4\rt1\!(\pi^2-\eta)\Big]\, \, \bigcup\, \Big(4\rt1\!(\pi^2-\eta) + \R_+\Big) .
\end{equation}
Lemma \ref{lem.LiftSqRtxPsi} allows to change the contour  in \ 
 ${\ds  \Re\Bigg\{\! \int_{0}^\infty\! \exp\!\left[{ {\frac{\e - y}{\e^2}} }\rt1\! x \right] \Phi(x)\, dx \Bigg\}}$,  into $\,\Gamma_\eta\,$ of (\ref{df.nContourb}), in order to obtain the first reduction announced in Remark \ref{rem.util2temps}, namely the following.  
\bpro \label{pro.ChangContz=0>yb} \  For any $\;y\le0\,$ and small $\,\e,\, \eta>0\,$ we have 
$$ \II_\e(y,0) = \, \exp\!\left[{ {\frac{y-\e}{\e^2}} }\,4(\pi^2-\eta)\right]\times 2\; \Re\Bigg\{\! \int_{0}^\infty\! \exp\!\left[{ {\frac{\e - y}{\e^2}} }\rt1\! x \right]  \Phi\big(4\eta - 4\pi^2, x\big)\, dx \Bigg\} \, ,  $$ 
where \  $\Phi(\chi, x)  \equiv \Phi\big(x-\rt1\!\chi\big)\,$ was defined in (\ref{df.PhiProlAnal}). 
\epro   
\ub{Proof} \quad The above-mentioned change of contour (\ref{df.nContourb}) is justified by 
$$ \limsup_{R\to\infty} \left|\int_{0}^{4\pi^2-\eta}\! \exp\!\left[{ {\frac{y-\e}{\e^2}} }\,(x -\rt1\!R)\right] \Phi\big(-x,R\big)\, dx\,\right|\,\le  \limsup_{R\to\infty} \int_{4\eta-4\pi^2}^{0}\! \left| \Phi\big(\chi,R\big)\right| d\chi\, , $$ 
as Lemma \ref{lem.LiftSqRtxPsi} ensures that the latter vanishes\,: actually it shows that $\,\left| \Phi\big(\chi,R\big)\right| = \O\big(R\,e^{-\sqrt{R/8}\,}\big)$, uniformly for bounded $\chi\,$.  Thus performing this change yields\,: 
$$ \int_{0}^\infty \exp\!\left[{ {\frac{\e - y}{\e^2}} }\rt1\! x \right] \Phi(x)\, dx \, =  \int_{\Gamma_\eta} \exp\!\left[{ {\frac{\e - y}{\e^2}} }\rt1\! x \right] \Phi(x)\, dx $$
$$ = \rt1\!\! \int_{0}^{4\pi^2-\eta} \exp\!\left[{ {\frac{y-\e}{\e^2}} }\,x \right] \Phi(-x,0)\, dx\, $$
$$ +\, \exp\!\left[{ {\frac{y-\e}{\e^2}} }\,4(\pi^2-\eta)\right]\times \int_{0}^\infty \exp\!\left[{ {\frac{\e - y}{\e^2}} }\rt1\! x \right]  \Phi\big(4\eta - 4 \pi^2, x\big)\, dx\, .  $$ 

Then by Lemma \ref{lem.ProlAnalEspCompl}, for $\,\chi<0\,$, letting $x\sea 0\,$, we find\,: \parsn 
$*$  \  $\,\f(\chi,0)=0\,$ separately for $\,b=\sqrt{-\chi}\, \in [0,\pi/2]\,$ and for $\,b=\sqrt{-\chi}\,\in\,]\pi/2, \pi[\,$ ; \parn 
$*$  \  $\,\f(\chi,0)=-\frac{\pi}{2}\,$ separately for $\,b=\sqrt{-\chi}\, \in\,]\pi,3\pi/2[\,$ and for $\,b=\sqrt{-\chi}\,\in\,]3\pi/2, 2\pi[\,$. \parn 
On the other hand, we saw in Lemma \ref{lem.racinf} that $\,\frac{-x}{f(-x,0)}\in\R_+\,$ for $\,0\le x\le \pi^2$ and $\,\frac{-x}{f(-x,0)}\in\R_-\,$ for $\,\pi^2\le x< 4\theta_1^2\,$. \quad Hence, according to (\ref{df.PhiProlAnal}), for all   $\,x\in \big[0,4\theta_1^2\big[\,$ we have \vspace{-2mm} 
$$ \Phi(-x,0) = e^{\rt1 \f(-x,0)} \,\psi(-x,0)\sqrt{\ts\frac{- 2\pi \,x}{f(-x,0)}}\,\in\R\, \raise0pt\hbox{,} \vspace{-2mm} $$ 
so that \vspace{-2mm} 
$$ \rt1\!\! \int_{0}^{4\pi^2-\eta} \exp\!\left[{\ts {\frac{y-\e}{\e^2}} }\,x \right] \Phi(-x,0)\, dx\, \in\rt1 \R\, $$ 
and then disappears when taking the real part.  $\diamond$ \parb 

%

   We have now to carefully analyze the integral appearing in Proposition \ref{pro.ChangContz=0>yb}, and then  
first to approach the term $\,\Phi\big(4\eta - 4\pi^2, x\big)$ near the singularity $4\rt1\!\pi^2$. This is the step  where the choice of $\,4\eta - 4\pi^2$ begins to prove to be the right one. We take \  $\,{\ds \eta = \frac{\e^2}{4(\e-y)}\,}$\raise1.9pt\hbox{.}
\blem \label{lem.AnalPhiprSing} \  For $\,y\le 0\,$, $\,r':=\frac{3}{2} -\frac{3}{4}\,1_{\{y=0\}} \le r \le \frac{9}{10}(2-1_{\{y=0\}})$ and $\,0\le x\le\e^{r-2}(\e-y)\,$, with real $\O()$ and uniformly with respect to $\,y\le 0$, \   we have 
$$ \Phi\left({\frac{\e^2}{\e-y}} - 4\pi^2, {\frac{\e^2}{\e - y}}\,x\right) =\, \frac{ \big(8\pi^{5/2} + \O({\frac{\e^2}{\e-y}})\big)\, \sqrt{\e-y}} {(x^2+1)^{1/4}\; \e}\,\times e^{\rt1 \left(-\frac{3\pi}{8}- \frac{1}{2}\,\arctg x + \O(\e^{r/3}) \right)}\, . \vspace{-2mm} $$ 
\elem 
\ub{Proof} \quad Using the notation (\ref{f.racinexkhi})\,:  $\sqrt{\rt1 x- 4(\pi^2-\eta)}\, =\, a+\rt1 b\,$, \  for small  $\,x >\eta\,$ we successively have\,: \vspace{-2mm} 
{
$$ \sqrt{\chi^2 + x^2}\, =\, a^2+b^2 =\, 4\pi^2 -4\eta + \frac{x^2}{8\pi^2} + \O(\eta\,x^2)\, ;  \quad  a\, = \sqrt{ \frac{x^2}{16\pi^2} +\O(\eta\, x^2)}\, = \frac{x}{4\pi}+ \O(\eta\,x)\; ; \vspace{-1mm}  $$ } 
$$  b\, = \sqrt{4\pi^2-4\eta + \frac{x^2}{16\pi^2} +\O(\eta\, x^2)}\, =\, 2\pi - \frac{\eta}{\pi} + {\frac{x^2 - 16\eta^2}{64\,\pi^3}} + \O(\eta\, x^2 \if{+\eta^3}\fi ) \; ; $$
$$ \sh a\, = \frac{x}{4\pi}+ \O(\eta\,x)\, =\, \th a\;  ; \quad   \ch a\, =\, 1+ \frac{x^2}{32\pi^2} + \O(\eta\, x^2)\; ;  $$
$$ \sin b\, =\, - \frac{\eta}{\pi} + {\frac{x^2 - 16\eta^2}{64\,\pi^3}} + \O(\eta\, x^2 \if{+\eta^3}\fi ) \; ; \quad \cos b\, =\, 1- {\frac{\eta^2}{2\pi^2}} + \O(\eta\, x^2) \; ; $$ 
whence \vspace{-3mm} 
$$ \ch a +  \cos b\, =\, 2+ \O(\if{\eta^2+}\fi x^2)\; ;  \quad \ch a - \cos b\, =\, {\frac{x^2 + 16\,\eta^2}{32\,\pi^2}} + \O(\eta\,x^2)\; ; $$
$$ a\,\sh a + b\,\sin b\, = -2\eta+  {\ts\frac{9x^2 +16\eta^2}{32\,\pi^2}}+ \O(\eta\,x^2 \if{+\eta^3}\fi ) . 
\vspace{1mm} $$ 
Therefore, according to Lemma \ref{lem.LiftSqRtxPsi} we obtain\,: \vspace{-1mm} 
$$ \Phi(\chi, x)\, =\, \frac{ \big(4\pi^2 + \O(\eta+x^2)\big)^{5/4}\, \sqrt{8\pi}\, \times  e^{\rt1\! \left(\f(\chi,x)+\tilde\f(\chi,x)/2\right)}} {\big[ 16\pi^2+8\eta+\O(\if{\eta^2+}\fi x^2)\big]^{1/4}\times \big(x^2+16\,\eta^2\big)^{1/4}  } $$
$$ =\, \frac{ 8\pi^{5/2} + \O(\eta+x^2)} {\big(x^2+16\,\eta^2\big)^{1/4}}\times  e^{\rt1\! \left(\f(4\eta - 4\pi^2,\,x)+\tilde\f(4\eta - 4\pi^2,\,x)/2\right)} \, \raise0pt\hbox{} $$  
on the one hand, and on the other hand\,: \vspace{-2mm} 
$$ \tilde\f\!\left(4\eta - 4\pi^2,x\right) =\, \frac{\pi}{4}- \frac{3x}{8\pi^2} + \O(\if{\eta^2+}\fi x^2) + \int_0^{4\eta - 4\pi^2} {\ts\frac{(b\, \sh a - a\,\sin b)(a\,\sh a+b\,\sin b-\ch a + \cos b)}{\sqrt{\chi^2+x^2}\; B}}\, d\chi\, , $$ 
and owing to Lemma \ref{lem.ProlAnalEspCompl}\,:  \vspace{-1mm} 
$$ \f\!\left(4\eta - 4\pi^2,x\right) =\,\5\,\arctg(b/a)-\pi - \5\, \arctg\!\big[\tg(b-2\pi)\times \coth a\big] $$
$$ = \,-\5\,\arctg(a/b)-\pi + \5\, \arctg\!\big[\cotg(b-2\pi)\times \th a\big] $$ 
\begin{equation} \label{f.approxphi} 
=\, -\pi - \frac{x}{16\pi^2} + \O(\eta\,x) -\,\frac{1}{2}\,\arctg\!\left(\frac{x+ \O(\eta\,x)}{4\eta +\frac{16\eta^2 -x^2}{16\pi^2} +\O(\if{\eta^3+}\fi \eta\,x^2)}\right) . \vspace{1mm}
\end{equation} 
Then with \  $\,\eta = \frac{\e^2}{4(\e-y)}\,$, for  $\,\frac{3}{4} (2-1_{\{y=0\}}) \le r \le  \frac{9}{10}(2-1_{\{y=0\}})$ we have \  $\e^{r-2}(\e-y)\,\to \infty\,$ as $\,\e\sea 0\,$, \   and $\,\e^{2r} \le \eta^{3/2}\,$. \  Thus for $\,0\le x\le\e^{r-2}(\e-y)$, by the above we have\,: \vspace{-1mm} 
$$ \Phi\left({\frac{\e^2}{\e-y}} - 4\pi^2, {\frac{\e^2}{\e - y}}\,x\right)  =\,\frac{\sqrt{\e-y}}{\e}\times  \frac{ 8\pi^{5/2} + \O\big(\frac{\e^2}{\e-y}\big)} {\big(x^2+1\big)^{1/4}}  \vspace{-2mm} $$
$$ \times\,   \exp\!\left[\rt1\! \left(\f\big({\ts\frac{\e^2}{\e-y}} - 4\pi^2,\,{\ts\frac{\e^2}{\e-y}}\,x\big) + \5\, \tilde\f\big({\ts\frac{\e^2}{\e-y}} - 4\pi^2,\,{\ts\frac{\e^2}{\e-y}}\,x\big)\right)\right]  , \vspace{1mm} $$ 
with by (\ref{f.approxphi}) \big(using $\,\eta \ge \e^{4r/3}\,$, hence $\,x^2\eta\le \e^{2r}/\eta\le \e^{2r/3}\,$\big)\,:  \vspace{-1mm} 
$$ \f\left({\frac{\e^2}{\e-y}} - 4\pi^2\, \raise0.6pt\hbox{,}\, {\frac{\e^2}{\e-y}}\,x\right) = \O(\e^r)  -\pi - \frac{1}{2}\,\arctg\!\left[\frac{x\,[1+\O(\eta)]}{1+\O(x^2\eta) } \right]  $$ 
$$ = \O(\e^r) -\pi  - \,\5\,\arctg\!\big(x\,\big[1+(1+x^2)\, \O(\eta)\big]\big) = \O(\e^r) -\pi  - \,\5\,\arctg\!\big(x\,\big[1+\O(\e^{2r/3})\big]\big) $$
$$ = \O(\e^r) -\pi  - \,\5\,\arctg x + \5\,\arctg\!\left(\frac{x\,\O(\e^{2r/3})}{1+x^2\,\O(\e^{2r/3})}\right)  
=\, -\pi - \,\5\,\arctg x + \O(\e^{r/3}) \, . $$
Moreover, by (\ref{f.partReelleb}) near 0 we have\,: \vspace{-1mm} 
\begin{equation}   \label{f.DL0} 
\frac{\rt1 x}{f(0,x)} \, = \,6 + \rt1\frac{6}{5}\,x - \frac{x^2}{140} + \O(x^3)\, ,  
\end{equation}
and then by Lemma \ref{lem.LiftSquaRoot}, again for $\,0\le x\le\e^{r-2}(\e-y)$\,: 
$$ \tilde\f\big({\ts\frac{\e^2}{\e-y}} - 4\pi^2,\,{\ts\frac{\e^2}{\e-y}}\,x\big) $$
$$ =\, \arg\!\left[{\frac{\rt1\! x}{f(0,x)}}\right]  + {\frac{5\pi}{4}} - {\frac{3}{2}}\, \arctg\!\left[{\frac{\O(\e^r) }{4\pi^2}}\right] + \int_0^{\frac{\e^2}{\e-y}- 4\pi^2} {\ts \frac{(b\, \sh a - a\,\sin b)(a\,\sh a+b\,\sin b-\ch a + \cos b)}{\sqrt{\chi^2+x^2}\; B}}\, d\chi $$  
$$ =\, \frac{5\pi}{4} + \O(\e^r) + \O(1) \int_0^{\frac{\e^2}{\e-y}- 4\pi^2} {\frac{(b\, \sh a - a\,\sin b)(a\,\sh a+b\,\sin b-\ch a + \cos b)}{\sqrt{\chi^2+x^2}}}\, d\chi\, , $$
$$ \hbox{with (for $\,0\ge \chi> -4\pi^2$)  } \qquad a = \sqrt{\ts\frac{\sqrt{\chi^2+\O(\e^{2r})}+\chi}{2}}\,\, , \quad  b = \sqrt{\ts\frac{\sqrt{\chi^2+\O(\e^{2r})}-\chi}{2}}\,\raise1.6pt\hbox{.} $$ 
Let us distinguish in the above integral between $\,4\pi^2> |\chi|\ge \e^{r/2}$ and $\,0\le |\chi|< \e^{r/2}$. \parn
In the former case we have 
$$ a = \O(\e^{r}/|\chi|)\, , \  b = \sqrt{|\chi|} + \O\big(\e^{2r}/|\chi|^{3/2}\big) = \sqrt{|\chi|} + \O(\e^{r})\, , \  \sqrt{\chi^2+x^2} = |\chi| + \O\big(\e^{3r/2}\big)\, , $$ 
whence 
$$ \int_{-\e^{r/2}}^{\frac{\e^2}{\e-y}- 4\pi^2} {\frac{(b\, \sh a - a\,\sin b)(a\,\sh a+b\,\sin b-\ch a + \cos b)}{\sqrt{\chi^2+x^2}}}\, d\chi\, $$ 
{\small 
$$ = \int_{\e^{r/2}}^{4\pi^2} {\frac{\O(\e^{r})\big(\O\big(\e^{2r}/\chi^{2}\big)+\O(\chi)+\O(\e^{3r/2}) \big)}{\chi + \O\big(\e^{3r/2}\big)}}\, d\chi\, = \O(\e^{r})\! \int_{\e^{r/2}}^{4\pi^2} \big[\O(1)+\O(\e^{r}/\chi)\big]\, d\chi \, =\, \O(\e^{r})\, . $$ }
As to the remaining part, using that \  $a\,,b\, \le (\chi^2+x^2)^{1/4} = \O(\e^{r/4})\,$, \  we have\,: 
$$ \int_{0}^{-\e^{r/2}} {\frac{(b\, \sh a - a\,\sin b)(a\,\sh a+b\,\sin b-\ch a + \cos b)}{\sqrt{\chi^2+x^2}}}\, d\chi\, =\, \O(\e^{r})\, . \vspace{-1mm}  $$ 
Therefore we obtain 
$$ \tilde\f\big({\ts\frac{\e^2}{\e-y}} - 4\pi^2,\,{\ts\frac{\e^2}{\e-y}}\,x\big) =\, \frac{5\,\pi}{4} + \O(\e^r)\, , \vspace{-2mm} $$ 
whence \vspace{-2mm} 
$$  \f\big({\ts\frac{\e^2}{\e-y}} - 4\pi^2,\,{\ts\frac{\e^2}{\e-y}}\,x\big) + \5\, \tilde\f\big({\ts\frac{\e^2}{\e-y}} - 4\pi^2,\,{\ts\frac{\e^2}{\e-y}}\,x\big) =\,-\, \frac{3\,\pi}{8} -\,\frac{1}{2}\,\arctg x + \O(\e^{r/3})  \, ,    \hbox{} $$ 
whence the claim (of Lemma \ref{lem.AnalPhiprSing}) follows at once. $\;\diamond$  \parm 

   Lemma \ref{lem.AnalPhiprSing} allows to handle the part of the integral of Proposition \ref{pro.ChangContz=0>yb} which is close to the singularity (saddle-point). Indeed, cutting that integral at $\,\e^r\,$ we get\,: 
$$ \int_{0}^\infty\! \exp\!\left[{ {\frac{\e - y}{\e^2}} }\rt1\! x \right]  \Phi\big(4\eta - 4\pi^2, x\big)\, dx\, = \int_{\e^r}^\infty\! \exp\!\left[{ {\frac{\e - y}{\e^{2}}} }\rt1\! x \right]  \Phi\big(4\eta - 4\pi^2, \,x\big)\, dx\, $$ 
\begin{equation}  \label{f.CutIntegepsr}
+\, {\frac{\e^2}{\e - y}} \int_{0}^{\e^{r-2}(\e-y)}\! \exp\!\left[\rt1\! x \right]\times \Phi\big(4\eta - 4\pi^2, {\ts\frac{\e^2\,x}{\e - y}}\big)\, dx\, .  \vspace{2mm} 
\end{equation} 

   From Lemma \ref{lem.AnalPhiprSing} we deduce the following behavior of the above main contribution to the wanted heat kernel. It will remain to ensure that it is indeed the main one. 
\bpro \label{pro.IntegPartPrinc}  \  For  $\,y\le 0\,$, $\e>0\,$ and $\,r':= \frac{3}{2} -\frac{3}{4}\,1_{\{y=0\}} \le r\le \frac{9}{10}(2-1_{\{y=0\}})$, we have 
$$ \Re\Bigg\{{\frac{\e^2}{\e - y}} \int_{0}^{\e^{r-2}(\e-y)}\! \exp\!\left[\rt1\! x \right]\times   \Phi\!\left({\frac{\e^2}{\e-y}} - 4\pi^2, {\frac{\e^2}{\e - y}}\,x\right)  dx \Bigg\} $$ 
$$ =\, \frac{8\pi^{5/2}\, \e}{\sqrt{\e-y}} \left(1+ \O(\e^{r/6}) + \O\big({\ts\frac{\e^{2-r}}{\e-y}}\big)^{3/2}\right)\! \int_{0}^{\infty}\, \frac{\sin \!\left(\frac{\pi}{8} + x - \frac{1}{2}\,\arctg  x\right)}{(x^2+1)^{1/4}}\, dx\, $$ 
$$ =\, \frac{8\pi^{5/2}\, \e}{\sqrt{\e-y}} \left(1+ \O\big(\e^{\frac{1}{8}(2-1_{\{y=0\}})}\big)\right)\! \int_{0}^{\infty}\, \frac{\sin \!\left(\frac{\pi}{8} + x - \frac{1}{2}\,\arctg  x\right)}{(x^2+1)^{1/4}}\, dx\, $$ 
by eventually taking $\,r= r'$. This holds uniformly with respect to $\,y\le 0$. 
\epro 
\ub{Proof} \quad  Replacing $\,\Phi\big(4\eta - 4\pi^2, {\ts\frac{\e^2}{\e - y}}\,x\big)\,$ according to Lemma \ref{lem.AnalPhiprSing}, for the above value under consideration we obtain\,:
$$ \frac{\e}{\sqrt{\e-y}} \int_{0}^{\e^{r-2}(\e-y)}\, \frac{8\pi^{5/2} + \O(\frac{\e^{2}}{\e-y})} {(x^2+1)^{1/4}}\times \cos\!\left[x -\frac{3\,\pi}{8} - \frac{1}{2}\,\arctg  x + \O(\e^{r/3})\right] dx $$ 
$$ =\, \frac{8\pi^{5/2}\, \e + \O(\frac{\e^{3}}{\e-y}) }{\sqrt{\e-y}} \int_{0}^{\e^{r-2}(\e-y)} \frac{\sin\!\left(\frac{\pi}{8} + x - \frac{1}{2}\,\arctg x + \O(\e^{r/3}) \right)}{(x^2+1)^{1/4}}\, dx  $$ 
$$ =\,  \frac{8\pi^{5/2}\, \e + \O(\frac{\e^{3}}{\e-y}) }{\sqrt{\e-y}} \int_{0}^{\e^{r}/\eta}\, \frac{\sin\!\left(\frac{\pi}{8} + x - \frac{1}{2}\,\arctg x \right) + \O(\e^{r/3})}{(x^2+1)^{1/4}}\, dx\, \quad \hbox{\big(recall $\,{\ts \eta = \frac{\e^2}{4(\e-y)}}\big)$} $$ 
$$ =\,  \frac{8\pi^{5/2}\, \e + \O(\frac{\e^{3}}{\e-y}) }{\sqrt{\e-y}} \left[ \int_{0}^{\e^{r}/\eta}\, \frac{\sin\!\left(\frac{\pi}{8} + x - \frac{1}{2}\,\arctg x \right) }{(x^2+1)^{1/4}}\, dx + \O(\e^{r/6}\,) \right]  \quad  \hbox{(as $\, \eta \ge \e^{4r/3}$)} $$ 
$$ =\,  \frac{8\pi^{5/2}\, \e}{\sqrt{\e-y}} \left(1+ \O(\e^{r/6}) + \O\big(\e^{-r}\eta\big)^{3/2}\right)\! \int_{0}^{\infty}\, \frac{\sin\!\left(\frac{\pi}{8} + x - \frac{1}{2}\,\arctg  x\right)}{(x^2+1)^{1/4}}\, dx\, $$ 
by Proposition \ref {pro.cstte>0} below.  $\;\diamond $

\bpro  \label{pro.cstte>0} \  $(i)$ \  The constant \  ${\ds \sigma := \int_0^\infty\, \frac{\sin\!\left(\frac{\pi}{8} + x - \frac{1}{2}\,\arctg  x\right)}{(x^2+1)^{1/4}}\, dx\,}$ \   is $\,> \frac{1}{10}\,$\raise1.9pt\hbox{.} \parsn
\if{ 
For any $\,n\in\N$, setting $\,\sin^{(-2n)} = (-1)^n\sin\,$ and $\,\sin^{(-2n+1)} = (-1)^n\cos\,$, we have \  
$$ \sigma =\, \frac{(2n)!}{4^n\,n!} \int_0^{\pi/2} \sin^{(-n)}\!\left[(n+\5)\,\theta -{\ts\frac{\pi}{8}} + \cotg\theta\right] (\sin\theta)^{n-\frac{3}{2}}\, d\theta\, - \cos\!\left[{\ts\frac{\pi}{8}}\right]\times \sum_{k=0}^{n-1} \frac{(2k)!}{4^{k}\,k!}\, \raise2pt\hbox{.} \vspace{-1mm} $$ 
}\fi 
$(ii)$ \  For large $R$ we have \quad 
${\ds \int_{R}^{\infty}\, \frac{\sin\!\left(\frac{\pi}{8} + x - \frac{1}{2}\,\arctg x\right)}{(x^2+1)^{1/4}}\, dx\, =\, \O\big(R^{-3/2}\big) . \vspace{-2mm}  }$ 
\epro 
\ub{Proof} \quad $(i)$ \  Let \  $\theta(x) := \frac{\pi}{8} + x - \frac{1}{2}\,\arctg x = x - \frac{\pi}{8}  + \frac{1}{2}\,\arctg \frac{1}{x}\,$, which increases (from $\frac{\pi}{8}$ to infinity), as well as its derivative $\theta'$ and then $\,\theta'\circ\theta\1$.  
\if{ 
Let $\,x_k := \theta\1(k\pi)$, for $\,k\in\N^*$. \quad 
As $\,\theta(x) > x- \frac{\pi}{8}\,$ we have \  $x_k< \frac{\pi}{8} + k\pi\,$. \parn   
As $\,\theta(x) < x - \frac{\pi}{8} + \frac{1}{2x}\,$ \  we have \  $\,x_k + \frac{1}{2x_k} > \frac{\pi}{8} + k\pi\,$ \  and then, since  $\,\frac{2}{(\frac{\pi}{8} + k\pi)^2} \le \frac{2}{(\frac{\pi}{8} + \pi)^2} <2(\sqrt{2} -1)$, we have  \vspace{-2mm}  
$$  k\pi + \frac{\pi}{8} >\, x_k\, > \left(\frac{\pi}{8} + k\pi\right) \frac{1+\sqrt{1-\frac{2}{(\frac{\pi}{8} + k\pi)^2}}}{2}\, > \left(\frac{\pi}{8} + k\pi\right)\!\left(1- \frac{1/2}{(\frac{\pi}{8} + k\pi)^2} - \frac{1/2}{(\frac{\pi}{8} + k\pi)^4} \right) $$
$$ = \,\frac{\pi}{8} + k\pi - \frac{1/2}{\frac{\pi}{8} + k\pi} - \frac{1/2}{(\frac{\pi}{8} + k\pi)^3} >  \frac{\pi}{8} + k\pi - \frac{13/12}{\frac{\pi}{8} + k\pi} >  \frac{\pi}{8} + k\pi - \frac{\pi}{10} =  k\pi + \frac{\pi}{40}\, \raise2pt\hbox{.} $$
}\fi 
\parn
Then changing the variable we have\,: \quad 
${\ds  \sigma\, = \int_{\pi/8}^\infty \frac{\sin t\; dt}{\left(1+\theta\1(t)^2\right)^{1/4}\times \theta'\circ\theta\1(t)} \, \raise1.5pt\hbox{} }$ 
$$ > \int_{\pi/8}^{\pi/2} \frac{\sin t\; dt}{\left(1+\theta\1(\frac{\pi}{2})^2\right)^{1/4}\times \theta'\circ\theta\1(\frac{\pi}{2})} + \int_{\pi/2}^{\pi} \frac{\sin t\; dt}{\left(1+\theta\1(\pi)^2\right)^{1/4}\times \theta'\circ\theta\1(\pi)} $$
$$ - \int_{\pi}^{3\pi/2} \frac{|\sin t|\; dt}{\left(1+\theta\1(\pi)^2\right)^{1/4}\times \theta'\circ\theta\1(\pi)} - \int_{3\pi/2}^{2\pi} \frac{|\sin t|\; dt}{\left(1+\theta\1(\frac{3\pi}{2})^2\right)^{1/4}\times \theta'\circ\theta\1(\frac{3\pi}{2})} $$ 
{\small 
$$ +\sum_{k\ge 1} \left[ \int_{2k\pi}^{(2k+1)\pi}\! {\ts\frac{\sin t\; dt}{\left(1+\theta\1((2k+1)\pi)^2\right)^{1/4}\, \theta'\circ\theta\1((2k+1)\pi)}} - \int_{(2k+1)\pi}^{(2k+2)\pi}\! {\ts\frac{|\sin t|\; dt}{\left(1+\theta\1((2k+1)\pi)^2\right)^{1/4}\, \theta'\circ\theta\1((2k+1)\pi)}} \right] $$}
$$ = \int_{\pi/8}^{\pi/2} \frac{\sin t\; dt}{\left(1+\theta\1(\frac{\pi}{2})^2\right)^{1/4}\times \theta'\circ\theta\1(\frac{\pi}{2})} - \int_{3\pi/2}^{2\pi} \frac{|\sin t|\; dt}{\left(1+\theta\1(\frac{3\pi}{2})^2\right)^{1/4}\times \theta'\circ\theta\1(\frac{3\pi}{2})} $$ 
$$ =\, \frac{\cos(\pi/8)}{\left(1+\theta\1(\frac{\pi}{2})^2\right)^{1/4}\times \theta'\circ\theta\1(\frac{\pi}{2})} - \frac{1}{\left(1+\theta\1(\frac{3\pi}{2})^2\right)^{1/4}\times \theta'\circ\theta\1(\frac{3\pi}{2})} \, \raise1.5pt\hbox{.} $$ 
Now observing the following\,: \quad $*$ \  $ \cos(\pi/8) = \sqrt{\frac{2+\sqrt{2}}{2}}\, > \frac{92}{100}\,$ \raise1.5pt\hbox{;} \parn 
$*$ \  as $\,\theta(x) > x- \frac{\pi}{8}\,$\raise1.1pt\hbox{,}  we have \  $\theta\1(\frac{\pi}{2}) < \frac{5\pi}{8} \,$; \quad 
$*$ \  as $\,\theta(x) < x- \frac{\pi}{8} + \frac{1}{2x}\,$\raise1.1pt\hbox{,}  we have \  $\theta\1(x+\frac{1}{2x+\frac{\pi}{4}}) > x+ \frac{\pi}{8} \,$ \  and then \quad 
$ \frac{3\pi}{2} > \frac{471}{100} < \frac{9}{2}+\frac{1}{9+\frac{\pi}{4}}\, \Rightarrow\, \theta\1(\frac{3\pi}{2}) > \frac{9}{2}+\frac{\pi}{8} > \frac{3\pi}{2} \,$; \parn
$*$ \  $\theta'(x) =  1- \frac{1}{2(1+x^2)} \,$\raise1.5pt\hbox{,}  \   and then \  $\left(1+\theta\1(t)^2\right)^{1/4}\times \theta'\circ\theta\1(t) = \frac{1+2\,\theta\1(t)^2}{2\,\left(1+\theta\1(t)^2\right)^{3/4}} \,$\raise1.7pt\hbox{,} \parn 
this yields\,: 
$$ \sigma\, >\, \frac{92\times 2\times (1+\frac{25\pi^2}{64})^{3/4} }{100 \times (1+\frac{25\pi^2}{32}) } - \frac{2\times (1+\frac{9\pi^2}{4})^{3/4} }{1+\frac{9\pi^2}{2}}  \, \raise1.9pt\hbox{.}  $$ 
Using \  $\frac{25\pi^2}{64} < \frac{386}{100}\,$ and $\,\frac{9\pi^2}{4} > \frac{222}{10}\,$\raise1.2pt\hbox{,} \  this gives 
$$ \sigma\, >\, \frac{184\times (4.86)^{3/4} }{872} - \frac{2\times (23.2)^{3/4} }{45.4}\,>\, \frac{1.8\times 27}{88} - \frac{22}{45.4}\,>\,  \frac{1}{10} \, \raise1.9pt\hbox{.}  $$ 
$(ii)$ \   Performing the change \  $x=\cotg\theta=\tg(\frac{\pi}{2}-\theta)\,$ we get\,: \vspace{-2mm} 
$$ \int_{R}^{\infty}\, \frac{\sin\!\left(\frac{\pi}{8} + x - \frac{1}{2}\,\arctg  x\right)}{(x^2+1)^{1/4}}\, dx\, =
\int_0^{\arctg(1/R)} \sin\!\left[\,\5\,\theta -{\ts\frac{\pi}{8}} + \cotg\theta\right] (\sin\theta)^{-\frac{3}{2}}\, d\theta \, \vspace{-2mm} $$ 
$$ =\, \5 \int_0^{\arctg\frac{1}{R}}\! \sin\!\left[{\ts\frac{\theta}{2}} - {\ts\frac{\pi}{8}} + \cotg\theta\right] \sqrt{\sin\theta}\, d\theta\, + \int_0^{\arctg\frac{1}{R}}\! {\ts\frac{d}{d\theta}}\! \left\{\cos\!\left[{\ts\frac{\theta}{2}} - {\ts\frac{\pi}{8}} + \cotg\theta\right] \right\}  \sqrt{\sin\theta}\, d\theta $$
$$ =\, \O\big(R^{-3/2}\big) + \5 \int_0^{\arctg\frac{1}{R}}\! \Big( \sin\!\left[{\ts\frac{\theta}{2}} - {\ts\frac{\pi}{8}} + \cotg\theta\right]\!\sin\theta - \cos\!\left[{\ts\frac{\theta}{2}} - {\ts\frac{\pi}{8}} + \cotg\theta\right]\! \cos\theta\Big) \frac{d\theta}{\sqrt{\sin\theta}} \vspace{-1mm}  $$ 
$$ =\, \O\big(R^{-3/2}\big) - \5\int_0^{\arctg\frac{1}{R}} \cos\!\left[{\ts\frac{3\theta}{2}} - {\ts\frac{\pi}{8}} + \cotg\theta\right] \frac{d\theta}{\sqrt{\sin\theta}}   $$ 
$$ =\, \O\big(R^{-3/2}\big) - \5 \int_0^{\arctg\frac{1}{R}}\! \Big( {\ts\frac{3}{2}} \cos\!\left[{\ts\frac{3\theta}{2}} - {\ts\frac{\pi}{8}} + \cotg\theta\right] - {\ts\frac{d}{d\theta}}\! \left\{ \sin\!\left[{\ts\frac{3\theta}{2}} - {\ts\frac{\pi}{8}} + \cotg\theta\right]\right\}\!\! \Big) (\sin\theta)^{3/2}\,d\theta  $$ 
$$ =\, \O\big(R^{-3/2}\big) -{\ts\frac{3}{4}} \int_0^{\arctg\frac{1}{R}} \sin\!\left[{\ts\frac{3\theta}{2}} - {\ts\frac{\pi}{8}} + \cotg\theta\right] \cos\theta\, \sqrt{\sin\theta}\,d\theta\, =\,\O\big(R^{-3/2}\big)  .  \;\;\diamond $$  
\if{  
Then for any $\,n\in\N$ we similarly have\,: \vspace{-2mm} 
$$ \int_0^{\pi/2} \sin^{(-n)}\!\left[(n+\5)\,\theta -{\ts\frac{\pi}{8}} + \cotg\theta\right] (\sin\theta)^{n-\frac{3}{2}}\, d\theta\, \vspace{-2mm}  $$
$$ = \int_0^{\pi/2} \sin^{(-n)}\!\left[(n+\5)\,\theta -{\ts\frac{\pi}{8}} + \cotg\theta\right] \left((n+\5) - \big[(n+\5) - \sin\2\!\theta\big] \right) (\sin\theta)^{n+\frac{1}{2}}\, d\theta\, $$
$$ =\, (n+\5) \int_0^{\pi/2} \sin^{(-n)}\!\left[(n+\5)\,\theta -{\ts\frac{\pi}{8}} + \cotg\theta\right] (\sin\theta)^{n+\frac{1}{2}}\, d\theta\, $$
$$ - \int_0^{\pi/2} {\frac{d}{d\theta}}\! \left\{\sin^{(-n-1)}\!\left[(n+\5)\,\theta -{\ts\frac{\pi}{8}} + \cotg\theta\right] \right\}  (\sin\theta)^{n+\frac{1}{2}}\, d\theta $$
$$ =\, (n+\5)\! \int_0^{\pi/2} \sin^{(-n)}\!\left[(n+\5)\,\theta -{\ts\frac{\pi}{8}} + \cotg\theta\right] (\sin\theta)^{n+\frac{1}{2}}\, d\theta\, - \sin^{(-n-1)}\!\left[(n+\5){\ts\frac{\pi}{2}} -{\ts\frac{\pi}{8}}\right] $$
$$ +\, (n+\5)\int_0^{\pi/2} \sin^{(-n-1)}\!\left[(n+\5)\,\theta -{\ts\frac{\pi}{8}} + \cotg\theta\right] \cos\theta\,  (\sin\theta)^{n-\frac{1}{2}}\, d\theta $$
$$ =\, (n+\5)\! \int_0^{\pi/2} \sin^{(-n-1)}\!\left[(n+{\ts\frac{3}{2}})\,\theta -{\ts\frac{\pi}{8}} + \cotg\theta\right] (\sin\theta)^{n-\frac{1}{2}}\, d\theta \, - \sin^{(-n-1)}\!\left[{\ts\frac{(4n+1)\pi}{8}}\right] $$ 
$$ =\, (n+\5)\! \int_0^{\pi/2} \cos^{(-n-1)}\!\left[(n+{\ts\frac{3}{2}})\,\theta -{\ts\frac{\pi}{8}} + \cotg\theta\right] (\sin\theta)^{n-\frac{1}{2}}\, d\theta\, - \cos\!\left[{\ts\frac{\pi}{8}}\right] . $$ 
This shows the claim by induction. $\;\diamond$  \parmn \noindent
\ub{Rque} \quad  Setting \  $p_k:= \frac{(2k)!}{4^{k}\,k!}\,$\raise1pt\hbox{,} \  we have \  ${\sigma/p_n = I_n - \cos\!\left[{\ts\frac{\pi}{8}}\right] \sum_{k=0}^{n-1}\limits p_k/p_n}\,$, \  whence \parn  
$I_n = v_n  \cos\!\left[{\ts\frac{\pi}{8}}\right]  +\O(p_n\1)$, \  with \  $(n-\5)\, v_n = 1+v_{n-1}\, $ and \  $v_n < \frac{1}{n-\frac{1}{2}} + \frac{1}{(n-\frac{1}{2})(n-\frac{3}{2})} \left(1+ \frac{n-2}{n-\frac{5}{2}}\right)$ 
$< \frac{1}{n-\frac{1}{2}} \times \frac{2n+1}{2n-3} \, $\raise1pt\hbox{,} and then \  $ v_n= \frac{1}{n} + \O(n\2)$, \  
$I_n = \frac{\cos({\pi}/{8})}{n} +\O(n\2)\,$... ?? \parb 
}\fi 
\parb 
   To conclude the first sub-case $\big[y\le0,\, z=0\big]$ we are dealing with, we have now to handle the remaining integral in (\ref{f.CutIntegepsr}), namely \  
\if{ 
\  ${ \ds {\frac{\e^2}{\e - y}} \int_{\e^{r-2}(\e-y)}^\infty  e^{\rt1\! x}\, \Phi\big(4\eta - 4\pi^2, {\ts\frac{\e^2}{\e - y}}\,x\big)\, dx\,}$.  \parn
Let us begin with the part 
${ \ds {\frac{\e^2}{\e - y}} \int_{\e^{r-2}(\e-y)}^{\e^{r'-2}(\e-y)} e^{\rt1\! x}\, \Phi\big(4\eta - 4\pi^2, {\ts\frac{\e^2}{\e - y}}\,x\big)\, dx\,}$, \  recalling the notation $\,r':=\frac{3}{2} -\frac{3}{4}\,1_{\{y=0\}}\,$ in Lemma \ref{lem.AnalPhiprSing} and Proposition \ref{pro.IntegPartPrinc}, and that in Proposition \ref{pro.IntegPartPrinc} we eventually took the optimal $\,r = \frac{6r'}{5}\,$\raise1pt\hbox{.} \  Using Lemma \ref{lem.AnalPhiprSing} as in the proof of Proposition \ref{pro.IntegPartPrinc}, and then noticing that $\,\eta \asymp \e^{4r'/3}$\raise1pt\hbox{,} \  we obtain 
$$ \Re\Bigg\{{\frac{\e^2}{\e - y}} \int_{\e^{r-2}(\e-y)}^{\e^{r'-2}(\e-y)} e^{\rt1\! x}\, \Phi\big(4\eta - 4\pi^2, {\ts\frac{\e^2}{\e - y}}\,x\big)\, dx\Bigg\} $$
$$ =\, \frac{\O( \e) }{\sqrt{\e-y}}  \int_{\e^{r}/\eta}^{\e^{r'}/\eta}\, \frac{\sin\!\left(\frac{\pi}{8} + x - \frac{1}{2}\,\arctg x \right) + \O(\e^{r'/3})}{(x^2+1)^{1/4}}\, dx\,  $$ 
$$ =\, \frac{\O( \e) }{\sqrt{\e-y}}  \int_{\e^{r}/\eta}^{\infty}\, \frac{\sin\!\left(\frac{\pi}{8} + x - \frac{1}{2}\,\arctg x \right)}{(x^2+1)^{1/4}}\, dx\, + \frac{\O( \e^{1+r'/3}) }{\sqrt{\e-y}}  \times \frac{\e^{r'}/\eta}{(\e^{r}/\eta)^{1/2}}\,   $$ 
$$ =\, \frac{\O( \e) }{\sqrt{\e-y}}\times \O\big(\e^{-r}\eta\big)^{3/2} + \O( \e^{4r'/3-r/2}) \, =\,\frac{\O( \e) }{\sqrt{\e-y}}\times \left[\O\big(\e^{2r' -9r'/5}\big) + \O\big( \e^{11r'/15- 2r'/3}\big) \right]   $$ 
$$ =\, \frac{\O( \e) }{\sqrt{\e-y}}\times \O\big( \e^{r'/15}\big)\,  =\, \frac{o( \e) }{\sqrt{\e-y}}\, \raise1.9pt\hbox{.}   $$ 
{\bf NB} \  Ce morceau peut \^etre int\'egr\'e directement dans la Proposition \ref{pro.IntegPartPrinc} en prenant d'embl\'ee $\,r=r'$... \parsn 
The last contribution in (\ref{f.CutIntegepsr}) is 
}\fi 
${\ds \int_{\e^{r'}}^{\infty} e^{\rt1\! \frac{\e - y}{\e^2}\,x}\, \Phi\big(4\eta - 4\pi^2, x\big)\, dx\, }$ (recall that $\,r'$ was defined in Proposition \ref{pro.IntegPartPrinc}).
\blem \label{lem.ControlIntPhi} \  We uniformly have \quad ${\ds \int_{\e^{r'}}^{\infty} \exp\!\left[{ {\frac{\e - y}{\e^2}} }\rt1\! x \right]  \Phi\big(4\eta - 4\pi^2, x\big)\, dx\, }$ \vspace{-2mm} 
$$ =\, \frac{\e}{\sqrt{\e-y}} \times o\big(\e^{\frac{1}{20}}\big) + {\frac{\e^2\rt1\!}{\e - y}} \int_{\e^{r'}}^{9 \log^2\!\e}\! \exp\!\left[{ {\frac{\e - y}{\e^2}} }\rt1\! x \right] \frac{d \Phi}{dx}\big({\ts\frac{\e^2}{\e-y}} - 4\pi^2, x\big)\, dx \, . \vspace{-3mm}  $$ 
\elem
\ub{Proof} \quad  Lemma \ref{lem.LiftSqRtxPsi} provides the uniform estimate\,: $\,{\ds \Phi\big({\ts\frac{\e^2}{\e-y}} - 4\pi^2, x\big) = {\ts\frac{\ds \O(x)}{{\sqrt{\ds\ch\!\sqrt{x/2}}}}}\,}$ for large $\,x\,$, so that we have \vspace{-2mm} 
$$ \int_{R}^{\infty} e^{\rt1\! \frac{\e - y}{\e^2}\,x}\, \Phi\big(4\eta - 4\pi^2, x\big)\, dx\, =  \O\!\left[\int_{\sqrt{R/2}}^{\infty} \frac{x^3}{\sqrt{\ch x}}\, dx\right]  =  \O\!\left[\int_{\sqrt{R/8}}^{\infty} x^3 e^{-x}\, dx\right] $$
$$  =\,   \O\!\left[  R^{3/2} e^{-\sqrt{R/8}}\,\right] , \vspace{-2mm} $$ 
and then \vspace{-1mm}
$$ \int_{9 \log^2\!\e}^{\infty} e^{\rt1\! \frac{\e - y}{\e^2}\,x}\, \Phi\big(4\eta - 4\pi^2, x\big)\, dx\, =\,\O\!\left[  \e^{3/2\sqrt{2}}\,\log^3({\ts\frac{1}{\e}})\,\right] = o\big(\e^{1+\frac{1}{20}}\big) = \frac{\e}{\sqrt{\e-y}} \times o\big(\e^{\frac{1}{20}}\big) . $$ 
It remains to control \  ${\ds \int_{\e^{r'}}^{9 \log^2\!\e} e^{\rt1\! \frac{\e - y}{\e^2}\,x}\, \Phi\big(4\eta - 4\pi^2, x\big)\, dx\, . }$\quad  Integrating by parts we have\,: 
$$ \int_{\e^{r'}}^{9 \log^2\!\e}\! \exp\!\left[{ {\frac{\e - y}{\e^2}} }\rt1\! x \right]  \Phi\big({\ts\frac{\e^2}{\e-y}} - 4\pi^2, x\big)\, dx \, =\, {\frac{\e^2\rt1\!}{\e - y}}  \exp\!\left[{ {\frac{\e - y}{\e^{2-r'}}} }\rt1\!\right]  \Phi\big({\ts\frac{\e^2}{\e-y}} - 4\pi^2, \e^{r'}\big)\, $$
$$ - \, {\frac{\e^2\rt1\!}{\e - y}}  \exp\!\left[{ {\frac{\e - y}{\e^{2}}} }\rt1 9 \log^2\!\e\right]  \Phi\big({\ts\frac{\e^2}{\e-y}} - 4\pi^2, 9 \log^2\!\e\big) $$ 
\begin{equation*}  
+\,  {\frac{\e^2\rt1\!}{\e - y}} \int_{\e^{r'}}^{9 \log^2\!\e}\! \exp\!\left[{ {\frac{\e - y}{\e^2}} }\rt1\! x \right]  \frac{d \Phi}{dx}\big({\ts\frac{\e^2}{\e-y}} - 4\pi^2, x\big)\, dx \, . 
\end{equation*}
Now, according to Lemma \ref{lem.AnalPhiprSing} we uniformly have \  
${\ds \Phi\big({\ts\frac{\e^2}{\e-y}} - 4\pi^2, \,\e^{r'}\big) = \O(\e^{-r'/2}) }\,$,  
so that (noticing that $\,\eta \asymp \e^{4r'/3}$) the first term of the right hand side is $\,{\ds \frac{\e\times \O\big(\e^{r'/6}\big)}{\sqrt{\e-y}}}\,$\raise1pt\hbox{.} \  Moreover by Lemma \ref{lem.LiftSqRtxPsi} we have \  ${\ds \Phi\big({\ts\frac{\e^2}{\e-y}} - 4\pi^2, 9 \log^2\!\e\big) = \O\!\left[ \e^{3/\sqrt{2}}\,\log^2({\ts\frac{1}{\e}})\right] = o(\e^2)}\,$, so that the second term of the right hand side is $\,{\ds \frac{\e\times o\big(\e^{2}\big)}{\sqrt{\e-y}}}\,$\raise1.5pt\hbox{.} \  Finally we have $\,\O\big(\e^{r'/6}\big) = o\big(\e^{\frac{1}{20}}\big) \,$. $\;\diamond$ 
\parsn
Then Lemma \ref{lem.LiftSqRtxPsi} allows to estimate $\,{\ds\frac{d \Phi}{dx}\big({\ts\frac{\e^2}{\e-y}} - 4\pi^2, x\big)}$. The following will be sufficient. 
\blem \label{lem.ControlPhi'} \  For any $\,y\le 0\,$ and $\,x\ge \e^{r'}$, uniformly as $\,\e \sea 0$  \  we have 
$$ \frac{d \Phi}{dx}\big({\ts\frac{\e^2}{\e-y}} - 4\pi^2, x\big) =\, 1_{\{\e^{r'}\le x< \e^{2r'/3}\}} \frac{\O( \e^{4r'/3})}{x^{5/2}} + 1_{\{\e^{2r'/3}\le x< 1\}} \frac{\O(1)}{\sqrt{x}} + 1_{\{x\ge1\}} \O\!\left( x\, e^{-\sqrt{x/8}}\,\right) . \vspace{-2mm} $$ 
\elem
\ub{Proof} \quad As already used above, Lemma \ref{lem.AnalPhiprSing} ensures that \  
${\ds \Phi\big({\ts\frac{\e^2}{\e-y}} - 4\pi^2, \,x\big) = \O(x^{-1/2}) }\,$ for small $\,x\,$ and Lemma \ref{lem.LiftSqRtxPsi} ensures that \   ${\ds \Phi\big({\ts\frac{\e^2}{\e-y}} - 4\pi^2, x\big) = \O\!\left[ x\, e^{-\sqrt{x/8}}\,\right]}$ for large $\,x\,$, uniformly with respect to $\,y\le 0$.  \  
Then using Lemma \ref{lem.LiftSqRtxPsi} again and (\ref{f.racinexkhi}) we obtain 
$$ \frac{d \log\Phi}{dx}(\chi, x)\, =\, \rt1 \frac{d(\f+\frac{1}{2}\tilde\f)}{dx}(\chi, x) + \frac{5\,x}{4\sqrt{\chi^2+x^2}} + \frac{a\sin b- b\, \sh a}{8 (a^2+b^2)(\ch a-\cos b)} $$
$$ +\,  \frac{a\sin b - b\, \sh a}{8 \big((a^2+b^2)(\ch a+\cos b) -4(a\, \ch a +b\sin b)+4 (\ch a-\cos b)\big)}\, $$ 
$$ =\, \rt1 \frac{d(\f+\frac{1}{2}\tilde\f)}{dx}(\chi, x)  +\O(1)\,. $$ 
Then by (\ref{f.racinexkhi}) we have \quad  
${\ds a\, \frac{\partial a}{\partial x} = b\, \frac{\partial b}{\partial x} = \frac{x}{4\sqrt{\chi^2+x^2}} =\O(1) , \, \frac{\partial a}{\partial x} = \frac{b/2}{a^2+b^2}\,\raise1pt\hbox{,}\, \frac{\partial b}{\partial x} = \frac{a/2}{a^2+b^2}\,\raise1pt\hbox{,}}$ \parn
and by Lemma \ref{lem.ProlAnalEspCompl} we have \vspace{-1mm} 
$$ \frac{d \f}{dx}(\chi, x) =\,\5\, \frac{d }{dx}\big[\arctg(b/a) -\f_a(b)\big] = \frac{\chi}{4\,(\chi^2+x^2)} - \frac{a\,\sh\!(2a)- b\, \sin(2b)}{4(a^2+b^2) (\ch\!(2a)-\cos (2b))} $$
which is bounded for $\,x\ge 1\,$; while for small $\,x$, proceeding as in the proof of Lemma \ref{lem.AnalPhiprSing} we successively have \  $\ch\!(2a) - \cos(2b) \sim \frac{x^2}{8\pi^2}\,$\raise1pt\hbox{,} $\,a\,\sh\!(2a)- b\, \sin(2b) \sim 2\eta+ \frac{3x^2}{32\pi^2}\,$\raise1pt\hbox{,} \  whence 
$$ \frac{d \f}{dx}(\chi, x) = \O(1 + \eta\,x\2) =\, 1_{\{\e^{r'}\le x< \e^{2r'/3}\}}\O( \e^{4r'/3}/x^2) + 1_{\{x\ge \e^{2r'/3}\}} \O(1)\, . $$
Then by Lemma \ref{lem.LiftSquaRoot} and (\ref{df.fnf}) we have 
$$ \frac{d \tilde\f}{dx}(\chi, x) =\, 
\frac{d }{dx} \arg \!\left[\frac{\chi + \rt1\! x}{f(\chi,x)}\right] = \frac{d }{dx}\, \arctg\frac{x}{\chi} - \frac{d }{dx} \arg \!\left[1- {\ts\frac{a\, \sh a + b\, \sin b + \rt1 (a\, \sin b - b\,\sh a)}{(a^2+b^2)\, (\ch a + \cos b)\big/ 2}}\right] $$ 
$$ =\, \frac{\chi}{\chi^2+x^2} -  \frac{d }{dx}\, \arctg \!\left[ { \frac{2\, (a\, \sin b - b\,\sh a)}{2(a\, \sh a + b\, \sin b) - (a^2+b^2)\, (\ch a + \cos b)}}\right] $$ 
$$ = \frac{\chi}{\chi^2+x^2} - {\ts \frac{\big(b\,\partial_a (a\, \sin b - b\,\sh a)+ a\,\partial_b (a\, \sin b - b\,\sh a)\big)\left[2(a\, \sh a + b\, \sin b) - (a^2+b^2)\, (\ch a + \cos b)\right] }{ (a^2+b^2) \left(\big[2 (a\, \sin b - b\,\sh a)\big]^2+\big[2(a\, \sh a + b\, \sin b) - (a^2+b^2)\, (\ch a + \cos b)\big]^2\right)}}  $$ 
$$ + {\ts \frac{\big(b\,\partial_a[2(a\, \sh a + b\, \sin b) - (a^2+b^2)\, (\ch a + \cos b)]+a\,\partial_b[2(a\, \sh a + b\, \sin b) - (a^2+b^2)\, (\ch a + \cos b)]\big) (a\, \sin b - b\,\sh a) }{ (a^2+b^2) \left(\big[2 (a\, \sin b - b\,\sh a)\big]^2+\big[2(a\, \sh a + b\, \sin b) - (a^2+b^2)\, (\ch a + \cos b)\big]^2\right)}}  $$ 
{\small 
$$ = \frac{\chi}{\chi^2+x^2} + { \frac{\O\big((a^2+b^2)^2(\ch a + \cos b)^2\big) }{ (a^2+b^2)^2 (\ch a + \cos b)\big[(a^2+b^2)\, (\ch a + \cos b)-4 (a\, \sh a + b\, \sin b) + 4\, (\ch a - \cos b)\big]}}  $$ }
which is bounded. \  The claim follows. $\;\diamond$ \parm 

   The following actually states that the main contribution to the wanted heat kernel was indeed given by the part considered in Proposition \ref{pro.IntegPartPrinc} ; in other words, that the dominant saddle point of the oscillatory integral $\,\II_\e(y,0)$ was indeed located near $\, 4\rt1\!\! \left(\pi^2-\frac{\e^2}{4(\e-y)}\right)$. \vspace{-1mm} 
\blem \label{lem.ControlIntPhif} \  We uniformly have \   ${\ds \int_{\e^{r'}}^{\infty} e^{\rt1\! \frac{\e - y}{\e^2}\,x}\, \Phi\big(4\eta - 4\pi^2, x\big)\, dx\, =\, \frac{\e}{\sqrt{\e-y}} \times o\big(\e^{\frac{1}{20}}\big)\, . \vspace{-0mm}  }$ 
\elem
\ub{Proof} \quad By Lemmas \ref{lem.ControlIntPhi} and \ref{lem.ControlPhi'} (recalling $\,\eta \asymp \e^{4r'/3}$) we have\,:  
$$ \int_{\e^{r'}}^{\infty} e^{\rt1\! \frac{\e - y}{\e^2}\,x}\, \Phi\big(4\eta - 4\pi^2, x\big)\, dx\, -\, \frac{\e}{\sqrt{\e-y}} \times o\big(\e^{\frac{1}{20}}\big) $$
$$ =\, {\frac{\e\,\O(\e^{2r'/3})}{\sqrt{\e - y}}} \left[\int_{\e^{r'}}^{\e^{2r'}}\! \frac{\O( \e^{4r'/3})}{x^{5/2}}\, dx + \int_{\e^{2r'}}^{1}\! \frac{\O(1)}{\sqrt{x}}\, dx + \int_{1}^{9 \log^2\!\e}\! \O\!\left( x\, e^{-\sqrt{x/8}}\,\right) dx \right]  $$ 
$$ =\, {\frac{\e\,\O(\e^{2r'/3})}{\sqrt{\e - y}}} \left[ \e^{4r'/3-3r'/2} + 1 \right]  =\, {\frac{\e}{\sqrt{\e - y}}}\times \,\O(\e^{r'/2}) = \, \frac{\e}{\sqrt{\e-y}} \times o\big(\e^{\frac{1}{20}}\big)\, . \;\; \diamond $$ \parm 

   We can now conclude this section, by the following exact equivalent of the oscillatory integral $\,\II_\e(y,0)$ arising in the case  $\,y\le 0=z\,$. 
\bpro  \label{pro.ControlIntPhi} \    For $\,y\le 0\,$,  uniformly as $\,\e \sea 0$  \   we have 
$$ \Re\Bigg\{\!\int_{0}^\infty\! \exp\!\left[{ {\frac{\e - y}{\e^2}} }\rt1\! x \right] \times   \Phi\!\left({\frac{\e^2}{\e-y}} - 4\pi^2 , \,x\right)  dx \Bigg\} =\, \frac{8\pi^{5/2}\,\sigma\, \e}{\sqrt{\e-y}} \left(1+o\big(\e^{\frac{1}{20}}\big)\right) \vspace{-0mm}  $$ 
(recall $\,\sigma\,$ is that of Proposition \ref{pro.cstte>0}) \   and 
$$ \II_\e(y,0) = \, \exp\!\left[4\pi^2\,\frac{y-\e}{\e^2}\right]\times\frac{\e}{\sqrt{\e-y}} \times 16\,\pi^{5/2}\,e\, \sigma  \times\left(1+o\big(\e^{\frac{1}{20}}\big)\right) . \vspace{-2mm} $$ 
\epro 
\ub{Proof} \quad  The first claim follows directly from (\ref{f.CutIntegepsr}) with $\,\eta = \frac{\e^2}{4(\e-y)}\,$ and $r=r'= \frac{3}{2} -\frac{3}{4}\,1_{\{y=0\}}$, \parn Lemma \ref{lem.ControlIntPhif} and Proposition \ref{pro.IntegPartPrinc}. With Proposition \ref{pro.ChangContz=0>yb}, it entails  the second one. $\;\diamond$ \parm 

   Propositions \ref{pro.equivIntPe0} and \ref{pro.ControlIntPhi} together give the following wanted small time equivalent, which is the content of Theorem \ref{th.mainres}$(ii)$. 
\bcor \label{cor.w=0=z>=y} \  For $\,y\le 0\,$, uniformly as $\,\e\sea 0\,$ we have 
$$ p_\e\big(0\, ; (0,y, 0)\big) \sim\, \exp\!\left[4\pi^2\,\frac{y-\e}{\e^2}\right]\times\frac{2\sqrt{2}\,\,e}{\e^{3}\,\sqrt{\e-y}}\times \int_0^\infty\, \frac{\sin\!\left(\frac{\pi}{8} + x - \frac{1}{2}\,\arctg  x\right)}{(x^2+1)^{1/4}}\, dx\, . \vspace{-3mm} $$ 
\ecor    
   
\subsection{Second sub-case\,: $\, y> 0\,,\, z=0$} \label{sec.znot=0<y} \indf 
       Having dealled in Section \ref{sec.z=0>y1} with the sub-case  $\,y\le 0\,$, we deal here with the sub-case $\,y> 0\,$, in order (recall Proposition \ref{pro.exprAvChgCont}) to handel
$$ \II_\e(y,0) =\, 2\; \Re\Bigg\{\! \int_{0}^\infty\! \exp\!\left[{ {\frac{\e - y}{\e^2}} }\rt1\! x \right] \Phi(x)\, dx \Bigg\} . $$   

   The dominant saddle point of the above oscillatory integral, not easy to compute, happens to be now close to $\, -\rt1\pi^2$ \Big(very roughly\,: because $\pi^2 $ is the first positive zero of {\small $\big(\ch\!\sqrt{-2\rt1\! x}-\cos\!\sqrt{-2\rt1\! x}\,\big)$}\Big), as will be confirmed later, by Proposition \ref{pro.ControlIntPhi>}. Since it is located under the real axis, the situation will be partly simpler as in Section \ref{sec.z=0>y1} (Corollary \ref{cor.ValFourIntomt} lets guess that the case of a positive $\chi$ should be a priori easier). We however follow the same route, using part of the work already made in this former case. \par  
   To use that $\,-\rt1\pi^2\,$ is indeed close to a saddle point, we need to expand $\Phi$ partially. \parn 
   Here is the analogue of Lemma \ref{lem.AnalPhiprSing}. 
\blem \label{lem.AnalPhiprSing>} \  For small $\,h\,$, \   we have 
\quad $ \Phi\big(\pi^2 , h\big) \equiv \,\Phi\big(-\rt1\!\pi^2 + h\big) = \, \vspace{-2mm} $ 
$$  \frac{2^{3/4}\,\pi^{11/4}\, e^{5\rt1\!\pi/16}}{ \sqrt{\pi -2\,\th\!{\frac{\pi}{2}}}\, ( \sh\pi)^{1/4}}\times h^{-1/4} \times  \exp\!\left[ \rt1\! \!\left[{\ts \frac{5/4 - 2\pi\,\coth\pi - \pi\,\th\pi}{8\,\pi^2} + \frac{3\pi -6\,\th\!{\frac{\pi}{2}} -\pi\,\th\!^2{\frac{\pi}{2}}}{4\pi^2 \left(\pi -2\,\th\!{\frac{\pi}{2}}\right)} }\right]\! h + \O(h^2)\right] .  $$  
\elem 
\ub{Proof} \quad  To begin, let us directly use the expression for $\,\Phi\,$ given in Proposition \ref{pro.exprAvChgCont}, applied first with $\,\sqrt{\frac{x}{2}} = (1-\rt1\!)\frac{\pi}{2} + u\,$, \  i.e., \  with \  $x= -\rt1\!\pi^2+ 2(1-\rt1\!)\,\pi\,u +u^2 $ (by parity, considering $-\sqrt{\frac{x}{2}} \,$ instead of $\sqrt{\frac{x}{2}} \,$ would yield the same). \parn
For that we need the following approximations. \parsn 
{\small 
$$ \ch\!\big[(1-\rt1\!){\ts\frac{\pi}{2}}\big] = -\rt1\! \sh\!{\ts\frac{\pi}{2}} \,\; ; \  \cos\!\big[(1-\rt1\!){\ts\frac{\pi}{2}}\big] = \rt1\! \sh\!{\ts\frac{\pi}{2}}\;\, ; \  
 \sh\!\big[(1-\rt1\!){\ts\frac{\pi}{2}}\big] = -\rt1\! \ch\!{\ts\frac{\pi}{2}} \,\; ; $$
 $$  \sin\!\big[(1-\rt1\!){\ts\frac{\pi}{2}}\big] =\, \ch\!{\ts\frac{\pi}{2}}\;\, ; \quad  \ch\!\sqrt{\ts\frac{x}{2}} = -\rt1\!\sh\!{\ts\frac{\pi}{2}}\, \big(1 + u\, \coth\!{\ts\frac{\pi}{2}} + \5u^2 +\O(u^3)\big) \,;$$
$$ \sh\!\sqrt{\ts\frac{x}{2}} = -\rt1\!\ch\!{\ts\frac{\pi}{2}}\, \big(1 + u\, \th\!{\ts\frac{\pi}{2}} + \5u^2 +\O(u^3)\big)\, ; \  \sin\!\sqrt{\ts\frac{x}{2}} = \ch\!{\ts\frac{\pi}{2}}\, \big(1 + u\rt1 \th\!{\ts\frac{\pi}{2}} - \5u^2 +\O(u^3)\big) \,;$$    
$$ \cos\!\sqrt{\ts\frac{x}{2}} = \rt1\!\sh\!{\ts\frac{\pi}{2}}\, \big(1 + u\rt1 \coth\!{\ts\frac{\pi}{2}} - \5u^2 +\O(u^3)\big) \,; $$
$$ \ch\!\sqrt{2x} - \cos\!\sqrt{2x}\, = 2\big(\ch\!^2\sqrt{\ts\frac{x}{2}}-\cos^2\!\sqrt{\ts\frac{x}{2}}\big) = 2(\rt1\!-1)\, u\, \sh\pi -4u^2 \ch\pi +\O(u^3) \; ; $$    
$$ \ch\!\sqrt{\ts\frac{x}{2}} + \cos\!\sqrt{\ts\frac{x}{2}}\, = -(1+\rt1\!)\,u\,\ch\!{\ts\frac{\pi}{2}} - u^2\rt1\!\sh\!{\ts\frac{\pi}{2}} +\O(u^3)\;;$$
$$ \sh\!\sqrt{\ts\frac{x}{2}} -\sin\!\sqrt{\ts\frac{x}{2}}\, = -\rt1\!\ch\!{\ts\frac{\pi}{2}}\, \big((1-\rt1\!) + 2u\,\th\!{\ts\frac{\pi}{2}} + \5(1+\rt1\!)\, u^2+\O(u^3)\big) \;; $$
$$ \sh\!\sqrt{\ts\frac{x}{2}} +\sin\!\sqrt{\ts\frac{x}{2}}\, = -\rt1\!\ch\!{\ts\frac{\pi}{2}}\, \big((1+\rt1\!) + \5(1-\rt1\!)\, u^2+\O(u^3)\big) \, .  $$ } 
Therefore on the one hand we have 
$$ \left[{\frac{2\,x}{\ch\!\sqrt{2x}  -\cos\!\sqrt{2x}}}\right]^{1/4} = \left[\frac{(\rt1\!-1)\,\sh\pi\, u  -2\,\ch\pi\,u^2 +\O(u^3)}{-\rt1\!\pi^2+ 2(1-\rt1)\pi\,u +u^2}\right]^{-1/4} $$ 
$$ = \left[-\,\frac{1+\rt1}{\pi^2}\times\frac{\sh\pi\, u +(1+\rt1)\,\ch\pi\,u^2 +\O(u^3)}{1+ 2(1+\rt1)\,u/\pi + \rt1 u^2/\pi^2}\right]^{-1/4} $$ 
$$  =\, \sqrt{\pi}\, 2^{-1/8}\, e^{3\rt1\!\pi/8}\, (u\, \sh\pi)^{-1/4} \Big[1 - {\ts\frac{1+\rt1}{4}} \big(\th\pi - {\ts\frac{2}{\pi}}\big)\,u +\O(u^2)\Big]  $$ 
$$  =\, \sqrt{\pi}\, 2^{-1/8}\, e^{3\rt1\!\pi/8}\, (u\, \sh\pi)^{-1/4} \,\exp\!\Big[- {\ts\frac{1+\rt1}{4}} \big(\th\pi - {\ts\frac{2}{\pi}}\big)\,u +\O(u^2)\Big] . $$ 
Then on the other hand we get\,: 
{
$$ \frac{\rt1\! x}{f(0,x)}\, =\,  {\ts \frac{\big[-\rt1\!\pi^2+ 2(1-\rt1)\pi\,u +u^2\big]\big[(1-\rt1\!)\frac{\pi}{2} + u\big] \big[-(1+\rt1\!)\,u\,\ch\!{\frac{\pi}{2}} - u^2\rt1\!\sh\!{\frac{\pi}{2}} +\O(u^3)\big]}
{-2\rt1\!\sh\!{\frac{\pi}{2}}\, u+(1-\rt1\!)\,\ch\!{\frac{\pi}{2}}\, u^2+\O(u^3)-\rt1\!\big[(1-\rt1\!)\frac{\pi}{2} + u\big] \big[-(1+\rt1\!)\,u\,\ch\!{\frac{\pi}{2}} - u^2\rt1\!\sh\!{\frac{\pi}{2}}\big]} } $$ }
{ 
$$ =  {\ts\frac{\big[\pi^2 + 2(1+\rt1)\pi\,u +\rt1\! u^2\big]\big[(1+\rt1\!)\frac{\pi}{2} + \rt1\!u\big] \big[(1+\rt1\!)\,\ch\!{\frac{\pi}{2}} + u\rt1\!\sh\!{\frac{\pi}{2}} +\O(u^2)\big]}
{-2\rt1\!\sh\!{\frac{\pi}{2}} + (1-\rt1\!)\,\ch\!{\frac{\pi}{2}}\, u + \O(u^2) + \big[(1+\rt1\!)\frac{\pi}{2} + \rt1\!u\big] \big[(1+\rt1\!)\,\ch\!{\frac{\pi}{2}} + u\rt1\!\sh\!{\frac{\pi}{2}}\big]} } $$ }
$$ =  {\frac{\pi^2 \big[ \pi\,\ch\!{\frac{\pi}{2}} + (1+\rt1\!)\big(3\,\ch\!{\frac{\pi}{2}}+ {\frac{\pi}{2}}\,\sh\!{\frac{\pi}{2}}\big)\, u +\O(u^2)\big]} {\big(\pi\,\ch\!{\frac{\pi}{2}}-2\,\sh\!{\frac{\pi}{2}}\big) + (1+\rt1\!)\, {\frac{\pi}{2}}\,\sh\!{\frac{\pi}{2}}\, u +\O(u^2) } } $$ 
$$ =  {\frac{\pi^3} {\pi -2\,\th\!{\frac{\pi}{2}}}} \times \left[1+ (1+\rt1\!)\, \frac{3\pi -6\,\th\!{\frac{\pi}{2}} -\pi\,\th\!^2{\frac{\pi}{2}}}{\pi \big(\pi -2\,\th\!{\frac{\pi}{2}}\big)}\, u + \O(u^2)\right] , \vspace{-2mm} $$ 
so that 
$$ \sqrt{\frac{2\pi\rt1\! x}{f(0,x)}}\, = \frac{\sqrt{2}\,\pi^2 }{\sqrt{\pi -2\,\th\!{\frac{\pi}{2}}}}\times \exp\!\left[(1+\rt1\!)\, \frac{3\pi -6\,\th\!{\frac{\pi}{2}} -\pi\,\th\!^2{\frac{\pi}{2}}}{2\pi \big(\pi -2\,\th\!{\frac{\pi}{2}}\big)}\, u + \O(u^2)\right] . $$ 
Moreover, according to Lemma \ref{lem.ProlAnalEspCompl} and setting $\,v:= 2\pi (1-\rt1)\, u +u^2\,$  we have 
{\small 
$${  \f(x) = \frac{1}{2}\, \arctg\sqrt{\frac{\sqrt{\pi^4+v^2}-\pi^2}{\sqrt{\pi^4+v^2}+\pi^2}}\, - \frac{1}{2} \,\arctg\!\left(\tg\!\sqrt{\frac{\sqrt{\pi^4+v^2}-\pi^2}{2}}\times \coth\!\sqrt{\frac{\sqrt{\pi^4+v^2}+\pi^2}{2}}\,\right) } $$ 
$$ = \frac{1}{2}\, \arctg\!\Big[\frac{v}{2\pi^2} -\frac{v^3}{8\pi^6}+\O(v^5)\Big] - \frac{1}{2}\,\arctg\!\!\left[\tg\!\bigg[\frac{v}{2\pi} -\frac{v^3}{16\pi^5}+\O(v^5)\bigg] \coth\!\bigg[\pi +\frac{v^2}{8\pi^3} \if{- \frac{5\,v^4}{128\pi^7}}\fi 
+ \O(v^4)\bigg]\right] $$ } 
$$ = \frac{1}{2}\, \arctg\Big[\frac{v}{2\pi^2} -\frac{v^3}{8\pi^6}+\O(v^5)\Big] - \frac{1}{2}\,\arctg\!\left[ \frac{v}{2\pi}\,\coth\pi +\frac{(2\pi^2-3)\,\ch\pi\,\sh\pi - 3\pi}{48\,\pi^5\,\sh\!^2\pi}\, v^3 +\O(v^5)\right] $$
$$ =\, \frac{1-\pi\,\coth\pi}{4\pi^2}\,v - \frac{(2\pi^2-3)\,\pi\,\ch\pi\,\sh\pi - 3\pi^2 + (8- 2\pi^3)\,\sh\!^2\pi}{96\,\pi^6\,\sh\!^2\pi}\, v^3 + \O(v^5) $$
$$ =\, \frac{1\!-\pi\,\coth\pi}{2\pi}(1-\rt1\!)\, u + \frac{1\!-\pi\,\coth\pi}{4\pi^2}\, u^2 \vspace{-2mm} $$
$$ +\, \frac{(2\pi^2-3)\,\pi\,\ch\pi\,\sh\pi - 3\pi^2 + (8- 2\pi^3)\,\sh\!^2\pi}{6\,\pi^3\,\sh\!^2\pi}(1+\rt1\!)\, u^3 + \O(v^4) .$$
Hence, according to Proposition \ref{pro.exprAvChgCont}, for \  $x= -\rt1\!\pi^2+ 2(1-\rt1)\pi\,u +u^2 $ \  we have\,: 
$$ \Phi(x) = \, e^{\rt1 \f(x)}\times  \left[{\frac{2\,x}{\ch\!\sqrt{2x}  -\cos\!\sqrt{2x}}}\right]^{1/4} \times \sqrt{\frac{2\pi\rt1\! x}{f(0,x)}}\, \raise1.9pt\hbox{} $$ 
$$ = {\ts\frac{2^{3/8}\,\pi^{5/2}\, e^{3\rt1\!\pi/8}}{ \sqrt{\pi -2\,\th\!{\frac{\pi}{2}}}\, (u\, \sh\pi)^{1/4}}} \times \exp\!\left[(1\!+\!\rt1\!) \!\left({\ts \frac{1\!-\pi\,\coth\pi}{2\pi} + \frac{3\pi -6\,\th\!{\frac{\pi}{2}} -\pi\,\th\!^2{\frac{\pi}{2}}}{2\pi \left(\pi -2\,\th\!{\frac{\pi}{2}}\right)} - \frac{\pi\,\th\pi - 2}{4\pi} }\right)\! u + \O(u^2)\right] $$ 
{\small 
$$ = \frac{2^{3/8}\,\pi^{5/2}\, e^{3\rt1\!\pi/8}}{ \sqrt{\pi -2\,\th\!{\frac{\pi}{2}}}\, (u\, \sh\pi)^{1/4}}\, \exp\!\left[(1+\rt1\!) \!\left({\ts \frac{4 - 2\pi\,\coth\pi - \pi\,\th\pi}{4\pi} + \frac{3\pi -6\,\th\!{\frac{\pi}{2}} -\pi\,\th\!^2{\frac{\pi}{2}}}{2\pi \left(\pi -2\,\th\!{\frac{\pi}{2}}\right)} } \right) u + \O(u^2)\right] . $$ }\parn 
Equivalently, taking $\,h= 2(1-\rt1)\pi\,u +u^2\,\LRa u= \frac{1+\rt1}{4\pi}\big(1-\frac{\rt1h}{8\pi^2}+\O(h^2)\big) h\,$, this yields
$$ \Phi\big(\pi^2,\, h\big) = \, \frac{2^{3/4}\,\pi^{11/4}\, e^{5\rt1\!\pi/16}}{ \sqrt{\pi -2\,\th\!{\frac{\pi}{2}}}\, ( \sh\pi)^{1/4}}\times \bigg(\frac{1+\frac{\rt1h}{32\pi^2}+\O(h^2)}{h^{1/4}}\bigg) \times   \vspace{-2mm}  $$ 
$$ \qquad \times \exp\!\left[ \rt1 \!\left({\ts \frac{4 - 2\pi\,\coth\pi - \pi\,\th\pi}{8\,\pi^2} + \frac{3\pi -6\,\th\!{\frac{\pi}{2}} -\pi\,\th\!^2{\frac{\pi}{2}}}{4\pi^2 \left(\pi -2\,\th\!{\frac{\pi}{2}}\right)} } \right) h + \O(h^2)\right] =  $$  
$$ \frac{2^{3/4}\,\pi^{11/4}\, e^{5\rt1\!\pi/16}}{ \sqrt{\pi -2\,\th\!{\frac{\pi}{2}}}\, ( \sh\pi)^{1/4}}\, h^{-1/4} \times  \exp\!\left[ \rt1\! \!\left[{\ts \frac{5/4 - 2\pi\,\coth\pi - \pi\,\th\pi}{8\,\pi^2} + \frac{3\pi -6\,\th\!{\frac{\pi}{2}} -\pi\,\th\!^2{\frac{\pi}{2}}}{4\pi^2 \left(\pi -2\,\th\!{\frac{\pi}{2}}\right)} }\right]\! h + \O(h^2)\right]\! . \;\; \diamond $$  
\parb 

     For the analogue of Proposition \ref{pro.ChangContz=0>yb}, let us consider the new contour 
\begin{equation} \label{df.nContourb'} 
\Gamma := \Big[0\,, -\rt1\pi^2\Big]\, \, \bigcup\, \Big(-\rt1\pi^2 + \R_+\Big) . \vspace{-2mm} 
\end{equation}
\bpro \label{pro.ChangContz=0>yb'} \  For any $\;y>0\,$ and small $\,\e >0\,$ we have 
$$ \II_\e(y,0) = \, \exp\!\left[{ -{\frac{y-\e}{\e^2}} }\, \pi^2 \right] \times 2\; \Re\Bigg\{\! \int_{0}^\infty\! \exp\!\left[{ {\frac{\e - y}{\e^2}} }\rt1\! x \right]  \Phi\big(\pi^2\hbox{,}\, x\big)\, dx \Bigg\} \, ,  $$ 
where \  $\Phi(\chi, x)  \equiv \Phi\big(x-\rt1\!\chi\big)\,$ was defined in (\ref{df.PhiProlAnal}). 
\epro   
\ub{Proof} \quad The change of contour (\ref{df.nContourb'}) is justified as for Proposition \ref{pro.ChangContz=0>yb}, using  $\,\left| \Phi\big(\chi,R\big)\right| = \O\big(R\,e^{-\sqrt{R/8}\,}\big)$ as well.  Performing this change yields\,: 
$$ \int_{0}^\infty \exp\!\left[{ {\frac{\e - y}{\e^2}} }\rt1\! x \right] \Phi(x)\, dx \, =  \int_{\Gamma} \exp\!\left[{ {\frac{\e - y}{\e^2}} }\rt1\! x \right] \Phi(x)\, dx\, =  $$
$$  \exp\!\left[{{\frac{\e-y}{\e^2}} }\,\pi^2\right]\! \times\! \int_{0}^\infty \exp\!\left[{ {\frac{\e - y}{\e^2}} }\rt1\! x \right]  \Phi\big(\pi^2\hbox{,}\, x\big)\, dx\, - \rt1\!\! \int_{0}^{\pi^2} \exp\!\left[{ {\frac{\e-y}{\e^2}} }\,x \right] \Phi(x,0)\, dx\,.  $$ 
Then for any $\,\chi>0\,$, we have $\,\f(\chi,0)=0\,$ by Lemma \ref{lem.ProlAnalEspCompl}, \   and $\,f(\chi,0)\in\R$ by (\ref{df.fnf}), and then $\,\Phi(\chi,0)\in\R$ by (\ref{df.PhiProlAnal}). Thus the second term of the right hand side 
disappears when taking the real part.  $\diamond$ \parb 

The analogue of (\ref{f.CutIntegepsr}) is the following. 
$$ \int_{0}^\infty\! \exp\!\left[{ {\frac{\e - y}{\e^2}} }\rt1\! x \right]  \Phi\big(\pi^2\hbox{,}\, x\big)\, dx\, = \int_{\e^r}^\infty\! \exp\!\left[{ {\frac{\e - y}{\e^{2}}} }\rt1\! x \right]  \Phi\big(\pi^2\hbox{,}\, x\big)\, dx\, $$ 
\begin{equation}  \label{f.CutIntegepsr>}
+\, \frac{\e^2}{y-\e} \int_{0}^{\e^{r-2}(y-\e)}\! \exp\!\left[-\rt1\! x \right]\times \Phi\left(\pi^2\hbox{,}\, {\frac{\e^2\, x}{y-\e}}\right)\, dx\, .  \vspace{2mm} 
\end{equation} 

We now apply Lemma \ref{lem.AnalPhiprSing>}  to the last integral in (\ref{f.CutIntegepsr>}), to obtain the following analogue of Proposition \ref{pro.IntegPartPrinc}. 
\bpro \label{pro.IntegPartPrinc>0}  \  For  $\,y> 0\,$, $\e>0\,$ and $\,0< r<2\,$, we have 
$$ \Re\Bigg\{\! \int_{0}^{\e^{r-2}(y-\e)} \exp\!\left[-\rt1\! x \right]\times \Phi\left(\pi^2 \if{+{\ts\frac{\e^2}{4 (y-\e)}}\,\raise0.5pt}\fi \hbox{,}\, {\frac{\e^2\,x}{y-\e}}\right)\, dx\Bigg\} $$
$$ =\, \frac{2^{3/4}\,\pi^{11/4}\, y^{1/4}}{ \sqrt{\pi -2\,\th\!{\frac{\pi}{2}}}\,(\sh\pi)^{1/4}}\times \frac{1+\O\big(\e^r + \e^{(2-r)/4}\big)}{\sqrt{\e}} \int_0^\infty\, \frac{\sin\!\left(\frac{3\pi}{16} + x\right)}{x^{1/4}}\, dx\,  $$ 
$$ =\, \frac{2^{3/4}\,\pi^{11/4}\, y^{1/4}}{ \sqrt{\pi -2\,\th\!{\frac{\pi}{2}}}\,(\sh\pi)^{1/4}}\times \frac{1+\O\big(\e^{2/5}\big)}{\sqrt{\e}} \int_0^\infty\, \frac{\sin\!\left(\frac{3\pi}{16} + x\right)}{x^{1/4}}\, dx \quad \hbox{by eventually taking $\,r= 2/5$.} $$ 
 \epro 
\ub{Proof} \quad  Applying apply Lemma \ref{lem.AnalPhiprSing>} we indeed have\,: 
$$ \Re\Bigg\{\! \int_{0}^{\e^{r-2}(y-\e)} \exp\!\left[-\rt1\! x \right]\times \Phi\big(\pi^2 \if{+{\ts\frac{\e^2}{4 (y-\e)}}\,\raise0.5pt}\fi \hbox{,}\, {\ts\frac{\e^2}{y-\e}}\,x\big)\, dx\Bigg\} \, 
 =\, \frac{2^{3/4}\,\pi^{11/4}}{ \sqrt{\pi -2\,\th\!{\frac{\pi}{2}}}\, ( \sh\pi)^{1/4}}\, \times \vspace{-2mm} $$
$$ \qquad \times \int_{0}^{\e^{r-2}(y-\e)} \left({\frac{\e^2 x}{y-\e}} \if{- {\ts\frac{\rt1 \e^2}{4 (y-\e)}}}\fi \right)^{-1/4} \Re\bigg\{\exp\!\left[{-\rt1 x+\rt1\frac{5\pi}{16}} \if{ + \rt1 C\, \frac{\e^2 x}{y-\e} \left({\ts\frac{\e^2 x}{y-\e}} - {\ts\frac{\rt1 \e^2}{4 (y-\e)}}\right)}\fi +\O(\e^{r})\right]\! \bigg\}\, dx $$ 
$$ =\, \frac{2^{3/4}\,\pi^{11/4}\, (y-\e)^{1/4}}{ \sqrt{\pi -2\,\th\!{\frac{\pi}{2}}}\,(\sh\pi)^{1/4}}\times \frac{1+\O(\e^r)}{\sqrt{\e}} \int_{0}^{\e^{r-2}(y-\e)}  \cos\!\left[ x-\frac{5\pi}{16}\right]   \,x^{-1/4} \, dx $$
$$ =\, \frac{2^{3/4}\,\pi^{11/4}\, y^{1/4}}{ \sqrt{\pi -2\,\th\!{\frac{\pi}{2}}}\,(\sh\pi)^{1/4}}\times \frac{1+\O(\e^r+\e)}{\sqrt{\e}} \left[ \int_{0}^{\infty}  \sin\!\left[ x+\frac{3\pi}{16}\right]   \,x^{-1/4} \, dx + \O\big(\e^{(2-r)/4}\big)\right] $$ 
$$ =\, \frac{2^{3/4}\,\pi^{11/4}\, y^{1/4}}{ \sqrt{\pi -2\,\th\!{\frac{\pi}{2}}}\,(\sh\pi)^{1/4}}\times \frac{1+\O\big(\e^r + \e^{(2-r)/4}\big)}{\sqrt{\e}} \int_{0}^{\infty}  \sin\!\left[ x+\frac{3\pi}{16}\right]   \,x^{-1/4} \, dx \,  $$ 
by the following (easy analogue of Proposition \ref{pro.cstte>0}). $\;\diamond$ 
\blem  \label{lem.cstte>0} \  The constant \quad ${\ds \sigma' := \int_0^\infty\, \frac{\sin\!\left(\frac{3\pi}{16} + x\right)}{x^{1/4}}\, dx\,}$ \   is positive ($\,> \frac{1}{10}\,$). 
\elem
\ub{Proof} \quad  We indeed have \vspace{-1mm} 
$$  \sigma'  = \int_{\frac{3\pi}{16} }^\infty \frac{\sin x\; dx}{\left(x-\frac{3\pi}{16}\right)^{1/4}}\, = \int_{\frac{3\pi}{16}}^{\frac{\pi}{2}} \frac{\sin x\; dx}{\left(x-\frac{3\pi}{16}\right)^{1/4}} + \int_{\frac{\pi}{2}}^\pi\, \frac{\sin x\; dx}{\left(x-\frac{3\pi}{16}\right)^{1/4}} + \int_{\pi}^{\frac{3\pi}{2}}\, \frac{\sin x\; dx}{\left(x-\frac{3\pi}{16}\right)^{1/4}} \vspace{-1mm} $$ 
$$ + \int_{\frac{3\pi}{2}}^{2\pi} \frac{\sin x\; dx}{\left(x-\frac{3\pi}{16}\right)^{1/4}} + \sum_{n\ge 1} \left[
\int_{2n\pi}^{(2n+1)\pi} \frac{\sin x\; dx}{\left(x-\frac{3\pi}{16}\right)^{1/4}} + \int_{(2n+1)\pi}^{2(n+1)\pi} \frac{\sin x\; dx}{\left(x-\frac{3\pi}{16}\right)^{1/4}} \right] $$ 
$$ > {\ts  \left(\frac{5\pi}{16}\right)^{-1/4}\! \cos\!\left(\frac{3\pi}{16}\right) + \left(\frac{13\pi}{16}\right)^{-1/4}\! - \left(\frac{13\pi}{16}\right)^{-1/4} - \left(\frac{21\pi}{16}\right)^{-1/4} =\,  \left(\frac{16}{5\pi}\right)^{1/4} \left[ \cos\!\left(\frac{3\pi}{16}\right) - \left(4+\frac{1}{5}\right)^{-1/4}  \right] } $$  
$$ >\,  {\ts \left(\frac{16}{5\pi}\right)^{1/4} \left[ 1- \5 \left(\frac{3\pi}{16}\right)^2 - 2^{-1/2}\right] } > \frac{1}{10}\, \raise1.9pt\hbox{.} \;\;\diamond $$ 
\parm 
  
  We now control the remaining integral of (\ref{f.CutIntegepsr>}), namely  ${\ds \int_{\e^r}^\infty\! \exp\!\left[{\ts{\frac{\e - y}{\e^{2}}} }\rt1\! x \right]  \Phi\big(\pi^2\hbox{,}\, x\big)\, dx\,}$. For that we hall resort to integration by parts, so that we need to estimate $\,{\ds\frac{d \Phi}{dx}\big(\pi^2, x\big)}$. This will be made by using Lemma \ref{lem.LiftSqRtxPsi}. The following analogue of Lemma \ref{lem.ControlPhi'} will actually be sufficient. 
\blem \label{lem.ControlPhi'>} \  For any $\,y> 0\,$ and $\,x\ge \e^{r}$, uniformly as $\,\e \sea 0$  \  we have 
$$ \big| \Phi\big(\pi^2, x\big)\big| +\bigg| \frac{d \Phi}{dx}\big(\pi^2, x\big)\bigg| =\, 1_{\{\e^{r}\le x< 1\}}\, \O(x^{-1/4}) + 1_{\{x\ge1\}}\, \O\!\left( x\, e^{-\sqrt{x/8}}\,\right) . \vspace{-3mm} $$ 
\elem
\ub{Proof} \quad  Lemma \ref{lem.LiftSqRtxPsi} provides the uniform estimate\,: ${\ds \big| \Phi\big(\pi^2, x\big) \big|= \O\!\left[x\times \Big(\ch\!\sqrt{x/2}\Big)^{-1/2}\right]}$  ${\ds = \O\!\left[ x\, e^{-\sqrt{x/8}}\,\right]}$ for large $\,x\,$. Moreover Lemma \ref{lem.AnalPhiprSing>} ensures that \  ${\ds \Phi\big(\pi^2, \,x\big) = \O(x^{-1/4}) }\,$ for small $\,x\,$. This settles the case of $\Phi$.   \  
Then using Lemma \ref{lem.LiftSqRtxPsi} again and (\ref{f.racinexkhi}), as in the proof of Lemma \ref{lem.ControlPhi'} we obtain\,:  \vspace{-1mm} 
$$ \frac{d \log\Phi}{dx}(\pi^2, x)\, 
=\, \rt1 \frac{d(\f+\frac{1}{2}\tilde\f)}{dx}(\pi^2, x)  +\O(1)\,. $$ 
Note that the present case $\,\chi=\pi^2$ is simpler than the case $\,\chi\approx -4\pi^2$ of Lemma \ref{lem.ControlPhi'}, since now we have $\,a\ge \pi\,$, hence $\,\ch a\ge \ch\pi\,$, so that the value of $\,b\,$ does not matter so much in the estimates needed here. \  
Using again  \  
${\ds \, \frac{\partial a}{\partial x} = \frac{b/2}{a^2+b^2}\,\raise1pt\hbox{,}\, \frac{\partial b}{\partial x} = \frac{a/2}{a^2+b^2}\,\raise1pt\hbox{}}$ \  and Lemma \ref{lem.ProlAnalEspCompl},  we thus have \vspace{-1mm} 
$$ \frac{d \f}{dx}(\pi^2, x) =\,\5\, \frac{d }{dx}\big[\arctg\big({\ts\frac{b}{a}}\big) -\f_a(b)\big] = \frac{\pi^2}{4\,(\pi^4+x^2)} - \frac{a\,\sh\!(2a)- b\, \sin(2b)}{4(a^2+b^2) (\ch\!(2a)-\cos (2b))} = \O(1). $$
Then by Lemma \ref{lem.LiftSquaRoot} (again as in the proof of Lemma \ref{lem.ControlPhi'}) and (\ref{df.fnf}) we have 
$$ \frac{d \tilde\f}{dx}(\pi^2, x) =\, 
\frac{\pi^2}{\pi^4+x^2} -  \frac{d }{dx}\, \arctg \!\left[ { \frac{2\, (a\, \sin b - b\,\sh a)}{2(a\, \sh a + b\, \sin b) - (a^2+b^2)\, (\ch a + \cos b)}}\right] $$ 
{\small 
$$ = \frac{\pi^2}{\pi^4+x^2} + { \frac{\O\big((a^2+b^2)^2(\ch a + \cos b)^2\big) }{ (a^2+b^2)^2 (\ch a + \cos b)\big[(a^2+b^2)\, (\ch a + \cos b)-4 (a\, \sh a + b\, \sin b) + 4\, (\ch a - \cos b)\big]}}  $$ }
which is bounded. \  The claim follows. 
\if{ 
\pars \centerline{Ou bien (preuve alternative ; laiss\'e incompl\`ete) :} \pars 
   We compute $\,{\ds\frac{\Phi'}{\Phi}\,}$ in the domain $\{\Re\, x>0\}$.  According to Proposition \ref{pro.exprAvChgCont}, we have\,: 
$$ \frac{\Phi'}{\Phi}(x)\, =\, \rt1 \f'(x) + \frac{3}{4x} - \frac{1}{4\sqrt{2x}}\times \frac{\sh\!\sqrt{2x} +\sin\!\sqrt{2x}}{\ch\!\sqrt{2x} - \cos\!\sqrt{2x}} - \frac{f'(x)}{2\,f(x)}\, \raise1.9pt\hbox{}  $$
$$ =\, - \frac{\rt1}{2\sqrt{2x}}\times \frac{\sh\!\sqrt{2x} - \sin\!\sqrt{2x}}{\ch\!\sqrt{2x} - \cos\!\sqrt{2x}}  + \frac{3}{4x} - \frac{1}{4\sqrt{2x}}\times \frac{\sh\!\sqrt{2x} + \sin\!\sqrt{2x}}{\ch\!\sqrt{2x} - \cos\!\sqrt{2x}} + \frac{1}{4x} \, \raise1.9pt\hbox{}  $$ 
{\small 
$$ +\, \frac{1}{4\sqrt{2x}} \left[\frac{\sh\!\sqrt{\frac{x}{2}} - \sin\!\sqrt{\frac{x}{2}}}{\ch\!\sqrt{\frac{x}{2}} + \cos\!\sqrt{\frac{x}{2}}} - \frac{ \ch\!\sqrt{\frac{x}{2}} - \cos\!\sqrt{\frac{x}{2}} - \rt1 (\sh\!\sqrt{\frac{x}{2}} - \sin\!\sqrt{\frac{x}{2}}\,)\sqrt{\frac{x}{2}} }{ \sh\!\sqrt{\frac{x}{2}} - \sin\!\sqrt{\frac{x}{2}} + \rt1\!\big(\sh\!\sqrt{\frac{x}{2}} + \sin\!\sqrt{\frac{x}{2}} - (\ch\!\sqrt{\frac{x}{2}} + \cos\!\sqrt{\frac{x}{2}}\,)\sqrt{\frac{x}{2}}\,\big)}\right] $$ } 
%
$$ =\,  \frac{1}{x} - \frac{1}{2\sqrt{2x}}\times \frac{\ch\!\sqrt{\frac{x}{2}} \sin\!\sqrt{\frac{x}{2}} + \sh\!\sqrt{\frac{x}{2}} \cos\!\sqrt{\frac{x}{2}}+\rt1\! \big(\sh\!\sqrt{2x} - \sin\!\sqrt{2x}\,\big)}{\ch\!\sqrt{2x} - \cos\!\sqrt{2x}} $$
$$ -\,  \frac{1}{4\sqrt{2x}}\times \frac{ \ch\!\sqrt{\frac{x}{2}} - \cos\!\sqrt{\frac{x}{2}} - \rt1 (\sh\!\sqrt{\frac{x}{2}} - \sin\!\sqrt{\frac{x}{2}}\,)\sqrt{\frac{x}{2}} }{ \sh\!\sqrt{\frac{x}{2}} - \sin\!\sqrt{\frac{x}{2}} + \rt1\!\big(\sh\!\sqrt{\frac{x}{2}} + \sin\!\sqrt{\frac{x}{2}} - (\ch\!\sqrt{\frac{x}{2}} + \cos\!\sqrt{\frac{x}{2}}\,)\sqrt{\frac{x}{2}}\,\big)} \, \raise1.9pt\hbox{.} $$ 
Moreover, owing to the proof of Lemma \ref{lem.AnalPhiprSing>} we also have 
$$ \sh\!\sqrt{2x} - \sin\!\sqrt{2x}\, =\, \sh\pi \big(-1-\rt1+2(\rt1\!-1) u^2 +\O(u^3)\big) \; \vspace{-2mm} $$    
and
$$ \ch\!\sqrt{\ts\frac{x}{2}} - \cos\!\sqrt{\ts\frac{x}{2}}\, = -\rt1\!\sh\!{\ts\frac{\pi}{2}}\, \big(2 + u(1+\rt1\!)\,\coth\!{\ts\frac{\pi}{2}} +\O(u^3)\big) \, . $$
Hence   \vspace{-1mm} 
$$ \frac{\Phi'}{\Phi}\big(-\rt1\!\pi^2+ 2(1-\rt1)\pi\,u +u^2\big) =\, \frac{\rt1}{\pi^2}\big(1-2\pi(1+\rt1\!)\,u +\O(u^2)\big) 
$$ 
$$-\, \frac{\frac{1-\rt1}{2}\,\sh\pi + 2u\,\sh^2\frac{\pi}{2} + \frac{1+\rt1}{2}\,u^2\sh\pi +(1-\rt1)\sh\pi -2(1+\rt1)u^2\sh\pi +\O(u^3)}{2 \big[ (1-\rt1\!)\pi+2u\big] \big[2(\rt1\!-1)\, u\, \sh\pi -4u^2 \ch\pi +\O(u^3)\big] } $$ 
{
$$ +\, {\ts\frac{ \rt1\!\sh\!{\ts\frac{\pi}{2}}\, \big(2 + u(1+\rt1\!)\,\coth\!{\ts\frac{\pi}{2}} \big) +\ch\!{\ts\frac{\pi}{2}}\, \big((1-\rt1\!) + 2u\,\th\!{\ts\frac{\pi}{2}} + \5(1+\rt1\!)\, u^2\big)\big[(1-\rt1\!)\frac{\pi}{2} + u\big] +\O(u^3)}{ 4 \big[ (1-\rt1\!)\pi+2u\big] \big[ -2\rt1u\, \sh\!{\ts\frac{\pi}{2}} + (1-\rt1\!)\, u^2 \ch\!{\ts\frac{\pi}{2}} + \big((\rt1\!-1)\,u\,\ch\!{\ts\frac{\pi}{2}} + u^2 \sh\!{\ts\frac{\pi}{2}}\,\big)\big[(1-\rt1\!)\frac{\pi}{2} + u\big] \big] +\O(u^3)} } $$ }
$$ = \frac{\rt1}{\pi^2}\big(1-2\pi(1+\rt1\!)\,u +\O(u^2)\big) $$
$$ -\, \frac{3(1-\rt1)\,\sh\pi + 4u\,\sh\!^2\frac{\pi}{2} - 3(1+\rt1)\,u^2\sh\pi  +\O(u^3)}{ 16\,u \big[\pi\, \sh\pi\rt1 + (1-\rt1\!)\, (\pi\,\ch\pi - \sh\pi)\, u +\O(u^2)\big] } $$ 
{
$$ +\, {\frac{ \rt1\!\big[2\, \sh\!{\ts\frac{\pi}{2}} -\pi\,\ch\!{\ts\frac{\pi}{2}}\,\big]+ (1-\rt1\!)\pi\,\sh\!{\ts\frac{\pi}{2}}\,u +\big[2\, \sh\!{\ts\frac{\pi}{2}} +{\ts\frac{\pi}{2}}\,\ch\!{\ts\frac{\pi}{2}}\,\big] u^2 +\O(u^3)}{ 4 \big[ (1-\rt1\!)\pi+2u\big] \big[ \rt1\!\big[\pi\,\ch\!{\ts\frac{\pi}{2}}-2\,\sh\!{\ts\frac{\pi}{2}}\big] u + (1-\rt1\!)\big[\frac{\pi}{2}\,\sh\!{\ts\frac{\pi}{2}} \big] u^2 \big] +\O(u^3) } } $$ }
$$ = \frac{\rt1}{\pi^2}\big(1-2\pi(1+\rt1\!)\,u +\O(u^2)\big) $$
$$ +\, \frac{3(1+\rt1) + 2\rt1 \th\frac{\pi}{2}\,u + 3(1-\rt1)\, u^2  +\O(u^3)}{ 16\pi\,u \big[1 - (1+\rt1\!)\, (\coth\pi - 1/\pi)\, u +\O(u^2)\big] } $$ 
{
$$ -\, {\frac{ (1+\rt1\!)+ \rt1\frac{2\pi\,\sh\!{\frac{\pi}{2}}}{\pi\,\ch\!{\frac{\pi}{2}}-2\, \sh\!{\frac{\pi}{2}}}\,u - (1-\rt1\!)\frac{2\, \sh\!{\frac{\pi}{2}} +{\frac{\pi}{2}}\,\ch\!{\frac{\pi}{2}}}{\pi\,\ch\!{\frac{\pi}{2}}-2\, \sh\!{\frac{\pi}{2}}}\, u^2 +\O(u^3)}{ 8\pi \,u \big[ 1+ (1+\rt1\!) \frac{\ch\!{\frac{\pi}{2}}-(\pi/2+2/\pi)\,\sh\!{\frac{\pi}{2}} }{\pi\,\ch\!{\frac{\pi}{2}}-2\, \sh\!{\frac{\pi}{2}}}\, u +\O(u^2)\big] } } $$ }
$$ = \frac{\rt1}{\pi^2}\big(1-2\pi(1+\rt1\!)\,u +\O(u^2)\big) $$
$$ +\, \frac{(1+\rt1) + 2\rt1\!\big(\frac{4\,\ch\pi -1}{\sh\pi }-\frac{3}{\pi} - \frac{(3\pi+4/\pi)\,\sh\!{\frac{\pi}{2}}-2\,\ch\!{\frac{\pi}{2}}}{\pi\,\ch\!{\frac{\pi}{2}}-2\, \sh\!{\frac{\pi}{2}}}\big) u +\O(u^2)}{ 16\pi\,u  } $$ 
$$ =  \frac{1+\rt1}{16\pi\,u} + \rt1\!\left[\frac{4\,\ch\pi -1}{8\pi\,\sh\pi } + \frac{5}{8\pi^2} - \frac{(3\pi+4/\pi)\,\sh\!{\frac{\pi}{2}}-2\,\ch\!{\frac{\pi}{2}}}{8\pi\,\big(\pi\,\ch\!{\frac{\pi}{2}}-2\, \sh\!{\frac{\pi}{2}}\,\big)}\right]+\O(u)  $$ 
$$ =  \frac{1+\rt1}{16\pi\,u} + \rt1  \frac{(\ch\pi-1)\big[(\pi-14/\pi)\,\ch\!{\frac{\pi}{2}}-\sh\!{\frac{\pi}{2}}\big] - 3\pi\,\ch\!{\frac{\pi}{2}} + 4\sh\!{\frac{\pi}{2}}+4}{8\pi\,\sh\pi \big(\pi\,\ch\!{\frac{\pi}{2}}-2\, \sh\!{\frac{\pi}{2}}\,\big)} + \O(u)\, . $$ 
}\fi 
$\;\diamond$ \parm

   Now the analogue of Lemma \ref{lem.ControlIntPhif} is the following. 
\blem \label{lem.ControlIntPhif>} \  We have \quad   ${\ds \int_{\e^{r}}^{\infty} \exp\!\left[{\frac{\e - y}{\e^2}}\rt1 x\right]\!\times  \Phi\big(\pi^2, x\big)\, dx\, =\, \O\big(\e^{2-r/4}\big)\, . \vspace{-0mm}  }$ 
\elem
\ub{Proof} \quad Integrating by parts, since $\Phi(\pi^2,\cdot)$ vanishes at infinity, using Lemma \ref{lem.ControlPhi'>} we have\,: 
$$ \int_{\e^{r}}^{\infty}\! \exp\!\left[{ {\frac{\e - y}{\e^2}} }\rt1\! x \right]  \Phi\big(\pi^2, x\big)\, dx \, $$
$$ =\, {\frac{\e^2\rt1\!}{\e - y}}  \exp\!\left[{ {\frac{\e - y}{\e^{2-r}}} }\rt1\!\right]  \Phi\big(\pi^2, \e^{r}\big) + {\frac{\e^2\rt1\!}{\e - y}} \int_{\e^{r}}^{\infty}\! \exp\!\left[{ {\frac{\e - y}{\e^2}} }\rt1\! x \right]  \frac{d \Phi}{dx}\big(\pi^2, x\big)\, dx $$ 
$$ =\, \O\big(\e^{2-r/4}\big) + \O(\e^{2}) \int_{\e^{r}}^{1} \O(x^{-1/4})\,  dx + \O(\e^{2}) \int_{1}^{\infty} \O\!\left( x\, e^{-\sqrt{x/8}}\,\right) dx\, =\, \O\big(\e^{2-r/4}\big) . \;\;\diamond $$ \parm 

   We can now conclude this section, by the following exact equivalent of the oscillatory integral $\,\II_\e(y,0)$ arising in the case  $\,y> 0=z\,$. 
\bpro  \label{pro.ControlIntPhi>} \    For $\,y> 0\,$,  as $\,\e \sea 0$  \  (with $\,\sigma'$ as in Lemma \ref{lem.cstte>0}) we have 
$$ \Re\Bigg\{\!\int_{0}^\infty\! \exp\!\left[{ {\frac{\e - y}{\e^2}} }\rt1\! x \right] \times   \Phi\!\left(\pi^2, \,x\right)  dx \Bigg\} =\, \frac{2^{3/4}\,\pi^{11/4}\, \sigma'}{ \sqrt{\pi -2\,\th\!{\frac{\pi}{2}}}\,(\sh\pi)^{1/4}}\times \frac{\e^{3/2}}{y^{3/4}} \left(1+O\big(\e^{2/5}\big)\right) \vspace{-0mm}  $$ 
and \vspace{-2mm} 
$$ \II_\e(y,0) = \, \exp\!\left[-\pi^2\,\frac{y-\e}{\e^2}\right]\times \frac{\e^{3/2}}{y^{3/4}} \times \frac{(2\pi)^{11/4}\, \sigma'}{ \sqrt{\pi -2\,\th\!{\frac{\pi}{2}}}\,(\sh\pi)^{1/4}} \times\left(1+\O\big(\e^{2/5}\big)\right) . \vspace{-2mm} $$ 
\epro 
\ub{Proof} \quad  The first claim follows directly from (\ref{f.CutIntegepsr>}) with $\,r = \frac{2}{5} \,$, Lemma \ref{lem.ControlIntPhif>} and Proposition \ref{pro.IntegPartPrinc>0}. With Proposition \ref{pro.ChangContz=0>yb'}, it entails the second one. $\;\diamond$ \parm 

   Propositions \ref{pro.equivIntPe0} and \ref{pro.ControlIntPhi>} together give the following wanted small time equivalent, which is the content of Theorem \ref{th.mainres}$(iii)$ when $\,z=0$. 
\bcor \label{cor.w=0=z<y} \  For $\,y > 0\,$, as $\,\e\sea 0\,$ we have 
$$ p_\e\big(0\, ; (0,y, 0)\big) \sim\, \exp\!\left[-\pi^2\,\frac{y-\e}{\e^2}\right]\times \frac{\e^{-5/2}}{y^{3/4}} \times \frac{(2\pi/\sh\pi)^{1/4}}{ \sqrt{\pi -2\,\th\!{\frac{\pi}{2}}}} \int_0^\infty\, \frac{\sin\!\left(\frac{3\pi}{16} + x\right)}{x^{1/4}}\, dx\, . \vspace{-2mm} $$ 
\ecor    
\vspace{-3mm} 

\subsection{Third sub-case\,: $\, z\not=0$} \label{sec.znot=0<yz} \indf 
   Recall that according to Proposition \ref{pro.exprAvChgCont} we have 
$$ \II_\e(y,z) =\, 2\; \Re \Bigg\{\! \int_{0}^\infty\! \exp\!\left[{ {\frac{\e - y}{\e^2}} }\rt1\! x - \frac{z^2}{2\e^3}\, \frac{\rt1\! x}{f(0,x)}\right] \Phi(x)\, dx \Bigg\} \, . \vspace{-0mm} $$ 

   We shall proceed somehow as in the preceding case in Section \ref{sec.znot=0<y}, $\,f(x)$ being positive at the saddle point, and roughly try to substitute $\,{\ds \frac{y-\e}{\e^2} + \frac{z^2}{2\e^3 f(0,x)}\,}$ for $\,{\ds \frac{y-\e}{\e^2}\,}$\raise1.9pt\hbox{.} \  However the dependence with respect to the non-constant $f(0,x)$ demands additional care. \pars    
   First, since $\,\chi\,f(\chi,0) >0$ by (\ref{df.fnf}), Proposition \ref{pro.ChangContz=0>yb'} at once becomes the following. 
\bpro \label{pro.ChangContz=0>yb'z} \  For any $\;z\not=0\,,\,y\in\R\,$ and small $\,\e >0\,$ we have 
$$ \II_\e(y,z) = \, 2\; \Re\Bigg\{\! \int_{0}^\infty\!  \exp\!\left( -\left[\frac{y-\e}{\e^2} + \frac{z^2}{2\e^3 f(\pi^2,x)}\right]\! \big(\rt1\! x + \pi^2\big)\!\right)\!\times \Phi\big(\pi^2\hbox{,}\, x\big)\, dx \Bigg\} \, .  $$ 
\epro   
\if{ 
We use (\ref{df.fnf}), with $\,\chi = \theta^2\,$ and successively\,: \vspace{-1mm} 
$$ \sqrt{\theta^2+\rt1\!x}\, = \,\theta + \frac{\rt1\!x}{2\theta} + \frac{x^2}{8\,\theta^3} - \frac{\rt1\!x^3}{16\,\theta^5} +\O(|x|^4)\, ; $$ 
$$  \th\!\left[\sqrt{\theta^2+\rt1\!x}\Big/2\right] =\, \th\!\!\left[ \frac{\theta}{2} + \frac{\rt1\!x}{4\,\theta} + \frac{x^2}{16\,\theta^3} - \frac{\rt1\!x^3}{32\,\theta^5} + \O(|x|^4)\right] $$ 
$$ =\, \th{\ts\frac{\theta}{2}} + \frac{\frac{\rt1\!x}{4\,\theta} + \frac{x^2}{16\,\theta^3} - \frac{\rt1\!x^3}{32\,\theta^5} +\O(|x|^4)}{\ch\!^2\frac{\theta}{2}} - \frac{2\,\th\frac{\theta}{2}}{\ch\!^2\frac{\theta}{2}}\! \left[ {\frac{\rt1\!x}{4\,\theta} + \frac{x^2}{16\,\theta^3} - \frac{\rt1\!x^3}{32\,\theta^5} +\O(|x|^4)}\right]^2 $$ 
$$ =\, \th{\ts\frac{\theta}{2}} + \frac{\rt1\!x}{4\,\theta\,\ch\!^2\frac{\theta}{2}} + \frac{1+ 2 \theta\, \th\frac{\theta}{2}}{16\,\theta^3\,\ch\!^2\frac{\theta}{2}}\,x^2 + \O(|x|^3) \, ; $$ 
$$ f(\theta^2, x) \, = \, 1- \frac{ \th\!{\ts\frac{\theta}{2}} + \frac{\rt1\!x}{4\,\theta\,\ch\!^2\frac{\theta}{2}} + \frac{1+ 2 \theta\, \th\frac{\theta}{2}}{16\,\theta^3\,\ch\!^2\frac{\theta}{2}}\,x^2 + \O(|x|^3) } { \frac{\theta}{2} + \frac{\rt1\!x}{4\,\theta} + \frac{x^2}{16\,\theta^3} + \O(|x|^3) } $$
$$ =\, 1- {\ts\frac{2}{\theta}}\, \th{\ts\frac{\theta}{2}} \times \frac{ 1 + \frac{\rt1\!x}{2\,\theta\,\sh\theta} + \frac{1+ 2 \theta\, \th\frac{\theta}{2}}{8\,\theta^3\,\sh\theta}\,x^2 + \O(|x|^3) } { 1 + \frac{\rt1\!x}{2\,\theta^2} + \frac{x^2}{8\,\theta^4} + \O(|x|^3) } $$
$$ =\, 1- {\ts\frac{2}{\theta}}\, \th{\ts\frac{\theta}{2}} + \frac{\sh\theta - \theta}{2\,\theta^3\,\ch\!^2\frac{\theta}{2}}\,\rt1\!x  + \frac{3(\sh\theta - \theta) - 2\theta^2\, \th\!\frac{\theta}{2} }{8\,\theta^5\,\ch\!^2\frac{\theta}{2}}\,x^2 + \O(|x|^3) \, ; $$ 
$$ \frac{\theta^2+\rt1\!x}{f(\theta^2, x)} \,= \, \frac{\theta^2}{1- {\frac{2}{\theta}}\,\th\!{\ts\frac{\theta}{2}} }\times \frac{1+\frac{\rt1\!x}{\theta^2} }{ 1+ \frac{(\sh\theta - \theta)\,\rt1\!x}{2\,\theta^3\,\ch\!^2\frac{\theta}{2}\,(1- {\frac{2}{\theta}}\,\th\!{\frac{\theta}{2}})}  + \frac{3(\sh\theta - \theta) - 2\theta^2\, \th\!\frac{\theta}{2} }{8\,\theta^5\,\ch\!^2\frac{\theta}{2}\,(1- {\frac{2}{\theta}}\,\th\!{\frac{\theta}{2}})}\,x^2 + \O(|x|^3) } $$
$$ = \, \frac{\theta^2}{1- {\frac{2}{\theta}}\,\th\!{\ts\frac{\theta}{2}} } + \frac{2 +\theta\, \ch\theta - 3\, \sh\theta}{2\,\theta\,\ch\!^2\frac{\theta}{2}\,(1- {\frac{2}{\theta}}\,\th\!{\frac{\theta}{2}})^2}\,\rt1\!x + \frac{h(\theta)\,\,x^2}{8\,\theta^4\,\ch\!^4\frac{\theta}{2}\,(1- {\frac{2}{\theta}}\,\th\!{\frac{\theta}{2}})^3} + \O(|x|^3) \, , $$ 
with \quad ${ h(\theta):=\, \frac{1}{2}\,\theta\,\ch\theta\,\sh\theta - 3\,\sh\!^2\theta + \theta^3\,\sh\theta - \frac{5}{2}\,\theta^2\,\ch\theta + \frac{11}{2}\,\theta\,\sh\theta - \5\,\theta^2\,}$. \parn 
Le coeff de $\rt1\!x\,$ ne s'annule pas\,; mais ce n'est peut-être pas un probl\`eme, cf Prop. \ref{pro.partprincz}.  
Le signe du coeff. $\,h(\theta)$ de $x^2$ ($>0$ ? nettement selon Maple) est plus crucial.  \parn 
En 0 \  ${\ds  h(\theta) = \frac{\theta^6}{24} + \O(\theta^8)}\,$; puis \  $ h'(\theta) =  $
}\fi 

   Then we take \  $ \frac{3}{2} < r < 3\,$, to be specified later, and cut the above integral into\,: 
\begin{equation} \label{f.J01infty} 
 \5\,\II_\e(y,z)\, = \, \Re\Bigg\{\!\int_{0}^{\e^r} \cdots \Bigg\} +\; \Re\Bigg\{\!\int_{\e^r}^{\infty}  \cdots \Bigg\}  \; \equiv \, J_0^\e(y,z) + J_{\infty}^\e(y,z)\, . \vspace{1mm}
\end{equation}
   We first obtain the following (analogue to Proposition \ref{pro.IntegPartPrinc>0}) behaviour of the main contribution $\,J_0^\e(y,z)$ to $\,\II_\e(y,z)$.  \vspace{-1mm} 
\bpro  \label{pro.partprincz}  \  For any $(y,z)\in\R\times\R^*$, set  \ ${\ds C_\e(y,z) := \bigg[ \frac{\pi\,z^2}{2\left(\pi -2\,\th\!{\frac{\pi}{2}}\right)\e^3} + \frac{y-\e}{\e^2} \bigg]}$,  $\,{ K_\e(y,z) := \exp\!\Big(\! -\pi^2 \,C_\e(y,z)  \Big) }$, and \  ${\ds C^2:= \frac{\pi \left(\pi\,\ch\pi - 3\,\sh\pi + 2\pi\right)}{4 \left(\pi\,\ch\frac{\pi}{2} - 2\,\sh\frac{\pi}{2}\right)^2} }\,>0\,$\raise0pt\hbox{.} 
Then as $\e\sea 0\,$ we have 
$$ J_0^\e(y,z) = \, K_\e(y,z) \times \left[\frac{C^2 z^2}{\e^3} + \frac{y-\e}{\e^2}\right]^{-3/4}\! \times \frac{2^{3/4}\,\pi^{11/4}\, \sigma' }{ \sqrt{\pi -2\,\th\!{\frac{\pi}{2}}}\, ( \sh\pi)^{1/4}} \times\, \big[1+ \O\big(\e^{2r-3} + \e^{(3-r)/4}\big)\big] .\vspace{-2mm} $$ 
\epro 
\ub{Proof} \quad In the proof of Lemma \ref{lem.AnalPhiprSing>} we saw that 
$$ \frac{\rt1\! x}{f(0,x)}\, 
=\, {\frac{\pi^3} {\pi -2\,\th\!{\frac{\pi}{2}}}} \times \left[1+ (1+\rt1\!)\, \frac{3\pi -6\,\th\!{\frac{\pi}{2}} -\pi\,\th\!^2{\frac{\pi}{2}}}{\pi \big(\pi -2\,\th\!{\frac{\pi}{2}}\big)}\, u + \O(u^2)\right] $$  
at $\,x= -\rt1\!\pi^2+ 2(1-\rt1\!)\,\pi\,u +u^2 $, \  so that \big(taking $\,h= 2(1-\rt1)\,\pi\,u +u^2\,\Lra u= \frac{1+\rt1}{4\pi}\big(1-\frac{\rt1h}{8\pi^2}+\O(h^2)\big) h\,$\big)\,:
$$ \frac{(\pi^2+\rt1 h)}{f\big(0,-\rt1\!\pi^2+h\big)}\, =\, {\frac{\pi^3} {\pi -2\,\th\!{\frac{\pi}{2}}}} \times \left[1+ \rt1 \frac{3\pi -6\,\th\!{\frac{\pi}{2}} -\pi\,\th\!^2{\frac{\pi}{2}}}{2\pi^2 \big(\pi -2\,\th\!{\frac{\pi}{2}}\big)}\, h + \O(h^2)\right]  , \vspace{-2mm} $$ 
i.e., for small $\,h\,$: \vspace{-4mm} 
$$ \frac{1}{f\big(\pi^2,h\big)}\, 
=\, {\frac{\pi} {\pi -2\,\th\!{\frac{\pi}{2}}}} \times \left[1- \frac{\rt1(\sh\pi - \pi)\, h}{2\pi^2 \big(\pi -2\,\th\!{\frac{\pi}{2}}\big)\, \ch\!^2{\frac{\pi}{2}}} + \O(h^2)\right] . \vspace{-1mm} $$ 
Therefore we have\,: \vspace{-1mm} 
$$ J_0^\e(y,z)\, \equiv  \int_{0}^{\e^r}  \exp\!\left( -\left[\frac{y-\e}{\e^2} + \frac{z^2}{2\e^3 f(\pi^2,x)}\right]\! \big(\rt1\! x + \pi^2\big)\!\right)\!\times \Phi\big(\pi^2\hbox{,}\, x\big)\, dx\, \vspace{-1mm}  $$
$$ =\, \exp\!\left(\! -\pi^2 C_\e(y,z) \right)\! \int_{0}^{\e^r}\! \exp\!\left({\ts - \left[\frac{C^2 z^2}{\e^3} + \frac{y-\e}{\e^2}\right]\! \rt1\! x + \O\big(\e^{2r-3}\big) }\right) \Phi\big(\pi^2\hbox{,}\, x\big)\, dx\,  $$ 
$$ = \, \frac{K_\e(y,z) }{\frac{C^2 z^2}{\e^3} + \frac{y-\e}{\e^2}} \int_{0}^{ \left[\frac{C^2 z^2}{\e^3} + \frac{y-\e}{\e^2}\right] \e^r}\! \exp\!\left[ -\rt1\! x + \O\big(\e^{2r-3}\big) \right]\!\times \Phi\!\left(\pi^2\hbox{,}\, \left[{\ts \frac{C^2 z^2}{\e^3} + \frac{y-\e}{\e^2}} \right]\1 x\right) dx\, . $$
Recall that we took  $\, \frac{3}{2} < r < 3\,$. Now applying Lemma \ref{lem.AnalPhiprSing>} (as for Proposition \ref{pro.IntegPartPrinc>0}) yields\,:  \vspace{-1mm} 
$$ J_0^\e(y,z)\, = \, K_\e(y,z) \times \left[\frac{C^2 z^2}{\e^3} + \frac{y-\e}{\e^2}\right]^{-3/4}\! \frac{2^{3/4}\,\pi^{11/4}}{ \sqrt{\pi -2\,\th\!{\frac{\pi}{2}}}\, ( \sh\pi)^{1/4}} \vspace{-2mm} $$
$$ \qquad \times\, \big[1+ \O\big(\e^{2r-3}+\e^r\big)\big]  \int_{0}^{\e^r \left[\frac{C^2 z^2}{\e^3} + \frac{y-\e}{\e^2}\right]} \sin\!\left[ x + \frac{3\pi}{16}\right]\!\times x^{-1/4}\, dx\, . \vspace{-2mm} $$ 
And finally  \vspace{-1mm} 
$$  \int_{0}^{\e^r \left[\frac{C^2 z^2}{\e^3} - \frac{y-\e}{\e^2}\right]} \sin\!\left[ x + \frac{3\pi}{16}\right]\!\times x^{-1/4}\, dx\,  = \,\sigma'\, + \O\big(\e^{(3-r)/4}\big) = \big[1+\O\big(\e^{(3-r)/4}\big)\big]\,\sigma' . \;\;\diamond $$ 
\parm 

  To estimate the contribution $\,J_\infty^\e(y,z)$, we crucially need the following. 
\blem \label{lem.PositivRez} \  We have \quad 
${\ds \Re\left[\frac{\pi^2+\rt1\! x}{f(\pi^2,x)}\right] > \frac{\pi^2}{f(\pi^2,0)}\, \raise1pt\hbox{,} }$ \  for any positive $\,x\,$. \vspace{-0mm} 
\elem
\ub{Proof} \quad From (\ref{df.fnf}) we compute\,:  \  for any real $\,\chi > -\pi^2$ (we actually need $\,\chi =\pi^2$), 
$$ \Re\left[ \frac{\chi + \rt1\! x}{f(\chi,x)}\right] = \, \frac{ (a^4-b^4)(\ch a + \cos b) -2(a^2-3b^2)\,a\,\sh a - 2 (3a^2-b^2)\,b\, \sin b } {(a^2+b^2)(\ch a + \cos b) - 4(a\, \sh a + b\, \sin b) + 4(\ch a - \cos b) } $$
$$ = \, \frac{\chi (2a^2-\chi)(\ch a + \cos b) + 2(2 a^2-3\chi)\,a\,\sh a - 2 (2a^2+\chi)\,b\, \sin b } {(2a^2-\chi)(\ch a + \cos b) - 4(a\, \sh a + b\, \sin b) + 4(\ch a - \cos b) }\, \raise1.9pt\hbox{,} $$ 
and then, as \ ${\ds f(\pi^2,0)= 1-\frac{2\,\sh\pi}{\pi\,(\ch\pi +1)} = 1 - {\ts \frac{2}{\pi}}\,\th\!{\ts \frac{\pi}{2}} =: 1-\lambda }\,$ and $\,a^2-b^2 = \pi^2 $, we have\,: \vspace{-1mm} 
$$ \Re\left[\frac{\pi^2+\rt1\! x}{f(\pi^2,x)}\right] -\frac{\pi^2}{f(\pi^2,0)}\, =\, \frac{N}{D}\, $$ 
with \  $D = \big[(a^2+b^2)(\ch a + \cos b) - 4(a\, \sh a + b\, \sin b) + 4(\ch a - \cos b) \big]\! \big[1 - {\ts \frac{2}{\pi}}\,\th\!{\ts \frac{\pi}{2}}\big] >0\,$  \  and 
$$ \5\,N =\,\5 \big[\pi^2 (a^2+b^2)(\ch a + \cos b) + 2(2 a^2-3\pi^2)\,a\,\sh a - 2 (2a^2+\pi^2)\,b\, \sin b\big] \! \big[1 - {\ts \frac{2}{\pi}}\,\th\!{\ts \frac{\pi}{2}}\big] \vspace{-1mm} $$
$$ - \,\5 \big[(a^2+b^2)(\ch a + \cos b) - 4(a\, \sh a + b\, \sin b) + 4(\ch a - \cos b) \big] \pi^2 \vspace{1mm}  $$
$$ =\, \big[ 2(1-\lambda)\, a^2 + (3\lambda -1)\,\pi^2\big]\, a\,\sh a - \big[ 2(1-\lambda)\, a^2 - (\lambda +1)\,\pi^2\big]\,b\, \sin b $$
$$ -\, \lambda\,\pi^2 (a^2-\pi^2/2)(\ch a + \cos b) - 2\pi^2\,(\ch a - \cos b)\, $$ 
$$ =\, \big[ 2(1-\lambda)\, b^2 + (\lambda +1)\,\pi^2\big]\, \!\sqrt{\pi^2+b^2}\,\sh\!\sqrt{\pi^2+b^2} - \big[ 2(1-\lambda)\, b^2 - (3\lambda -1)\,\pi^2\big]\,b \sin b \vspace{-1mm} $$
$$ -\, \lambda\,\pi^2 (b^2+\pi^2/2)\big(\ch\!\sqrt{\pi^2+b^2} + \cos b\big) - 2\pi^2 \big(\ch\!\sqrt{\pi^2+b^2} - \cos b\big) $$ 
$$ = \big[ (\lambda +1)\,\pi^2+2(1-\lambda)\, b^2\big] \pi \sqrt{1+{\ts \frac{b^2}{\pi^2}}}\,\sh\!\!\left[\pi \sqrt{1+{\ts \frac{b^2}{\pi^2}}}\,\right] + \big[ (3\lambda -1)\,\pi^2 - 2(1-\lambda)\, b^2\big]\,b \sin b \vspace{-2mm}  $$
$$ - \big(\lambda\,{\ts\frac{\pi^2}{2}} + 2 + \lambda\,b^2 \big) \pi^2\, \ch\!\!\left[\pi \sqrt{1+{\ts \frac{b^2}{\pi^2}}}\,\right] - \big(\lambda\,{\ts\frac{\pi^2}{2}} - 2 + \lambda\,b^2 \big) \pi^2 \cos b $$ 
$$ = \Big[\big({\ts\frac{3}{2}}(1-\lambda) - \lambda {\ts\frac{\pi^2}{4}}\big)\pi\,\sh\pi + {\ts\frac{1-\lambda}{2}}\,\pi^2 \ch\pi - 2 (1-\lambda)\pi^2+ \lambda{\ts\frac{\pi^4}{4}}\Big] b^2+\O(b^4) \approx\, 
0.186\, b^2 +\O(b^4) . $$
$(i)$ \  Positivity of $N\,$ for $\,0< b\le 1.16\,$: in this range we successively have\,:  \parn 
$$ \frac{N}{2} > \big[ (\lambda +1)\,\pi^2+2(1-\lambda)\, b^2\big] \pi \big[1+{\ts \frac{b^2}{2\pi^2} -\frac{b^4}{8\pi^4}}\big] \sh\!\!\left[\pi \big[1+{\ts \frac{b^2}{2\pi^2} -\frac{b^4}{8\pi^4}}\big]\right] \vspace{-1mm}  $$
$$ + \big[ (3\lambda -1)\,\pi^2 - 2(1-\lambda)\, b^2\big]b \big(b-{\ts \frac{b^3}{6} + \frac{b^5}{125}}\big) \vspace{-1mm}  $$
$$ -\, \big(\lambda\,{\ts\frac{\pi^2}{2}} + 2 + \lambda\,b^2 \big) \pi^2\, \ch\!\!\left[\pi \big[1+{\ts \frac{b^2}{2\pi^2} -\frac{b^4}{8\pi^4} + \frac{b^6}{16\pi^6}}\big]\right] - \big(\lambda\,{\ts\frac{\pi^2}{2}} - 2 + \lambda\,b^2 \big) \pi^2 \big(1-{\ts \frac{b^2}{2} + \frac{b^4}{24} - \frac{b^6}{740}}\big)\, ; $$ 
$$ \sh\!\!\left[\pi +{\ts \frac{b^2}{2\pi} -\frac{b^4}{8\pi^3}}\right] \ge\, \sh\pi + {\ts \frac{b^2\ch\pi}{2\,\pi}} + {\ts \frac{b^4(\pi\,\sh\pi-\ch\pi)}{8\,\pi^3}} \, $$ 
\Big(indeed, if $\,f_0(b):=\sh\!\!\left[\pi +{\ts \frac{b^2}{2\pi} -\frac{b^4}{8\pi^3}}\right] - \sh\pi - {\ts \frac{b^2\ch\pi}{2\,\pi}} - {\ts \frac{b^4(\pi\,\sh\pi-\ch\pi)}{8\,\pi^3}}$\raise1pt\hbox{,} $\;f_1(b):= \frac{f_0'(\pi b)}{b-b^3/2}\,$, \ then \parn 
\centerline{$\,f_2(b):= \frac{(1-{b^2}/{2})^2\,f_1'(b)}{\pi b}\, =\, \sh\!\!\left[\pi +{\ts \frac{b^2}{2\pi} -\frac{b^4}{8\pi^3}}\right]\!(1-{\ts\frac{b^2}{2}})^3 -\sh\pi\,$,} \parn
$\,\frac{f'_2(b)}{b(1-b^2/2)^2} = \pi\, \ch\!\!\left[\pi +{\ts \frac{b^2}{2\pi} -\frac{b^4}{8\pi^3}}\right]\!(1-{\ts\frac{b^2}{2}})^2 - 3\,\sh\!\!\left[\pi +{\ts \frac{b^2}{2\pi} -\frac{b^4}{8\pi^3}}\right]$ decreases on $[0,\sqrt{2}\,]$ and is positive near 0, whence $\,f_2$ increases near 0 and then possibly decreases, whence the same for $f_1$ and $f_0\,$, and $\,f_0(1.2)>0$, whence finally $\,f_0>0$ on $]0, 1.2]$\Big)\! ; 
$$ \ch\!\!\left[\pi +{\ts \frac{b^2}{2\pi} -\frac{b^4}{8\pi^3}+ \frac{b^6}{16\pi^5}}\right] \le \,\ch\pi + {\ts \frac{b^2\sh\pi}{2\,\pi}} + {\ts \frac{b^4(\pi\,\ch\pi-\sh\pi)}{8\,\pi^3}} + {\ts \frac{b^6\,\ch\pi}{36\,\pi^4}} \, $$ 
\Big(indeed, if $\,g_0(b):=$ the left hand side minus the right hand side, $\,g_1(b):= \frac{g_0'(\pi b)}{b-b^3/2}\,$, \  then \parn 
$\,g_2(b):=\frac{(1-{b^2}/{2})^2\,g_1'(b)}{\pi b} = \ch\!\!\left[\pi +{\ts \frac{b^2}{2\pi} -\frac{b^4}{8\pi^3}}\right]\!(1-{\ts\frac{b^2}{2}})^3 -\frac{(6+4b^2-b^4)\ch\pi}{6}\,$, $\,g_3(b):= \frac{g_2'(b)}{b(1-b^2/2)^2} = \pi\, \sh\!\!\left[\pi +{\ts \frac{b^2}{2\pi} -\frac{b^4}{8\pi^3}+ \frac{b^6}{16\pi^5}}\right]\! (1-{\ts\frac{b^2}{2}})^2 - 3\,\ch\!\!\left[\pi +{\ts \frac{b^2}{2\pi} -\frac{b^4}{8\pi^3}+ \frac{b^6}{16\pi^5}}\right] - \frac{4\,\ch\pi }{3(1-b^2/2)} $ decreases on $[0,\sqrt{2}\,]$, which easily entails $\,g_2(b)<0$, whence then easily $\,g_0(b)<0$ on $]0, 1.2]$\Big). \parn 
Whence \vspace{-2mm}  
$$  {N}/{2}\, > \big[ (\lambda +1)\,\pi^2+2(1-\lambda)\, b^2\big]\! \big[\pi+{\ts \frac{b^2}{2\pi} -\frac{b^4}{8\pi^3}}\big] \!\left[\sh\pi + {\ts \frac{b^2\ch\pi}{2\,\pi}} + {\ts \frac{b^4(\pi\,\sh\pi-\ch\pi)}{8\,\pi^3}}\right] \vspace{-1mm}  $$
$$ + \big[ (3\lambda -1)\,\pi^2 - 2(1-\lambda)\, b^2\big]b \big(b-{\ts \frac{b^3}{6} + \frac{b^5}{125}}\big) - \big(\lambda\,{\ts\frac{\pi^2}{2}} - 2 + \lambda\,b^2 \big) \pi^2 \big(1-{\ts \frac{b^2}{2} + \frac{b^4}{24} - \frac{b^6}{740}}\big)\vspace{-1mm}  $$
$$ - \big(\lambda\,{\ts\frac{\pi^2}{2}} + 2 + \lambda\,b^2 \big) \pi^2\!\left[\ch\pi + {\ts \frac{b^2\sh\pi}{2\,\pi}} + {\ts \frac{b^4(\pi\,\ch\pi-\sh\pi)}{8\,\pi^3}} + {\ts \frac{b^6\,\ch\pi}{36\,\pi^4}}\right] \, = \, b^2\, P(b^2)\, , $$ 
where $P$ is polynomial of degree 4, which has the same behaviour as $\frac{N}{2b^2}$ near 0, increases first  and then decreases on $[0,1.16]$, with a positive value at 1.16. Therefore it is positive on $]0,1.16]$, and then $N$ as well. \parsn 
$(ii)$ \  Positivity of $N\,$ for $\,1.16\le b\le5\,$: \quad on the one hand, in this range we have  
$$\big[ (3\lambda -1)\,\pi^2 - 2(1-\lambda)\, b^2\big] b \sin b - \big(\lambda\,{\ts\frac{\pi^2}{2}} - 2 + \lambda\,b^2 \big) \pi^2 \cos b\, > 0\, $$ 
\Big(indeed, since $\,\frac{\pi}{3}<1.16<\frac{\pi}{2}< \sqrt{\frac{3\lambda -1}{2-2\lambda}}\,\pi<\pi<\frac{3\pi}{2}<5<2\pi\,$, denoting the left hand side by $f(b)$ we see at once that $\,f$ is positive on $\big[\frac{\pi}{2} ,\sqrt{\frac{3\lambda -1}{2-2\lambda}}\,\pi\big] \cup \big[\pi,\frac{3\pi}{2}\big]$, \   and that   \parn
$\,f'(b) = \big[ (3\lambda -3+ \lambda\,{\ts\frac{\pi^2}{2}})\pi^2 +(\lambda\,\pi^2- 6+6\lambda)\,b^2\big] \sin b - (1-\lambda)(\pi^2 +2\,b^2)\,b \cos b\;$ is positive on $\big[\frac{\pi}{2} ,\pi\big]$ and negative on  $\big[\frac{3\pi}{2},5\big]$. Since $\,f(5) \approx 20>0\,$, this shows the wanted positivity on $\big[\frac{\pi}{2},5\big]$.  Moreover for $\,1.16\le b\le \frac{\pi}{2}\,$ we have $\,\cos b < \sin b/\sqrt{3}\,$ and then $\,f'(b)> Q(b)\,\sin b\,$, with $\,Q(b):= (3\lambda -3+ \lambda\,{\ts\frac{\pi^2}{2}})\pi^2 +(\lambda\,\pi^2- 6+6\lambda)\,b^2 - (\frac{1-\lambda}{\sqrt{3}})(\pi^2 +2\,b^2)\,b \,$, which is easily seen to satisfy $Q''(b)>2$, and then to be increasing and positive\Big). \   Whence on the other hand\,:  \vspace{-1mm}
$$ \frac{N}{2\, \ch\!\sqrt{\pi^2+b^2}}\, \ge \big[ (\lambda +1)\,\pi^2+2(1-\lambda)\, b^2\big] \sqrt{\pi^2+b^2}\,\th\!\sqrt{\pi^2+b^2} - \big(\lambda\,{\ts\frac{\pi^2}{2}} + 2 + \lambda\,b^2 \big) \pi^2 \vspace{-1mm} $$
$$ \ge \big[ (\lambda +1)\,\pi^2+2(1-\lambda)\, b^2\big] {\ts\frac{\pi + b}{\sqrt{2}}}\,\th\pi - \big(\lambda\,{\ts\frac{\pi^2}{2}} + 2 + \lambda\,b^2 \big) \pi^2 $$
$$ \ge \big[ (\lambda +1)\,\pi^2+2(1-\lambda)(1.16)^2\big] {\ts\frac{\pi + 1.16}{\sqrt{2}}}\,\th\pi - \big(\lambda\,{\ts\frac{\pi^2}{2}} + 2 + \lambda(1.16)^2 \big) \pi^2 > 0 \, . \vspace{1mm} $$
$(iii)$ \  Positivity of $N\,$ for $\, b\ge 5\,$: \  estimating $\,|\cos b\,|\,$ and $\,|\sin b\,|\,$ by 1, we get \vspace{-1mm}
$$ \frac{N}{2\, \ch\!\sqrt{\pi^2+b^2}}\, \ge\, \big[ (\lambda +1)\,\pi^2+2(1-\lambda)\, b^2\big] \frac{\pi + b}{\sqrt{2}}\,\th\!\sqrt{\pi^2+25} \vspace{-2mm}  $$
$$ - \big(\lambda\,{\ts\frac{\pi^2}{2}} + 2 + \lambda\,b^2 \big) \pi^2 - \frac{\big(\lambda\,{\ts\frac{\pi^2}{2}} - 2 + \lambda\,b^2 \big) \pi^2 + 2(1-\lambda)\, b^3}{\ch\![(\pi+b)/\sqrt{2}\,]} $$ 
$$ \ge\, {\ts\frac{\th\!\sqrt{\pi^2+25}}{\sqrt{2}}}\big[ (\lambda +1)\,\pi^2+2(1-\lambda)\, b^2\big] (\pi + b) - \lambda\,\pi^2b^2 - 50 $$
$$ \ge \, {\ts\frac{\th\!\sqrt{\pi^2+25}}{\sqrt{2}}}\big[ (\lambda +1)\,\pi^2+50(1-\lambda)\big] (\pi + 5) - 25\lambda\,\pi^2 - 50\, > 0\, . \;\;\diamond  $$ 
\parm 

  To handle the remaining contribution $J_\infty^\e(y,z)$ in (\ref{f.J01infty}) we still need the following. \vspace{-1mm} 
\blem \label{lem.lemzeroz} \  The equation \ $(a+\rt1\! b)\,\ch\!(a+\rt1\! b) - 3\, \sh\!(a+\rt1\! b) + 2 (a+\rt1\! b) = 0\, $\parn  has no solution such that $\,b\ge 0,\; a=\sqrt{\pi^2+b^2}\,$. 
\elem
\ub{Proof} \quad This equation is equivalent to \vspace{-2mm} 
$$ (a\,\ch a -3\,\sh a) \cos b - b\,\sh a \sin b + 2a = 0 = b\, \ch a \cos b + (a\,\sh a -3\,\ch a) \sin b + 2b \,, \vspace{-1mm} $$ 
i.e., to \vspace{-1mm} 
{\small 
$$ \frac{\cos\!\sqrt{a^2-\pi^2}}{2} = \frac{3a\,\ch a - (2a^2-\pi^2)\,\sh a}{(a^2+\frac{9-\pi^2}{2})\,\sh\!(2a) - 3a\,\ch\!(2a)}\; ; \  \frac{\sin\!\sqrt{a^2-\pi^2}}{2} = \frac{3\sqrt{a^2-\pi^2}\,\sh a}{(a^2+\frac{9-\pi^2}{2})\,\sh\!(2a) - 3a\,\ch\!(2a)}\, \raise1.9pt\hbox{.} $$ } 
Now for $\,a\ge \pi\,$: \quad  $a^2\,\frac{d}{da}\big(3\,\ch a - (2a-\frac{\pi^2}{a})\,\sh a\big) = (a^2-\pi^2)\,\sh a -(2a^2-\pi^2)\,a\,\ch a \,$ \parn
$ < -(2a^3-a^2-\pi^2 a+\pi^2)\,\sh a \le -\pi^3\sh a < 0\,$, whence \ $3a\,\ch a - (2a^2-\pi^2)\,\sh a \le -3\pi/2 <0$. \parn
And since \  $(a^2+\frac{9-\pi^2}{2})\,\sh\!(2a) - 3a\,\ch\!(2a) = \frac{(a^2+3a+\frac{9-\pi^2}{2})\, e^{2a}}{2}\Big[\frac{a^2-3a+\frac{9-\pi^2}{2}}{a^2+3a+\frac{9-\pi^2}{2}} - e^{-4a} \Big] $ \parn 
$\ge \frac{(a^2+3a+\frac{9-\pi^2}{2}) e^{2a}}{2}\Big[\frac{\pi^2-3\pi+\frac{9-\pi^2}{2}}{\pi^2+3\pi+\frac{9-\pi^2}{2}} - e^{-4\pi} \Big] >0\,$,  \  we must have $\,\cos\!\sqrt{a^2-\pi^2}\,<0\,$, hence $\,a>\pi\frac{\sqrt{5}}{2}\,$\raise1.5pt\hbox{.} \parn 
This implies \vspace{-1mm} 
{\small 
$$ \frac{\sin\!\sqrt{a^2-\pi^2}}{2} <\, \frac{3\sqrt{a^2-\pi^2}\, e^{-a}}{a^2+3a+\frac{9-\pi^2}{2}}\! \left[\frac{\frac{5\pi^2}{4}-\frac{3\pi\sqrt{5}}{2}+\frac{9-\pi^2}{2}}{\frac{5\pi^2}{4}+\frac{3\pi\sqrt{5}}{2}+\frac{9-\pi^2}{2}} - e^{-2\pi\sqrt{5}}\right]\1\! < 16.44\,e^{-\pi\frac{\sqrt{5}}{2}}\frac{3\sqrt{a^2-\pi^2}\,}{a^2+3a+\frac{9-\pi^2}{2}} \, \raise1.9pt\hbox{,}  $$ }
whence \quad $\sin\!\sqrt{a^2-\pi^2} < \frac{3\sqrt{a^2-\pi^2}\,}{a^2+3a+\frac{9-\pi^2}{2}}\,$\raise1.9pt\hbox{.} \quad  Now the polynomial \parn
$\tilde P (a) := a^4+6a^3+(18-4\pi^2-c^2)a^2+3(9-\pi^2)a + (\frac{9-\pi^2}{2})^2 + \pi^2c^2\,$ is positive for $\,a> \pi\frac{\sqrt{5}}{2}\,$ and $c\le 5.8$ (actually, $\tilde P''(a)>0$ and $\tilde P'(a)>0$ in that range), \  so that we have $\,5.8\times \sqrt{a^2-\pi^2} < a^2+3a+\frac{9-\pi^2}{2}\,$, hence \   $\sin\!\sqrt{a^2-\pi^2} < \frac{3\sqrt{a^2-\pi^2}\,}{a^2+3a+\frac{9-\pi^2}{2}} <  \frac{3}{5.8}\,$\raise1.2pt\hbox{,} hence  $\, \sqrt{a^2-\pi^2} > \pi- \arcsin(\frac{3}{5.8})\,$,  and then $\, a> 4$. \  Finally $\, a> 4$ implies  \parn 
{
$$ {|\cos\!\sqrt{a^2-\pi^2}\,|} <\, \frac{2(2a^2-3a-\pi^2)}{a^2+3a+\frac{9-\pi^2}{2}}\, e^{-4}\! \left[\frac{4+\frac{9-\pi^2}{2}}{28+\frac{9-\pi^2}{2}} - e^{-16}\right]\1\! < 4\times 0.142 < 0.57\,, \vspace{-2mm}  $$ }
\hskip-1mm and similarly\,: \vspace{-2mm} 
$$ {|\sin\!\sqrt{a^2-\pi^2}\,|}  <\, \frac{3\sqrt{a^2-\pi^2}\,\, e^{-4}}{a^2+3a+\frac{9-\pi^2}{2}} \left[\frac{4+\frac{9-\pi^2}{2}}{28+\frac{9-\pi^2}{2}} - e^{-16}\right]\1\! <\,\frac{3}{5.8}\times 0.142 < 0.08\, \raise0pt\hbox{,} \vspace{-2mm} $$ 
which forbids any solution $\,a > 4$. $\;\diamond$ \parm 

   We finally have the following control of the remaining contribution $J_\infty^\e(y,z)$ in (\ref{f.J01infty}). \vspace{-1mm} 
\blem \label{lem.ControlIntPhifz} \  For any $(y,z)\in\R\times\R^*\,$ we have \  ${\ds J_\infty^\e(y,z) = J_0^\e(y,z)\times \O\big(\e^{\frac{3-r}{4}}\big) . \vspace{-0mm}  }$ 
\elem
\ub{Proof} \quad In the spirit of Lemma \ref{lem.ControlIntPhif>}, we integrate by parts. For that, consider the derivative of the exponent\,:  \vspace{-1mm} 
$$ F_\e(x) :=  \frac{z^2}{2\e^3}\times\! \left( {\frac{\rt1}{f(\pi^2,0)}-\frac{\rt1}{f(\pi^2,x)}} + { \frac{(\rt1\! x +\pi^2)\,\frac{df}{dx}(\pi^2,x)}{f(\pi^2,x)^2} } \right) - \rt1 C_\e(y,z) \, , $$ 
so that  
for any large $A$ we have\,: \vspace{-2mm} 
$$  \Re \Bigg\{\!\int_{\e^r}^{A} \exp\!\left( {  \frac{z^2}{2\,\e^3}\! \left[\frac{1}{f(\pi^2,0)}-\frac{1}{f(\pi^2,x)}\right]}\!\big(\rt1\! x +\pi^2\big) -C_\e(y,z) \rt1\! x \!\right)\!\times \Phi\big(\pi^2\hbox{,}\, x\big)\, dx \Bigg\}  $$
$$ = \,\exp\!\left( {  \frac{-z^2}{2\,\e^3}\! \left[ \Re\!\left(\frac{\pi^2+\rt1 A}{f(\pi^2,A)}\right) -\frac{\pi^2}{f(\pi^2,0)} \right] }\right)\!\times \O\left[\bigg|\frac{\Phi\big(\pi^2\hbox{,}\, A\big)}{F_\e(A)}\bigg|\right] $$ 
$$ + \,\exp\!\left( { \frac{-z^2}{2\,\e^3}\! \left[ \Re\!\left(\frac{\pi^2+\rt1 \e^r}{f(\pi^2,\e^r)}\right) -\frac{\pi^2}{f(\pi^2,0)} \right] }\right)\!\times \O\left[\bigg|\frac{\Phi\big(\pi^2\hbox{,}\, \e^r\big)}{F_\e(\e^r)}\bigg|\right] $$ 
$$ -\, \Re \Bigg\{\!\int_{\e^r}^{A}\! \exp\!\left( {\ts  \frac{z^2}{2\,\e^3}\! \left[\frac{1}{f(\pi^2,0)}-\frac{1}{f(\pi^2,x)}\right]}\!\big(\rt1\! x +\pi^2\big) -C_\e(y,z) \rt1\! x \!\right)\!\times \frac{d}{dx}\bigg[\frac{\Phi\big(\pi^2\hbox{,}\, x\big)}{F_\e(x)}\bigg] dx \Bigg\} $$
\begin{equation}  \label{f.controlJ1}
=\, \frac{\O\big(A\,e^{-\sqrt{A/8}}\,\big)}{\big|F_\e(A)\big|} + \frac{\O(\e^{-r/4})}{\big|F_\e(\e^r)\big|} + \int_{\e^r}^{A} \O\left[\bigg|\frac{d}{dx}\bigg[\frac{\Phi\big(\pi^2\hbox{,}\, x\big)}{F_\e(x)}\bigg]\bigg|\right] dx\, ,
\end{equation}
by Lemmas \ref{lem.PositivRez} and \ref{lem.ControlPhi'>}.  \quad 
We have to estimate $\,\big|F_\e(x)\big|\1$, for $\,x\ge \e^r$. \  Now, writing $\,f$ for $f(\pi^2,\cdot)$, we have  \vspace{-1mm} 
$$ F_\e(x)\1 \,=\, \frac{\O(\e^3)\, f(x)}{\frac{(\pi^2+\rt1\! x) f'(x)}{f(x)} +\rt1\! \left[\frac{f(x)}{f(0)} -1\right] - 2\rt1\! z\2\, \e^3\, C_\e(y,z) f(x)} \,\qquad \vspace{-1mm} $$ 
$$ \quad =\, \frac{\O(\e^3)\, f(x)}{\frac{(\pi^2+\rt1\! x) f'(x)}{\rt1 f(x)} - 1 - 2 z\2 (y-\e)\,\e \, f(x)} \ \hbox{\big(where $\O(\e^3)$ does not depend on $x$\big)}\raise0pt\hbox{.} $$ 
We shall verify below that $\,F_\e$ cannot vanish on $\R_+^*$. \par  

   In the proof of Proposition \ref{pro.partprincz} we saw that for small positive $x$ we have \vspace{-1mm} 
$$ \frac{1}{f(x)}\, 
=\, {\frac{\pi} {\pi -2\,\th\!{\frac{\pi}{2}}}} \times \left[1- \frac{\rt1(\sh\pi - \pi)\, x}{2\pi^2 \big(\pi -2\,\th\!{\frac{\pi}{2}}\big)\, \ch\!^2{\frac{\pi}{2}}} + \O(x^2)\right] . \vspace{-1mm} $$ 
Thus in particular we have \  $\frac{f'(0)}{f(0)} = \frac{\rt1(\sh\pi - \pi)}{2\pi^2 \big(\pi -2\,\th\!{\frac{\pi}{2}}\big)\, \ch\!^2{\frac{\pi}{2}}}\,$ and then \parsn 
$\frac{\pi^2f'(0)}{f(0)} - \frac{2\rt1\! \e^3}{z^2}\, C_\e(y,z) f(0) = \frac{\rt1(\sh\pi - \pi)}{2 (\pi -2\,\th\!{\frac{\pi}{2}}) \ch\!^2{\frac{\pi}{2}}} - \rt1\! + \O(\e) = \frac{\rt1\!(3\,\sh\pi - 2\pi - \pi\,\ch\pi)}{(1+\ch\pi)(\pi -2\,\th\!{\frac{\pi}{2}})} + \O(\e) \not=0. $ \parsn
Hence \  $\,F_\e(0)\not=0\,$ and $\,{\ds \frac{\O(\e^{-r/4})}{|F_\e(\e^r)|} = \O\big(\e^{3-r/4}\big)}$. \quad 
Then by (\ref{df.fnf}) we have \  $f(x) = 1+ \O\big(x^{-1/2}\big)$ \parn
for large $\,x$, which entails that $\,f\,$ is bounded on $\R_+$,  \  and also\,:  
$$ f'(x)\, = \, \frac{\rt1\! \big(\sh\!\sqrt{\pi^2+\rt1\! x} - \sqrt{\pi^2+\rt1\! x}\,\big)}{(\pi^2+\rt1\! x)^{3/2} \big(1+\ch\!\sqrt{\pi^2+\rt1\! x}\,\big)}\, =\, \frac{\big[b(3a^2+b^2)+\rt1\! a(a^2-3b^2)\big] }{(a^2+b^2)^3\times (\ch a+\cos b)^2 }\, \times \quad $$ 
$$ \hskip 9mm \times\big[ (\ch a+\cos b)(\sh a +\rt1\!\sin b) - (a+\rt1\! b)(1+\ch a\cos b - \rt1\! \sh a \sin b)\big] , $$ 
which shows that for large $x\,$ we have \  $f'(x) = \O\big((a^2+b^2)^{-3/2}\big) = \O\big(x^{-3/2}\big)$. \parn
Hence on the one hand  for large $\,x\,$ we have  \vspace{-1mm} 
$$ \frac{\O\big(A\,e^{-\sqrt{A/8}}\,\big)}{\big|F_\e(A)\big|}  = \frac{\O\big(\e^3 A\,e^{-\sqrt{A/8}}\,\big)}{1+\O(A^{-1/2} +\e)}\,\raise1pt\hbox{,} \quad   \hbox{  and } \quad    \frac{(\pi^2+\rt1\! x) f'(x)}{\rt1  f(x)} = \O\big(x^{-1/2}\big)\, . $$
This shows that for any small enough $\e$  $\,F_\e$ cannot vanish for large $\,x$. \  On the other hand, \vspace{-0mm} 
$$ \frac{(\pi^2+\rt1\! x) f'(x)}{\rt1 f(x)}  = \frac{ \sh\!\sqrt{\pi^2+\rt1\! x} - \sqrt{\pi^2+\rt1\! x} }{\big(\sqrt{\pi^2+\rt1\! x} - 2\,\th\!\big[\sqrt{\pi^2+\rt1\! x}/2\big]\big) \big(1+\ch\!\sqrt{\pi^2+\rt1\! x}\,\big)} \vspace{-1mm}  $$ 
$$ = \frac{ \sh\!(a+\rt1\! b) - (a+\rt1\! b) }{(a+\rt1\! b)\big(1+\ch\!(a+\rt1\! b)\big) - 2\,\sh\!(a+\rt1\! b) } $$ 
is equal to 1 \iff \vspace{-1mm} 
$$ (a+\rt1\! b)\,\ch\!(a+\rt1\! b) - 3\, \sh\!(a+\rt1\! b) + 2 (a+\rt1\! b) = 0\, , \vspace{-1mm} $$ 
which is forbidden by Lemma \ref{lem.lemzeroz}. This shows that $\,F_\e$ cannot vanish on $\R_+^*$, and moreover, \parn
that for some small positive $\e_0$ 
we have $\,\inf_{0<\e\le \e_0, x\ge 0}\limits \Big|\frac{(\pi^2+\rt1\! x) f'(x)}{\rt1 f(x)} - 1 - 2 z\2 (y-\e)\,\e \, f(x)\Big| >0\,$.  \pars  

   So far, according to (\ref{f.J01infty}) and (\ref{f.controlJ1}) (in which we let $A\to\infty$) and the above, we have\,: 
$$ \frac{J_\infty^\e(y,z)}{K_\e(y,z)}\, =\,  \O\big(\e^{3-\frac{r}{4}}\big) + \int_{\e^r}^{\infty} \bigg|\frac{\O(\e^3)\, f(x)\, \frac{d\Phi}{dx}(\pi^2\hbox{,}\, x)}{\frac{(\pi^2+\rt1\! x) f'(x)}{f(x)} -\rt1\! + \O(\e)}\bigg|\, dx\, \vspace{-2mm} \qquad $$ 
$$ \qquad +\, \O(\e^3) \int_{\e^r}^{\infty} \Bigg|\Phi\big(\pi^2\hbox{,}\, x\big)\! \times \frac{d}{dx}\Bigg[\frac{f(x)}{\frac{(\pi^2+\rt1\! x) f'(x)}{\rt1 f(x)} - 1 - 2 z\2 (y-\e)\,\e \, f(x)}\Bigg] \Bigg|\, dx\,  $$ 
$$ =\,\O\big(\e^{3-\frac{r}{4}}\big) + \O(\e^3) \int_{\e^r}^{\infty} \bigg|\frac{d\Phi}{dx}(\pi^2\hbox{,}\, x)\bigg|\, dx\, + \O(\e^3) \int_{\e^r}^{\infty} \Big|\Phi\big(\pi^2\hbox{,}\, x\big) f'(x)\Big|\, dx \vspace{-0mm} \qquad $$ 
$$ \qquad +\, \O(\e^3) \int_{\e^r}^{\infty} \bigg|\Phi\big(\pi^2\hbox{,}\, x\big)\! \times \frac{d}{dx}\bigg[\frac{(\pi^2+\rt1\! x) f'(x)}{\rt1 f(x)} - 1 - 2 z\2 (y-\e)\,\e \, f(x)\bigg] \bigg|\, dx\, $$ 
$$ =\,\O\big(\e^{3-\frac{r}{4}}\big) + \O(\e^3) \int_{\e^r}^{\infty} \bigg|\frac{d\Phi}{dx}(\pi^2\hbox{,}\, x)\bigg|\, dx\, + \O(\e^3) \int_{\e^r}^{\infty} \Big|\Phi\big(\pi^2\hbox{,}\, x\big) f'(x)\Big|\, dx \vspace{-0mm} \qquad $$ 
$$ \qquad +\, \O(\e^3) \int_{\e^r}^{\infty} \big|\Phi\big(\pi^2\hbox{,}\, x\big)\big| \times \Big[\big|(\pi^2+\rt1\! x) f''(x)\big|+\big|(\pi^2+\rt1\! x)f'(x)^2\big|\Big] dx\, . $$ 
Using Lemma \ref {lem.ControlPhi'>} this yields 
$$ \frac{J_\infty^\e(y,z)}{K_\e(y,z)}\, =\, \O\big(\e^{3-\frac{r}{4}}\big) + \O(\e^3)\! \int_{\e^r}^{1} \Big[ 1+ \big|f'(x)\big| + \big|f''(x)\big| + \big|f'(x)\big|^2\Big]\, x^{-1/4}\, dx \vspace{-1mm} $$ 
$$ \qquad +\, \O(\e^3) \int_{1}^{\infty} \Big[ 1+ x \big|f''(x)\big|+ x \big|f'(x)\big|^2 \Big]\, x\,e^{-\sqrt{x/8}}\, dx\, . \vspace{1mm} $$ 
We already saw above that \  $f'(x) = \O\big(x^{-3/2}\big)$ for large $x\,$. Lemma \ref{lem.racinf} shows that $\,f'$ and $f''$ are continuous  on $\R_+$. It remains to control $\,f''$ on $[1,\infty[\,$. \  Now from the expression of $f'$ displayed above we have 
{\small 
$$ f''(x)\, = \, \frac{3 \big(\sh\!\sqrt{\pi^2+\rt1\! x} - \sqrt{\pi^2+\rt1\! x}\,\big)}{4(\pi^2+\rt1\! x)^{5/2} \big(1+\ch\!\sqrt{\pi^2+\rt1\! x}\,\big)} - \frac{\ch\!\sqrt{\pi^2+\rt1\! x} - 1}{2(\pi^2+\rt1\! x)^{2} \big(1+\ch\!\sqrt{\pi^2+\rt1\! x}\,\big)}$$ }
$$ +\, \frac{\big(\sh\!\sqrt{\pi^2+\rt1\! x} - \sqrt{\pi^2+\rt1\! x}\,\big)\, \sh\!\sqrt{\pi^2+\rt1\! x}}{2(\pi^2+\rt1\! x)^{2}\, \big(1+\ch\!\sqrt{\pi^2+\rt1\! x}\,\big)^2} \, =\, \O\big( x\2 \big) . $$ 
Therefore we finally obtain\,: 
$$ \frac{J_\infty^\e(y,z)}{K_\e(y,z)} = \O\big(\e^{3-\frac{r}{4}}\big) + \O(\e^3)\!\left[ \int_{\e^r}^{1} x^{-1/4}\, dx + \int_{1}^{\infty}\! e^{-\sqrt{x/8}}\, dx\right] = \O\big(\e^{3-\frac{r}{4}}\big) =\, \frac{J_0^\e(y,z)}{K_\e(y,z)}\, \O\big(\e^{\frac{3-r}{4}}\big)  \vspace{0mm} $$ 
since \  $ J_0^\e(y,z) \asymp K_\e(y,z)\, \e^{9/4}\,$ by Proposition \ref{pro.partprincz}.  $\;\diamond$  \parm 
 
   We can now conclude this last section, by the following exact equivalent of the oscillatory integral $\,\II_\e(y,z)$ arising in the case  $\,z\not= 0\,$. \vspace{-1mm} 
\bpro  \label{pro.ControlIntPhiz} \    For any $(y,z)\in\R\times \R^*$,  with \   $\,{\ds C_\e(y,z) := \Big[ \frac{\pi\,z^2}{2\left(\pi -2\,\th\!{\frac{\pi}{2}}\right)\e^3} + \frac{y-\e}{\e^2} \Big]\,}$ and $\;\sigma'$ as in Lemma \ref{lem.cstte>0}, \ as $\,\e \sea 0$ \  we have 
$$ \II_\e(y,z) =\, \exp\!\left[-\pi^2\,C_\e(y,z)\right]\times C_\e(y,z)^{-3/4} \times \frac{(2\pi)^{11/4}\, \sigma'}{ \sqrt{\pi -2\,\th\!{\frac{\pi}{2}}}\,(\sh\pi)^{1/4}} \times\left(1+\O\big(\e^{1/3}\big)\right) . \vspace{-2mm} $$ 
\epro 
\ub{Proof} \quad  The first claim follows directly from (\ref{f.J01infty}) with $\,r = \frac{5}{3}\,$, Proposition \ref{pro.partprincz} and Lemma \ref{lem.ControlIntPhifz}. With Proposition \ref{pro.ChangContz=0>yb'z}, this gives the second one. $\;\diamond$ \parm 

   Propositions \ref{pro.equivIntPe0} and \ref{pro.ControlIntPhiz} together give the following wanted small time equivalent, which is the content of Theorem \ref{th.mainres}$(iii)$ when $\,z\not=0$, and actually extends Corollary \ref{cor.w=0=z<y}. 
\bcor \label{cor.w=0z} \  For any $(y,z)\in\R\times \R^*$, as $\,\e \sea 0$ we have 
$$ p_\e\big(0\, ; (0,y, z)\big) \sim\, \frac{\exp\!\left[-\pi^2\,C_\e(y,z)\right]}{\e^4\, C_\e(y,z)^{3/4}} \times \frac{(2\pi/\sh\pi)^{1/4}}{ \sqrt{\pi -2\,\th\!{\frac{\pi}{2}}}} \int_0^\infty\, \frac{\sin\!\left(\frac{3\pi}{16} + x\right)}{x^{1/4}}\, dx\, . \vspace{-2mm} $$ 
\ecor    
\vspace{-3mm}

\if{
$$ $$ 
\centerline{\bf * \'equation aux points-selle * \  (\`a partir de la Proposition \ref{pro.exprAvChgCont})}\parsn 
$$ {\frac{\e - y}{\e^2} }\rt1 - \frac{z^2}{2\e^3}\, \frac{\rt1}{f(0,x)} \left[1 - \frac{x\, f'(0,x)}{f(0,x)} \right] =\, - \frac{\Phi'(x)}{\Phi(x)}\,  $$ 
$$ =\, -\rt1 \f'(0,x) -\frac{3}{4x}Ê+ \frac{1}{4\sqrt{2x}}\times \frac{\sh\!\sqrt{2x} +\sin\!\sqrt{2x}}{\ch\!\sqrt{2x}  -\cos\!\sqrt{2x}} + \frac{f'(0,x)}{2\,f(0,x)}\, \raise1.9pt\hbox{}  $$
$$ =\, \frac{\rt1}{2\sqrt{2x}} \times {\frac{\sh\!\sqrt{2x} - \sin\!\sqrt{2x}}{\ch\!\sqrt{2x}  -\cos\!\sqrt{2x}}} -\frac{3}{4x}+ \frac{1}{4\sqrt{2x}}\times \frac{\sh\!\sqrt{2x} +\sin\!\sqrt{2x}}{\ch\!\sqrt{2x}  -\cos\!\sqrt{2x}} + \frac{f'(0,x)}{2\,f(0,x)}\, \raise1.9pt\hbox{.}  $$
Now 
$$  {f'(0,x)}\, = \, -\, \frac{d}{dx}\left[\frac{\sh\!\sqrt{x/2} + \sin\!\sqrt{x/2} - \rt1\!\big(\sh\!\sqrt{x/2} - \sin\!\sqrt{x/2}\,\big)}{(\ch\!\sqrt{x/2} + \cos\!\sqrt{x/2}\,)\sqrt{x/2}}\right]   \vspace{-0mm} $$ 
$$ = \, -\, \frac{\big[\ch\!\sqrt{x/2} + \cos\!\sqrt{x/2} - \rt1\!\big(\ch\!\sqrt{x/2} - \cos\!\sqrt{x/2}\,\big)\big] (\ch\!\sqrt{x/2} + \cos\!\sqrt{x/2}\,)\sqrt{x/2} }{(\ch\!\sqrt{x/2} + \cos\!\sqrt{x/2}\,)^2\, x\sqrt{2x}}\,  \vspace{-0mm} $$ 
$$ +\, {\ts \frac{\big[\sh\!\sqrt{x/2} + \sin\!\sqrt{x/2} - \rt1\!\big(\sh\!\sqrt{x/2} - \sin\!\sqrt{x/2}\,\big)\big] \big[\ch\!\sqrt{x/2} + \cos\!\sqrt{x/2} + (\sh\!\sqrt{x/2} - \sin\!\sqrt{x/2}\,)\sqrt{x/2} \,\big]}{(\ch\!\sqrt{x/2} + \cos\!\sqrt{x/2}\,)^2\, x\sqrt{2x}} }\,  \vspace{1mm} $$ 
{\small 
$$ = \, \frac{\sqrt{x/2}\, \big[ \sh\!^2\sqrt{x/2} - \sin^2\!\sqrt{x/2} - (\ch\!\sqrt{x/2} + \cos\!\sqrt{x/2}\,)^2\big]}{(\ch\!\sqrt{x/2} + \cos\!\sqrt{x/2}\,)^2\, x\sqrt{2x}} + {\frac{\big(\sh\!\sqrt{x/2} + \sin\!\sqrt{x/2}\, \big) }{(\ch\!\sqrt{x/2} + \cos\!\sqrt{x/2}\,)\, x\sqrt{2x}} }\,  \vspace{-0mm} $$ 
$$  +\, \frac{\rt1\!\big[\ch\!^2\sqrt{x/2} - \cos\!^2\sqrt{x/2} - \big(\sh\!\sqrt{x/2} - \sin\!\sqrt{x/2}\,\big)^2\,\big]\sqrt{x/2} }{(\ch\!\sqrt{x/2} + \cos\!\sqrt{x/2}\,)^2\, x\sqrt{2x}}\,  \vspace{-1mm} - { \frac{\rt1\!\big(\sh\!\sqrt{x/2} - \sin\!\sqrt{x/2}\,\big) }{(\ch\!\sqrt{x/2} + \cos\!\sqrt{x/2}\,)\, x\sqrt{2x}} }\,  \vspace{1mm} $$ } 
$$ = \, -\,\frac{1 + \ch\!\sqrt{x/2}\,\cos\!\sqrt{x/2}}{(\ch\!\sqrt{x/2} + \cos\!\sqrt{x/2}\,)^2\, x} + {\frac{\big(\sh\!\sqrt{x/2} + \sin\!\sqrt{x/2}\, \big)}{(\ch\!\sqrt{x/2} + \cos\!\sqrt{x/2}\,)\, x\sqrt{2x}} }\,  \vspace{-1mm} $$ 
$$  +\, \frac{\rt1 \sh\!\sqrt{x/2}\,\sin\!\sqrt{x/2} }{(\ch\!\sqrt{x/2} + \cos\!\sqrt{x/2}\,)^2\, x} - { \frac{\rt1\!\big(\sh\!\sqrt{x/2} - \sin\!\sqrt{x/2}\,\big)}{(\ch\!\sqrt{x/2} + \cos\!\sqrt{x/2}\,)\, x\sqrt{2x}} }\,  \vspace{1mm} $$ 
whence 
$$ \frac{f'(0,x)}{f(0,x)}\, =\, {\frac{\big(\sh\!\sqrt{x/2} + \sin\!\sqrt{x/2}\, \big) - \rt1\big(\sh\!\sqrt{x/2} - \sin\!\sqrt{x/2}\,\big) }{2x\,D(x)}} $$ 
$$ \hskip 20mm +\, \frac{ \rt1 \sh\!\sqrt{x/2}\,\sin\!\sqrt{x/2}  - \big(1+ \ch\!\sqrt{x/2}\,\cos\!\sqrt{x/2}\,\big)}{\sqrt{2x}\,(\ch\!\sqrt{x/2} + \cos\!\sqrt{x/2}\,)\, D(x)}\,  \vspace{1mm} $$ 
with 
$$ D(x) :=  \sqrt{x/2}\,\big(\ch\!\sqrt{x/2} + \cos\!\sqrt{x/2}\,\big) -\sh\!\sqrt{x/2} - \sin\!\sqrt{x/2} + \rt1\!\big(\sh\!\sqrt{x/2} - \sin\!\sqrt{x/2}\,\big) \raise1.9pt\hbox{.} $$ 

Now 
 \pars
\centerline{$\ch x + \cos x = 0\, \LRa \,\cos\!\big[(1\pm\rt1\!)\, x/2\big] = 0\, \LRa\, x= (1\pm\rt1)\big(\frac{\pi}{2}+k\pi\big)$,} \parn 
whence
$$ \ch\!\sqrt{x/2} + \cos\!\sqrt{x/2}\, = 0\,\LRa\, \sqrt{x/2}\, = (1\pm\rt1)\big({\ts\frac{\pi}{2}}+k\pi\big)\LRa\, x\in \pm \rt1\! \pi^2 (1+2\N)^2 . $$ 
}\fi 
\if{ 
\parb  \centerline{\bf Bribes semblant p\'erim\'ees et inutiles ; sauf peut-être pour $\,z\not= 0\,$ ??} \pars

Note that (\ref{df.fnf}) indeed defines an even function $f(x^2)$, which is holomorphic at any complex $\,x\notin (1\pm\rt1)\, \pi (\Z+\5)$, \  since \pars
\centerline{$\ch x + \cos x = 0\, \LRa \,\cos\!\big[(1\pm\rt1\!)\, x/2\big] = 0\, \LRa\,  e^{(\rt1\pm 1)\,x} = e^{\rt1\pi}$.} \parsn 
Similarly, using that $\,\sh\!^2 x + \sin^2\! x = \5 \sin\!\big[(1+\rt1\!)\, x\big] \sin\!\big[(1-\rt1\!)\, x\big]\,$ we have \pars 
\centerline{$\sh\!^2 x + \sin^2\! x = 0\, \LRa \,\sin\!\big[(1\pm\rt1\!)\, x\big] = 0\, \LRa\,  e^{2(\rt1\pm 1)\,x} = 1\, \LRa \,x\in \,(1\pm\rt1)\, \frac{\pi}{2}\,\Z\,$.} \pars 

   Recall from Corollary \ref{cor.calculQuadrF} that (\ref{df.fnf}) can be re-written as\,: \  ${\ds f\big(-2\rt1\! x^2\big) = 1-\frac{\th x}{x}\,}$\raise1.5pt\hbox{.} \parsn 
Now the equation  \  $r\,e^{\rt1 \theta} =\th(r\,e^{\rt1 \theta})\,$ has no solution for $\,r>0$ and $\,| \theta| < {\pi}/{2}\,$\raise1.5pt\hbox{.} \parn   
Indeed, this equation is equivalent to both \  $r= \cos\theta\, \th\!(r \cos\theta) + \sin\theta\, \tg(r \sin\theta)$ and  
\  $r\, \th\!(r \cos\theta)\,\tg(r \sin\theta) = \cos\theta\, \tg(r \sin\theta) - \sin\theta\,  \th\!(r \cos\theta) $, which by eliminating $r$ implies \  $\tg\theta = \frac{\sin(2r \sin\theta)}{\sh(2r \cos\theta)}< \tg\theta\;$ for $\,0<\theta<\pi/2\,$, a contradiction. By symmetry, the claim holds for negative $\theta$ as well, and also clearly for $\theta=0$. 
\parsn 
Therefore \  ${\ds \frac{\rt1 (-2\rt1\! x^2)}{f(-2\rt1\! x^2)} = \frac{2\,x^3}{x-\th x}\;}$ is holomorphic in $\{\,|\arg(x)|<\pi/2\}$. \  On the contrary, this does not hold for $\,\theta=\pi/2\,$, since  \  ${\ds \frac{\rt1 (2\rt1\! x^2)}{f(2\rt1\! x^2)} = \frac{2\,x^3}{\tg x-x}\;}$ has poles. \parsn 
Then near 0 we have\,: 
\begin{equation}   \label{f.DL'}  
\sqrt{\frac{\rt1\!\pi\, x^2}{f(x^2)}}\, =\, \sqrt{6\,\pi } \,\Big(1+{\ts\frac{\rt1 x^2}{10}}+\O(x^4)\Big) .  
\end{equation}
\blem \label{lem.minorRe} \  We have \  $\,{\ds \Re\!\left[\frac{\rt1 x^2}{f(x^2)} \right] \ge \, 6\, }$ \   for any $\,x\ge 0\,$. \vspace{-1mm} 
\elem 
\ub{Proof} \quad  By (\ref{f.partReelleb}) the claim is equivalent to \  $f(x) \ge 0\,$, \  where 
$$ f(x) :=\,\frac{\sh x-\sin x}{\ch x + \cos x}\times \frac{x^3}{6} - \left(x- \frac{\sh x+\sin x}{\ch x + \cos x}\right)^2 - 
\left(\frac{\sh x-\sin x}{\ch x + \cos x}\right)^2 $$ 
$$ =\,\frac{\sh x-\sin x}{\ch x + \cos x}\times \frac{x^3}{6} - x^2 + 2 \left(\frac{\sh x+\sin x}{\ch x + \cos x}\right) x - 2 
\left(\frac{\sh^2 x + \sin^2 x}{(\ch x + \cos x)^2}\right) = \,\frac{g(x)}{\ch x + \cos x}\, \raise1.9pt\hbox{,} $$ 
with \vspace{-2mm} 
$$ g(x) :=\, (\sh x-\sin x)\, \frac{x^3}{6} -(\ch x+\cos x)\,x^2 + 2 (\sh x+\sin x)\, x - 2 \left(\frac{\sh^2 x + \sin^2 x}{\ch x + \cos x}\right) . $$ 
Then we compute\,: \vspace{-2mm} 
$$ g'(x) =\, (\ch x-\cos x)\, \frac{x^3}{6} -(\sh x-\sin x)\, \frac{x^2}{2}\, =\, (\ch x-\cos x)\, \frac{x^3}{6}\, h(x)\,, \vspace{-2mm}  $$ 
with \vspace{-1mm} 
$$ h(x) :=\, x - 3 \left(\frac{\sh x - \sin x}{\ch x - \cos x}\right) , $$ 
and then\,: \quad $ h'(x) = (\ch x - \cos x)\2\, \ell(x)\,$, \quad with successively\,:   
$$ \ell(x) := \ch^2 x + \cos^2 x + 4\,\ch x\,\cos x -6\, , \  \ell'(x) = \sh\!(2x)- \sin(2x) + 4\,\sh x\,\cos x - 4\,\ch x\,\sin x\, ,  $$
$$ \ell''(x) = \, 4\,(\sh x-\sin x)^2\, >\, 0\,  \  \hbox{ for any $\,x> 0\,$},  \quad 0= \ell'(0) = \ell(0)= h(0)= g(0)= f(0)\, . $$
Hence for any $\,x> 0\,$ we successively have\,: \  
$$ \ell'(x)> 0\,, \  \ell(x)> 0\,, \ h'(x)> 0\,, \  h(x)> 0\,, \  g'(x)> 0\,, \   g(x)> 0\,, \ f(x)> 0\,. \;\;\diamond $$
}\fi 



\vspace {-3mm}  

\section{References} \label{Ref} 
\par \vskip 4mm 


\vbox{ \noindent 
{\bf [A]} \   Azencott R.  \quad {\it Densit\'e des diffusions en temps petit : d\'eveloppements asymptotiques.} \par \hskip 30mm   S\'em. Proba. XVIII, Lecture Notes n$^o\,$1059, 402-498, Springer 1984. }
\par \bigskip 

\vbox{ \noindent 
{\bf [BA1]} \  Ben Arous G.  \quad {\it D\'eveloppement asymptotique du noyau de la chaleur hypoelliptique \par \hskip 32mm   hors du cut-locus.}  \  Ann. sci. \'E.N.S., S\'er. 4, t. 21 n$^o\,$3, 307-331, 1988. }
\par \bigskip 

\vbox{ \noindent 
{\bf [BA2]} \  Ben Arous G.  \quad {\it D\'eveloppement asymptotique du noyau de la chaleur hypoelliptique \par \hskip 38mm   sur la diagonale.}  \  Ann. Inst. Fourier 39, n$^o\,$1, 73-99, 1989. }
\par \bigskip 

\vbox{ \noindent 
{\bf [BA3]} \  Ben Arous G.  \quad {\it M\'ethode de Laplace et  de la phase stationnaire \par \hskip 38mm   sur l'espace de Wiener.}  \  Stochastics vol. 25, 125-153, 1988. }
\par \bigskip 

\vbox{ \noindent 
{\bf [BA-G]} \  Ben Arous G., Gradinaru M.  \quad  {\it Singularities of Hypoelliptic Green Functions.}  \par\hskip 70mm   Potential Analysis 8, 217-258, 1998. }
\par \bigskip 


\vbox{ \noindent 
{\bf [B-O]} \  Beghin L., Orsingher E.   \ {\it  On the maximum of the generalized Brownian bridge.}  \par \hskip 54mm  Lithuanian Math. J. vol. 39, n$^o\,$2, 157-167, 1999. }
\par \bigskip 

\vbox{ \noindent 
{\bf [B-P-Y]} \  Biane P., Pitman J., Yor M. \  {\it  Probability laws related to the Jacobi theta }  \par 
\hskip 12mm   {\it  and Riemann zeta functions, and Brownian excursions.}  \par \hskip 54mm  Bull. Amer. Math. Soc. vol. 38, n$^o\,$4, 435-465, 2001. }
\par \bigskip 

\vbox{ \noindent 
{\bf [B]} \  Bismut J.M.  \quad {\it Large deviations and the Malliavin calculus.} \par \hskip 31mm  Progress in math. n$^o\,$45, Birkh\" auser 1984. }
\par \bigskip 


\vbox{ \noindent 
{\bf [Ca]} \  Castell F. \ \  {\it Asymptotic expansion of stochastic flows.}
\par \smallskip \hskip 28mm  Prob. Th. Rel. Fields 96, 225-239, 1993. }
\par \bigskip 

\vbox{ \noindent 
{\bf [CM-E]} \  Chaleyat-Maurel M., Elie L.  \quad  {\it  Diffusions gaussiennes.}   \par \hskip 44mm  
   in Ast\'erisque 84-85, \og g\'eod\'esiques et diffusions en temps petit\fg, \par \hskip 44mm  
   s\'eminaire de probabilit\'es de Paris VII, 255-279, 1981. }
\par \bigskip 

\vbox{ \noindent 
{\bf [D-M]} \  Delarue F., Menozzi S. \ \  {\it  Density estimates for a random noise propagating through\par \hskip 33mm  a chain  of differential equations. }  \hskip 3mm   J. F. A.   $n^o\,$259, 1577-1630, 2010. }
\par \bigskip 

\vbox{ \noindent 
{\bf [Du]} \  Dudley R.M. \ \  {\it Lorentz-invariant Markov processes in relativistic phase space.}
\par \smallskip \hskip 39mm  Arkiv f\"or Matematik 6, n$^o\,$14, 241-268, 1965. }
\par \bigskip 



\vbox{ \noindent 
{\bf [Fa]} \  Fang Shizan  \  {\it Introduction to Malliavin Calculus.} \par
\hskip 29mm Math. Series for Grad. Students, Springer 2003.  }
\par \bigskip 

\vbox{ \noindent 
{\bf [F]} \  Franchi J.  \  {\it Small time asymptotics for an example of strictly hypoelliptic heat kernel.} \par
\hskip 19mm S\'eminaire de Probabilit\'es, volume XLVI, 71-103, LNM 2123, Springer 2014.  }
\par \bigskip 

\vbox{ \noindent 
{\bf [F-LJ1]} \  Franchi J., Le Jan Y.  \ \  {\it Relativistic Diffusions and Schwarzschild Geometry.} \par
\hskip 50mm Comm. Pure Appl. Math., vol. LX, n$^o\,$2, 187-251, 2007.}
\par \bigskip 


\vbox{ \noindent 
{\bf [F-LJ2]} \   Franchi J., Le Jan Y.  \ \   {\it Hyperbolic Dynamics and Brownian Motion}. \par \hskip 12mm   Oxford Mathematical Monographs,  Oxford Science Publications, august 2012.  } 
\par \bigskip 

\vbox{ \noindent 
{\bf [H]} \  Hairer M.  \  {\it On Malliavin's proof of H\"ormander's Theorem.} \par
\hskip 30mm Bull. Sci. math. 135, 650-666, 2011.}
\par \bigskip 


\vbox{  \noindent 
{\bf [I-M]} \  Ikeda N.,  Manabe S. \  {\it  Asymptotic formulae for stochastic oscillatory integrals.}  \par \hskip 2mm  Asymp. prob. prob. th., Kyoto 1990 (Ed. Elworthy KD \& Ikeda N), Pitman Res. Notes Math. S. n$^o\,$284, Longman Sc;Tech. 1993, 1981. }
\par \bigskip  

\vbox{  \noindent 
{\bf [I-W]} \  Ikeda N.,  Watanabe S. \  {\it  Stochastic differential equations and diffusion processes.}  \par \hskip 52mm  North-Holland Kodansha, 1981. }
\par \bigskip  



%
\vbox{ \noindent 
{\bf [L]} \   L\'eandre R. \quad {\it  Int\'egration dans la fibre associ\'ee \`a une diffusion d\'eg\'en\'er\'ee.} \par \hskip 33mm    Prob. Th. Rel. Fields  76,  n$^o\,$3, 341-358, 1987.}
\par \bigskip 

\vbox{ \noindent 
{\bf [M]} \  Malliavin P.  \quad   {\it  Stochastic Analysis.}  \par \hskip 19mm  Grundlehren der mathematischen Wissenschaften vol. 313, Springer 1997. }
\par \bigskip 

\vbox{ \noindent 
{\bf [M-T]} \  Malliavin P.,  Taniguchi S.  \  {\it  Analytic Functions, Cauchy Formula, and \par \hskip 8mm  Stationary Phase on a Real Abstract Wiener Space.}  \hskip 3mm  J.F.A. n$^o\,$143, 470-528, 1997. }
\par \bigskip 



\vbox{ \noindent 
{\bf [R-Y]} \  Revuz D., Yor M. \quad {\it Continuous Martingales and Brownian Motion.} \par 
\hskip 8mm   Grundlehren der mathematischen Wissenschaften vol. 293,  Springer, Berlin 1999. }
\par \bigskip 

%
\vbox{ \noindent 
{\bf [T1]} \  Taniguchi S.  \quad  {\it  On the Exponential Decay of Oscillatory Integrals  \par \hskip 24mm  on an
Abstract Wiener Space.}  \hskip 3mm  J.F.A. n$^o\,$154, 424-443, 1998. }
\par \bigskip 

\vbox{ \noindent 
{\bf [T2]} \  Taniguchi S.  \quad  {\it  Exponential decay of stochastic oscillatory integrals  \par \hskip 20mm  on classical Wiener spaces.}  \hskip 3mm  J. Math. Soc. Japan vol. 55 n$^o\,$1, 59-79, 2003. }
\par \bigskip 


\vbox{ \noindent 
{\bf [V]} \  Varadhan S.R.S. \ \  {\it  Diffusion Processes in a Small Time Interval.}
\par \hskip 34mm  Comm. Pure Applied Math.   n$^o\,$20, 659-685, 1967. }
\par \bigskip 

\vbox{ \noindent 
{\bf [Y]} \  Yor M. \quad {\it Some aspects of Brownian motion. Part I. Some special functionals.}
\par  \hskip 21mm Lectures in Mathematics, ETH Z\"urich. Birkh\"auser Verlag, Basel, 1992. }
\par  \bigskip 
\centerline{------------------------------------------------------------------------------------------------------------}\par 
\parmn 
{\bf Keywords}\,: \  Strictly hypoelliptic diffusion, heat kernel, small time asymptotics, Brownian bridge, Malliavin calculus, oscillatory integral. \parsn 
{\bf AMS classification}\,:  \  primary \  58J65 ; \  secondary \  35K65, 60H07, 35K05, 60J65. 
\par  \parmn 
\centerline{------------------------------------------------------------------------------------------------------------}\par

\if{ 
$$ $$ 
\ub{\bf Questions restant en suspens} :  \parsn
$\bullet$ \  Comprendre si possible le lien entre $w=0$ et cut-locus. \parn
$\bullet$ \  Passer \`a la diffusion de Dudley ??  \parn  
$\bullet$ \  \'Elargir \`a un cas un peu g\'en\'erique... Avec $V_0$ born\'e ? ou incluant aussi Dudley (avec $F$ born\'ee, tronqu\'ee et liss\'ee hors d'un voisinage compact, comme dans [BA1] ?) ?? \parn 
$\bullet$ \  \'Elargir \`a un cas g\'en\'erique... avec $V_1$ non trivial, voire $V_1,\ldots , V_k\,$ au lieu du seul $V_1$  ? \parn 
$\bullet$ \ Dans quels cas plus sophistiqu\'es ou compliqu\'es la n\'ecessit\'e de $F$ (de [BA1]) r\'eappara\^\i tra-t-elle ??  \parn
}\fi 

\end{document}